%% file: shortened.tex
\documentclass{article}
\usepackage{graphicx, hyperref, mathtools, amsfonts, amssymb, blkarray, amsthm} % Required for inserting images

\usepackage{tikz-cd}                    % For commutative diagrams
\usetikzlibrary{arrows.meta}            % For hook arrows in commutative diagrams

\definecolor{darkgreen}{rgb}{0.0, 0.35, 0.0}

\usepackage[margin=2.5cm]{geometry}
\usepackage{float}

\usepackage{nicematrix}

\usepackage[ruled]{algorithm2e}

\usepackage{subcaption}

\usepackage{biblatex}
\makeatletter
\newcommand*\rel@kern[1]{\kern#1\dimexpr\macc@kerna}
\newcommand*\widebar[1]{%
  \begingroup
  \def\mathaccent##1##2{%
    \rel@kern{0.8}%
    \overline{\rel@kern{-0.8}\macc@nucleus\rel@kern{0.2}}%
    \rel@kern{-0.2}%
  }%
  \macc@depth\@ne
  \let\math@bgroup\@empty \let\math@egroup\macc@set@skewchar
  \mathsurround\z@ \frozen@everymath{\mathgroup\macc@group\relax}%
  \macc@set@skewchar\relax
  \let\mathaccentV\macc@nested@a
  \macc@nested@a\relax111{#1}%
  \endgroup
}
\makeatother

\usepackage{tikz}
\usetikzlibrary{decorations.pathreplacing}
\makeatletter
\def\underbrace#1{\@ifnextchar_{\tikz@@underbrace{#1}}{\tikz@@underbrace{#1}_{}}}
\def\tikz@@underbrace#1_#2{\tikz[baseline=(a.east)] {\node (a) {\(#1\)}; \draw[ultra thick,line cap=round,decorate,decoration={brace,amplitude=5pt}] (a.south east) -- node[below,inner sep=7pt] {\(\scriptstyle #2\)} (a.south west);}}
\makeatother
\date{\today}

\newtheorem{theorem}{Theorem}[section]

\newtheorem{lemma}[theorem]{Lemma}
\theoremstyle{definition}
\newtheorem{definition}{Definition}[section]
\newtheorem{example}{Example}[section]
\newtheorem{remark}{Remark}

\newcommand{\RR}[0]{\mathbb{R}}

\newcommand{\im}[0]{\text{Im}\,}

\DeclareMathOperator*{\rank}{\text{rank}}

\begin{document}
% \maketitle
\title{\fontsize{45}{8}\selectfont New persistent chain spaces: identifying and computing the relevant persistent Laplacians of a filtration}
\author{Joris Kirchner, Martina Vittorietti}
\maketitle

\begin{abstract}Persistent homology (PH) summarizes the topological features of data that persist across scales, but it is insensitive to many geometric properties. Persistent Laplacians (PLs) address this limitation: for a pair of simplicial complexes in a filtration, the kernel of the persistent Laplacian recovers the persistent Betti number, while its non-zero spectrum carries additional geometric information. However, PLs are defined for individual pairs of filtration values, so analysing a whole filtration requires evaluating many operators that largely contain the same information, and the result is hard to summarize and interpret.
    In this paper, we show that this redundancy has a precise structure. For each pair of filtration values, we decompose the space of chains underlying the persistent up-Laplacian into the contributions inherited from the neighboring pairs and a remaining part, which we call the new persistent chain space. We show that the new persistent chain spaces of all earlier pairs form a direct sum decomposition of this space, so that every persistent up-Laplacian in the filtration is determined by them. Moreover, the new persistent chain spaces are non-trivial only at pairs given by the persistence pairing in the same and in the next homological dimension, so that the persistence pairing determines exactly where the persistent up-Laplacian can change. Experiments on several types of point clouds show that evaluating persistent Laplacians only at these pairs reduces the computation time, both for our implementation and for existing software.%In this paper, we decompose the space of persistent chains underlying the persistent up-Laplacian into new persistent chain spaces, each of which records the information that first appears at a given pair of filtration values. We prove that the persistent up-Laplacian at every pair of filtration values can be reconstructed from these spaces, so that it suffices to evaluate it at a small set of locations determined by the persistence pairing; all other persistent up-Laplacians are recovered without loss of information. The new spaces can be computed with the standard persistent homology reduction algorithm, which leads to an algorithm that avoids redundant computations. Finally, the construction yields an extension of the classical persistence barcode, the persistent geometric barcode, which retains the structure of the classical barcode while including the geometric information captured by persistent Laplacians.
\end{abstract}

% barcodes with additional intervals and an added weight per interval. The new information will be obtained via the persistent Laplacian, a recently proposed operator that encodes the topological information of persistent homology in its kernel and geometric information in its non-zero spectra. Usage of the operator is difficult as it functions differently from standard Topological Data Analysis methods. By proposing the multiplicity operator, a new method is introduced to summarize information captured by the persistent Laplacian into a form that is more commonly used. This allows methods based on persistent homology to additionally describe geometric properties, opposed to only topological features. For the multiplicity operator the trace is analyzed and the features captured by it are discussed. In certain situations, individual eigenvalues are still deemed hard to explain. Therefore, an adjusted multiplicity operator is proposed that contains separately interpretable eigenvalues. Finally, the operators are used to classify handwritten digits from the MNIST dataset and to classify different generation processes of artificially made cross sections of crystalline structures. For the analysis of the cross sections no improvement was found over standard persistent homology. However, for the MNIST dataset, performance in terms of accuracies exceeded persistent homology and was similar to the performance of persistent Laplacians, indicating that most of its information is captured by the new operator.

\section{Introduction}
\label{ch:introduction}
\input{shortened_sections/introduction_v3}

\section{Mathematical Background}
\label{ch:mathematical_background}
\input{shortened_sections/mathematical_background_v3}

\section{Dissecting the Persistent Laplacian}
\label{ch:new_chain_spaces}
\input{shortened_sections/newchainspaces_v4}

\section{Algorithm}
\label{ch:algorithm}
\input{shortened_sections/algorithm_v2}

% \section{Applications}
% \label{ch:applications}
% \input{shortened_sections/applications}

\section{Comparison with other methods}
\label{ch:comparison}
\input{shortened_sections/comparison_v2}

\section{Conclusion and Discussion}
\label{ch:conclusion}
\input{shortened_sections/conclusion_v3}

\newpage
\printbibliography

\newpage
\appendix
% \section{Proofs}
\input{shortened_sections/supplementary_v3}

\end{document}

%% file: shortened_sections/introduction_v3.tex
% Revised introduction.
% TODO  = information that seems to be missing.
% CHECK = places where I had to guess the intended meaning.
% NEW   = sentences I added; remove them if you prefer.

Over the past two decades, topological data analysis (TDA) has developed into a major framework for studying the shape of data, combining ideas from algebraic topology, geometry, statistics, and computation~\cite{carlsson_topology_data,ghrist_barcodes,chazal_tda}. 
TDA aims to extract global and multiscale structure from data, often under the assumption that the observations are sampled from an underlying geometric object such as a manifold. Its central tool is persistent homology (PH), which studies a filtration of simplicial complexes built from the data and records how the homology groups evolve across scales~\cite{edelsbrunner_persistent_homology,computing_ph}. The resulting persistent Betti numbers and persistence barcodes or, equivalently, persistence diagrams, record at which scales connected components, loops, voids, and higher-dimensional holes appear and disappear. These summaries are stable under perturbations of the data~\cite{cohen_steiner_stability} and can be vectorized for statistical learning~\cite{persistence_landscapes,persistence_images}. As a result, PH has been successfully applied in many areas, including shape analysis \cite{persistent_homology_transform}, biology \cite{nicolau2011topology}, chemistry \cite{cang2017topologynet}, materials science \cite{hiraoka2016hierarchical,pv_tests}, neuroscience \cite{bendich2016persistent}, and image analysis \cite{carlsson2008local} (see also~\cite{tda} for an extensive review).

Despite these successes, classical persistent homology is intrinsically topological. It records the presence and persistence of homology classes, but it does not directly encode many geometric or metric features of the underlying data. For example, two filtrations may have the same persistent Betti numbers and persistence barcode, while differing substantially in the size, location, arrangement, or interaction of their cycles. This limitation has motivated several extensions of PH that retain its multiscale topological robustness while incorporating additional geometric information. One important example is the persistent homology transform (PHT), which represents an embedded shape by the collection of persistence diagrams of its directional height filtrations, and which is injective on certain classes of shapes~\cite{persistent_homology_transform}. Another recent direction is based on spectral invariants, in particular persistent Laplacians \cite{laplacians_spectral_graph,laplacians_survery}.

Persistent Laplacians (PLs) were introduced as a spectral extension of persistent homology~\cite{laplacians_spectral_graph}. They are rooted in the graph Laplacian and its higher-dimensional analogue, the combinatorial or Hodge Laplacian on simplicial complexes~\cite{combinatorial_laplacians,harmonische_funktionen}.
For a pair of complexes $K_s \subseteq K_t$ in a filtration, the persistent Laplacian is the sum of the down-Laplacian of the smaller complex $K_s$ and a persistent up-Laplacian, which is built from the chains of the larger complex $K_t$ whose boundary lies in $K_s$. The multiplicity of the zero eigenvalue of the persistent Laplacian equals the corresponding persistent Betti number, so the persistent Laplacian retains the topological information of PH, while its non-zero spectrum contains additional geometric and combinatorial information~\cite{laplacians_algorithms,botnan_cheeger_2026}. This makes PLs attractive for applications where topology alone is too coarse. Features derived from persistent Laplacians have already shown strong empirical performance in molecular and biological prediction tasks~\cite{pers_lap_covid,pers_lap_protein_ligand} and in image recognition problems such as MNIST classification~\cite{pers_lap_mnist}.

At the same time, persistent Laplacians introduce new challenges. In classical PH, a persistence barcode provides a compact global summary of an entire filtration, whereas PLs are typically computed for individual pairs of complexes. Thus, although a persistent Laplacian can compare geometric and topological properties between two filtration levels, there is not yet a standard barcode-like object that summarizes how the persistent Laplacian changes across the whole filtration. Moreover, PLs are expensive to compute. A persistent Laplacian is defined for every pair of filtration values, so their number grows quadratically with the length of the filtration. 
In addition, each individual PL can be costly for large complexes and in higher homological dimensions, since it requires the construction and spectral analysis of large matrices derived from the boundary maps. Even though progress has been made in this area~\cite{hermes,computationalefficiency}, the computational cost still limits the direct use of PLs on large filtrations, which motivates methods that avoid recomputing redundant information at every scale.

Recent work has begun to address the question of how to summarize and vectorize the information contained in persistent Laplacians. In particular, the authors of~\cite{persistent_laplacian_diagrams} introduced Persistent Laplacian Diagrams and Persistent Laplacian Images, which assign scalar spectral signatures to birth--death coordinates and then construct vectorized image representations analogous to persistence images. Their approach is designed primarily as a diagrammatic and vectorized feature representation for machine learning, in which a chosen signature, such as the spectral gap, the entropy, or an eigenvector-based statistic, is evaluated on the persistent Laplacians of selected filtration pairs.

In this paper, we take a different but complementary perspective. Rather than evaluating a chosen signature at selected pairs, we determine where in a filtration the persistent Laplacians change, and we use this to reduce the cost of computing the persistent Laplacians of a whole filtration. To this end, we decompose the space of persistent chains $C_{q+1}^{s,t}$, on which the persistent up-Laplacian $\Delta_{q,+}^{s,t}$ is built, into four parts: the contribution already present at the earlier end value $t-1$ (the vertical effect), the contribution already present at the earlier start value $s-1$ (the horizontal effect), the contribution already present at $(s-1,t-1)$ (diagonal), and a remainder.
This remainder, which we call the \emph{new persistent chain space} $N_{q+1}^{s,t}$, is a small subspace that captures exactly the information that first appears at $(s,t)$. We prove that the spaces $N_{q+1}^{s',t'}$ with $s' \leq s$ and $t' \leq t$ form a direct sum decomposition of $C_{q+1}^{s,t}$, so that every persistent up-Laplacian in the filtration can be reconstructed from them. The down part of the persistent Laplacian can be treated in the same way, since the non-zero eigenvalues of the down-Laplacian of $K_s$ coincide with those of the up-Laplacian of $K_s$ one dimension lower. %Since the positive spectrum of the down part of a persistent Laplacian is fully captured by the up part of the persistent Laplacian on the diagonal of one dimension lower, this determines all full persistent Laplacians of the filtration. 
Moreover, we show that $N_{q+1}^{s,t} \neq \{0\}$ only at pairs $(s,t)$ determined by the persistence pairing. 
% These pairs, together with the spaces $N_{q+1}^{s,t}$, form the barcode-like summary mentioned above. 
The new persistent chain spaces can be computed with the standard reduction algorithm for persistent homology, at a cost that depends on the number of simplices rather than on the number of pairs of filtration values. This avoids the redundant computations of approaches that evaluate persistent Laplacians independently for many pairs. In experiments on several types of point clouds, our method is faster than or on par with existing methods on filtrations with low-dimensional simplices, and the pairs it finds can also be used to speed up existing methods such as PETLS~\cite{computationalefficiency}.

The paper is organized as follows. Section~\ref{ch:mathematical_background} reviews the necessary background on persistent homology, barcodes, combinatorial Laplacians, and persistent Laplacians. Section~\ref{ch:new_chain_spaces} analyzes the redundant information contained in the persistent Laplacians of a filtration, introduces the horizontal and vertical operators and the new persistent chain spaces, and proves that these spaces allow the reconstruction of all persistent up-Laplacians. Section~\ref{ch:algorithm} shows how the new persistent chain spaces are obtained from the standard reduction algorithm and and how the persistent up-Laplacians are computed from them. In Section~\ref{ch:comparison}, we test the algorithm on several examples and compare its running time with that of existing methods. Conclusions and directions for future research are presented in Section~\ref{ch:conclusion}. All proofs are collected in the appendix. Code is available at \url{https://github.com/siroj99/NewPersistentChains}.

%% file: shortened_sections/mathematical_background_v3.tex
% Revised mathematical background.
% TODO  = information that seems to be missing.
% CHECK = places where I had to guess the intended meaning.
% NEW   = sentences or headings I added; remove them if you prefer.
% The commented-out older versions of paragraphs have been removed.

In this section, we introduce the notation and the mathematical background used throughout the paper.

\subsection{Homology theory}
Let $\mathcal{P}=\{p_0,\dots,p_m\}$ be a finite ordered set of points. A (finite, abstract) simplicial complex $K$ on $\mathcal{P}$ is a collection of non-empty subsets of $\mathcal{P}$, called simplices, that is closed under taking non-empty subsets; a simplex with $q+1$ elements is called a $q$-simplex. We work with ordered simplices: a $q$-simplex is written as the ordered tuple $[p_{i_0},\dots,p_{i_q}]$ with $0\leq i_0<\cdots<i_q\leq m$. We write $S_q^K$ for the set of $q$-simplices of $K$ and $n_q^K:=|S_q^K|$ for their number. The $q$-th chain group $C_q^K$ is the real vector space with basis $S_q^K$, equipped with the inner product for which $S_q^K$ is an orthonormal basis. The $q$-th boundary operator $\partial_q^K:C_q^K\to C_{q-1}^K$ is the linear map defined on a $q$-simplex $\sigma=[v_0,\dots,v_q]$ by
\begin{equation}
    \partial_q^K(\sigma):=\sum_{i=0}^q(-1)^i d_i\sigma,
\end{equation}
where $d_i\sigma:=[v_0,\dots,v_{i-1},v_{i+1},\dots,v_q]$ is the $(q-1)$-simplex obtained by omitting the vertex $v_i$ from $\sigma$. By convention, $\partial_0^K:=0$.
Since $\partial_{q}^K\partial_{q+1}^K=0$, the $q$-th homology group $H_q(K):=\ker\partial_q^K/\im\partial_{q+1}^K$ is well defined; its dimension $\beta_q^K:=\dim H_q(K)$ is the $q$-th Betti number of $K$.

\begin{figure*}[t!]
    \centering
    \begin{subfigure}[t]{0.7\textwidth}
        \centering
        \includegraphics[height=1.5in]{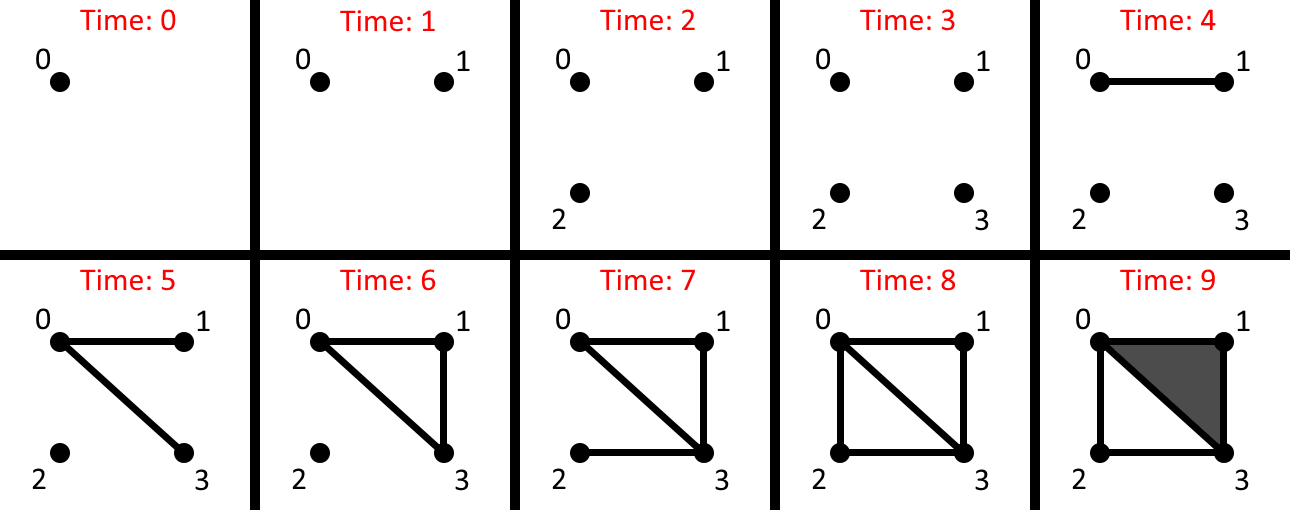}
        \caption{Filtration of simplicial complexes.}
    \end{subfigure}
    \begin{subfigure}[t]{0.24\textwidth}
        \centering
        \includegraphics[height=1.5in]{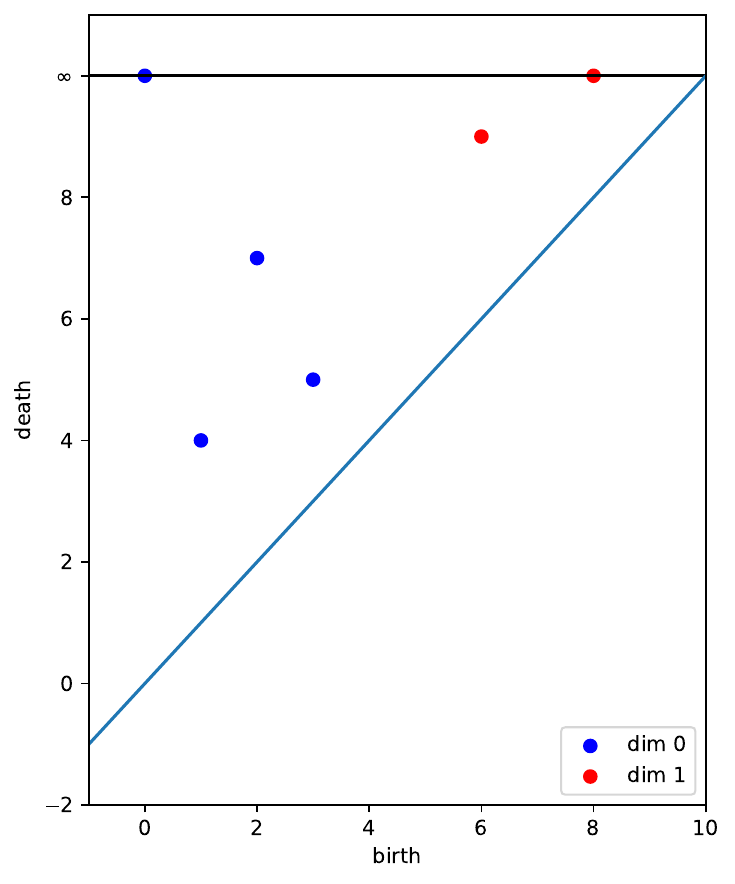}
        \caption{Persistence diagram.}
    \end{subfigure}
    \caption{(a) Example filtration of simplicial complexes. Points indicate 0-simplices, edges indicate 1-simplices, and the filled triangle indicates a 2-simplex. (b) Persistence diagram of the filtration in (a).\label{fig:filtration_example}}
\end{figure*}

Instead of analyzing the homology of a single simplicial complex, we analyze a filtration of simplicial complexes, that is, a nested sequence $K_0\subseteq K_1\subseteq\dots\subseteq K_T$; see Figure~\ref{fig:filtration_example}(a) for an example. Such filtrations are typically built from data, for example with the \v{C}ech, Vietoris--Rips or alpha constructions.
In practice, filtrations are often indexed by a continuous parameter, such as a radius. Since a filtration of a finite complex changes only at finitely many parameter values $r_0<r_1<\dots<r_T$, we may and do index it by the integers $0,1,\dots,T$, where $K_s$ denotes the complex at parameter value $r_s$. With this convention, $s-1$ denotes the filtration index preceding $s$.
For a filtration, we abbreviate $C_q^s:=C_q^{K_s}$, $\partial_q^s:=\partial_q^{K_s}$ and $n_q^s:=n_q^{K_s}$. We order the simplices compatibly with the filtration, so that the simplices of $K_{s-1}$ precede those of $K_s\setminus K_{s-1}$.

For a simplicial pair $K\subseteq L$, the persistent Betti number $\beta_q^{K,L}:=\rank\big(H_q(K)\to H_q(L)\big)$ is the rank of the map in homology induced by the inclusion $K\hookrightarrow L$; it counts the $q$-dimensional features that are present in $K$ and still present in $L$. For a filtration and $s\leq t$, we write $\beta_q^{s,t}:=\beta_q^{K_s,K_t}$.

The collection of all persistent Betti numbers is redundant: it can be summarized by recording only when features appear (are \emph{born}) and when they disappear (\emph{die}). This gives rise to the persistence barcode, a multiset of half-open intervals, called bars. For $0\leq s<t\leq T$, the multiplicity of the bar $[s,t)$ in the barcode is
\begin{equation}\label{eq:multiplicity_betti}
    \mu_q^{s,t}:=
    \left(\beta_q^{s,t-1}-\beta_q^{s,t}\right)
    -
    \left(\beta_q^{s-1,t-1}-\beta_q^{s-1,t}\right),
\end{equation}
and the multiplicity of the bar $[s,\infty)$, which corresponds to features that persist until the end of the filtration, is
\begin{equation}
    \mu_q^{s,\infty}:=\beta_q^{s,T}-\beta_q^{s-1,T}.
\end{equation}
Here we use the convention $\beta_q^{-1,t}:=0$ for all $t$.

Besides being drawn as a collection of bars, the barcode is often visualized as a persistence diagram: a scatter plot of the points $(s,t)$, with birth times on the horizontal axis and death times on the vertical axis; see Figure~\ref{fig:filtration_example}(b) for an example. The multiplicities are usually not shown in such a plot.

\subsection{Laplacians on simplicial complexes}
Before looking at simplicial complexes, we start with the standard Laplacian defined on graphs. For a simple undirected graph $G=(V,E)$ with vertex set $V=\{v_1,\dots,v_n\}$ and edge set $E$, consisting of two-element subsets of $V$, the graph Laplacian is defined as $L_G:=D-A$. Here $D$ is the $n\times n$ diagonal matrix whose $i$-th diagonal entry is the degree of $v_i$, i.e.\ the number of edges containing $v_i$, and $A$ is the adjacency matrix, the symmetric binary matrix with $A_{ij}=1$ if $\{v_i,v_j\}\in E$ and $A_{ij}=0$ otherwise. Viewing $G$ as a one-dimensional simplicial complex, one has $L_G=[\partial_1][\partial_1]^T$, where $[\partial_1]$ denotes the matrix of $\partial_1$ in the simplex bases. This definition extends to simplicial complexes of arbitrary dimension, where it is called the combinatorial or Hodge Laplacian.

\begin{definition}[Combinatorial Laplacian]\label{def:combinatorial_laplacian}
    For a simplicial complex $K$, the $q$-th combinatorial Laplacian $\Delta^{K}_q:C_q^K\to C_q^K$ is defined by
    \begin{equation}
        \Delta_q^{K}:=\underbrace{\partial_{q+1}^{K}(\partial_{q+1}^{K})^*}_{\Delta_{q,+}^{K}}+\underbrace{(\partial_q^K)^*\partial_q^K}_{\Delta_{q,-}^{K}},
    \end{equation}
    where $(\partial_{q+1}^K)^*:C_q^K\to C_{q+1}^K$ denotes the adjoint of $\partial_{q+1}^K$ with respect to the inner products on the chain groups. We call $\Delta^{K}_{q,+}$ and $\Delta_{q,-}^{K}$ the up-Laplacian and the down-Laplacian, respectively.
\end{definition}
With respect to the orthonormal simplex bases, the matrix of the adjoint is the transpose: $[(\partial_{q+1}^K)^*]=[\partial_{q+1}^K]^T$. Since $\partial_0^K=0$, we have $\Delta_0^K=\partial_1^K(\partial_1^K)^*$, whose matrix is the graph Laplacian of the one-dimensional skeleton of $K$.

Combinatorial Laplacians are relevant for TDA because of the following classical result, which shows that the Betti numbers are encoded in the kernel of the Laplacian.
\begin{theorem}[{Combinatorial Hodge theorem~\cite{harmonische_funktionen}}]\label{thm:combi_lap_betti}
    For each $q\geq 0$, $\beta_q^K=\dim\ker \Delta_q^K$.
\end{theorem}

The non-zero eigenvalues are often aggregated into a single value for each filtration index and each dimension~\cite{laplacians_spectral_graph,laplacians_survery}. The most common choice is the smallest non-zero eigenvalue, which has a proven connection to effective resistance in graph theory~\cite{effective_resistance}.
Some works instead consider the largest eigenvalue \cite{laplacians_spectral_graph}.
%To visualize these values over a filtration, they are often plotted together with the Betti numbers~\cite{laplacians_spectral_graph,laplacians_survery}; Figure~\ref{fig:combi_lap_min_eval_example_filtration} shows this for the smallest non-zero eigenvalue.
To visualize these values over a filtration, they are often plotted together with the Betti numbers as functions of the filtration index~\cite{laplacians_spectral_graph,laplacians_survery}. 
\subsubsection{Persistent Laplacian}
To apply the combinatorial Laplacian to a filtration rather than to a single complex, the persistent Laplacian was proposed~\cite{laplacians_spectral_graph}. It is defined for a simplicial pair $K\subseteq L$ and is applied to a filtration by considering all pairs $K_s\subseteq K_t$ with $s\leq t$. Consider the subspace
\begin{equation}
    C_q^{K,L}:=\{c\in C_q^L:\partial_q^L(c)\in C_{q-1}^K\}\subseteq C_q^L
\end{equation}
of $q$-chains of $L$ whose boundary lies in $C_{q-1}^K$, equipped with the inner product inherited from $C_q^L$, 
and let $\partial_q^{K,L}:C_q^{K,L}\to C_{q-1}^K$ denote the restriction of $\partial_q^L$ to $C_q^{K,L}$.

\begin{definition}[{Persistent Laplacian~\cite{laplacians_spectral_graph}}]\label{def:persistent_laplacian}
    For a simplicial pair $K\subseteq L$, the $q$-th persistent Laplacian $\Delta_q^{K,L}:C_q^K\to C_q^K$ is defined by
    \begin{equation}
        \Delta_q^{K,L}:=\underbrace{\partial_{q+1}^{K,L}(\partial_{q+1}^{K,L})^*}_{\Delta_{q,+}^{K,L}}+\underbrace{(\partial_q^K)^*\partial_q^K}_{\Delta_{q,-}^{K}}.
    \end{equation}
    We call $\Delta^{K,L}_{q,+}$ the persistent up-Laplacian; the second term $\Delta_{q,-}^{K}$ is the down-Laplacian of $K$ from Definition~\ref{def:combinatorial_laplacian}.
\end{definition}
The persistent Laplacian is self-adjoint and positive semidefinite, so its eigenvalues are real and non-negative. For a filtration, we write $C_q^{s,t}:=C_q^{K_s,K_t}$, $\partial_q^{s,t}:=\partial_q^{K_s,K_t}$, $\Delta_q^{s,t}:=\Delta_q^{K_s,K_t}$ and $\Delta_{q,+}^{s,t}:=\Delta_{q,+}^{K_s,K_t}$. We denote the eigenvalues of $\Delta_{q,+}^{s,t}$, counted with multiplicity and in increasing order, by $\lambda_{q,+,1}^{s,t}\leq\dots\leq\lambda_{q,+,n_q^s}^{s,t}$.

The persistent Laplacian satisfies the same relation with the persistent Betti numbers as the combinatorial Laplacian does with the Betti numbers.
\begin{theorem}[{Kernel of the persistent Laplacian, \cite[Theorem~2.7]{laplacians_algorithms}}]\label{thm:pers_lap_betti}
    For each $q\geq 0$, $\beta_q^{K,L}=\dim\ker \Delta_q^{K,L}$.
\end{theorem}
Consequently, the analysis of the persistent Laplacian has two parts. The multiplicity of the eigenvalue zero, the harmonic spectrum, describes the topological information. The non-zero eigenvalues, the non-harmonic spectrum, may contain additional information. Less is known about them, but they are believed to capture geometric rather than topological information~\cite{laplacians_survery}. Recent works relate the smallest non-zero eigenvalue to Cheeger-type constants on simplicial complexes~\cite{jost_cheeger_2024,botnan_cheeger_2026}. Although these bounds are not intuitively easy to interpret, they show that this eigenvalue is related to structural properties of the simplicial complex or simplicial pair.

\subsubsection{Computing persistent Laplacians}
The matrix of $\partial_q^K$ in the simplex basis $S_q^K$ is immediate. For the persistent Laplacian, however, we also need a matrix representation of $\partial_{q+1}^{K,L}$, and the subspace $C_{q+1}^{K,L}$ generally has no basis consisting of simplices: a chain in $C_{q+1}^{L}$ may need to be a linear combination of several simplices of $S_{q+1}^L$ for its boundary to lie in $C_q^K$. In~\cite{laplacians_algorithms}, the authors prove that the matrix of the persistent up-Laplacian $\Delta_{q,+}^{K,L}$ is a Schur complement of the matrix of the up-Laplacian of $L$.

\begin{definition}[{Generalized Schur complement~\cite{laplacians_algorithms}}]\label{def:schur_complement}
Let $M\in \RR^{n\times n}$ be a block matrix $M=\begin{pmatrix}
    A&B\\
    C&D
\end{pmatrix}$ with a square block $D\in\RR^{d\times d}$. The (generalized) Schur complement of $D$ in $M$ is
\begin{equation}
    M/D:=A-BD^\dagger C,
\end{equation}
where $D^\dagger$ denotes the Moore--Penrose pseudoinverse of $D$.
\end{definition}

To state the theorem, we need some more notation. For an integer $n$, we write $[n]:=\{1,2,\dots,n\}$. For a matrix $M\in\RR^{n\times m}$ and sets of row and column indices $I_r\subseteq [n]$ and $I_c\subseteq [m]$, we write $M(I_r,I_c)$ for the submatrix of $M$ consisting of the rows in $I_r$ and the columns in $I_c$. We order the basis of $C_q^L$ such that the $q$-simplices of $K$ come first, followed by the remaining $q$-simplices of $L$, and we write $I_q^{K,L}:=[n_q^L]\setminus[n_q^K]$ for the indices of the latter, with $I_q^{s,t}:=I_q^{K_s,K_t}$ for a filtration.

\begin{theorem}[{\cite[Theorem~4.6]{laplacians_algorithms}}]\label{thm:laplacian_schur}
    Let $K\subseteq L$ be a simplicial pair with $n_q^K<n_q^L$. Then
    \begin{equation}\label{eq:laplacian_schur}
       [\Delta_{q,+}^{K,L}]=[\Delta_{q,+}^L]\,/\,[\Delta_{q,+}^L](I_q^{K,L},I_q^{K,L}).
    \end{equation}
\end{theorem}
If $n_q^K=n_q^L$, then $C_q^K=C_q^L$, so the boundary of every $c\in C_{q+1}^L$ lies in $C_q^K$. Hence $C_{q+1}^{K,L}=C_{q+1}^L$, and the persistent up-Laplacian equals the up-Laplacian of $L$, that is, $\Delta_{q,+}^{K,L}=\Delta_{q,+}^L$.

Using this Schur complement representation, the authors of~\cite{laplacians_algorithms} show that the persistent up-Laplacian can be computed in $O\big((n_q^L)^3+n_{q+1}^L\big)$ time in the worst case. To obtain the persistent Betti number, the dimension of the kernel of the persistent Laplacian must then be determined, which can be done in $O\big((n_q^L)^3\big)$ time. This gives a total time complexity of $O\big((n_q^L)^3+n_{q+1}^L\big)$ for computing a persistent Betti number via the persistent Laplacian. They compare this with the standard algorithm for computing persistent Betti numbers, whose time complexity they state as $O\big((n_q^L)^2n_{q+1}^L+(n_{q-1}^K)^2n_q^K\big)$, and conclude that the Schur complement approach is faster whenever $n_q^L=O(n_{q+1}^L)$. They argue that this is often the case, for example for the \v{C}ech filtration.

On the other hand, the authors of~\cite{pers_lap_mnist} note that many applications require persistent Laplacians for many pairs of filtration indices and in multiple dimensions. Even when the Schur complement representation is used for each individual pair, repeating this computation over a full filtration can be time-consuming. In particular, for a filtration of $N$ complexes in which the first complex consists of a single simplex and one simplex is added at each step, they show that computing the persistent Betti numbers for all pairs via the Schur complement representation takes $O\big(N^4/(q+1)\big)$ time, where $q$ is the maximal dimension considered. This is asymptotically larger than the $O(N^3)$ time complexity of the standard matrix reduction algorithm for persistent homology.

\subsubsection{Monotonicity of the eigenvalues} % NEW heading
The authors of~\cite{laplacians_algorithms} also study the persistent up-Laplacian along a filtration and show that its eigenvalues are monotone in both filtration indices.
\begin{theorem}[{\cite[Theorem~5.2]{laplacians_algorithms}}]\label{thm:incr_eig}
    Let $K_0\subseteq K_1\subseteq\dots\subseteq K_T$ be a filtration and let $q\geq 0$. Then for all $0\leq t_1\leq t_2\leq t_3\leq T$ and $k=1,\dots,n_q^{t_1}$,
    \begin{equation}
        \lambda_{q,+,k}^{t_1,t_2}\leq \lambda_{q,+,k}^{t_1,t_3}
        \quad\text{and}\quad
        \lambda_{q,+,k}^{t_2,t_3}\leq \lambda_{q,+,k}^{t_1,t_3}.
    \end{equation}
\end{theorem}
In other words, the eigenvalues of the persistent up-Laplacian do not decrease when the end index increases or when the start index decreases. It is therefore natural to study how the persistent up-Laplacian changes between neighboring pairs of filtration indices; as we will see in Section~\ref{sec:duplicate_information}, many of these changes are zero and corresponding persistent Laplacians therefore capture the same information.

%% file: shortened_sections/newchainspaces_v4.tex
% Revised new chain spaces section, part 1.
% TODO  = information that seems to be missing.
% CHECK = places where I had to guess the intended meaning, or changed the mathematics.
% NEW   = sentences I added; remove them if you prefer.
% The commented-out older versions of paragraphs have been removed.

In Section~\ref{ch:mathematical_background}, the persistent Laplacian was defined for a simplicial pair. In practice, however, we want to analyze a whole filtration. One option is to compute the persistent Laplacian $\Delta_q^{s,t}$ for every pair $s\leq t$, but this is time-consuming and involves many repeated computations.

In this section, we introduce the new persistent chain spaces. They are designed to isolate the part of the space of persistent chains $C_{q+1}^{s,t}$ that is not already present at the neighboring pairs $(s,t-1)$, $(s-1,t)$ and $(s-1,t-1)$. In this way, they record where new information appears in the filtration, instead of storing a full basis of $C_{q+1}^{s,t}$ for every pair $(s,t)$. As a direct consequence, they show where the eigenvalues of the persistent up-Laplacian may change, and hence which persistent up-Laplacians are worth evaluating. Loosely speaking, the new persistent chain spaces summarize the changes in the geometric features captured by the persistent Laplacian.
We first motivate the definition by examining the duplicate information in the persistent Laplacians of a filtration. We then derive the new persistent chain spaces analytically and discuss their topological interpretation.

\subsection{Duplicate information in the persistent Laplacian}\label{sec:duplicate_information}
In this subsection, we give further evidence for the presence of duplicate information in the persistent Laplacians of a filtration. As noted before, the monotonicity of the eigenvalues (Theorem~\ref{thm:incr_eig}) suggests that the persistent Laplacians of neighboring pairs capture similar information. We illustrate this in two ways: visually, and by comparing the boundary operators of neighboring persistent Laplacians.

\subsubsection{Visualizing the persistent Laplacian}\label{sec:visual_pers_lap}
\begin{figure}
    \centering
    \includegraphics[width=0.7\linewidth]{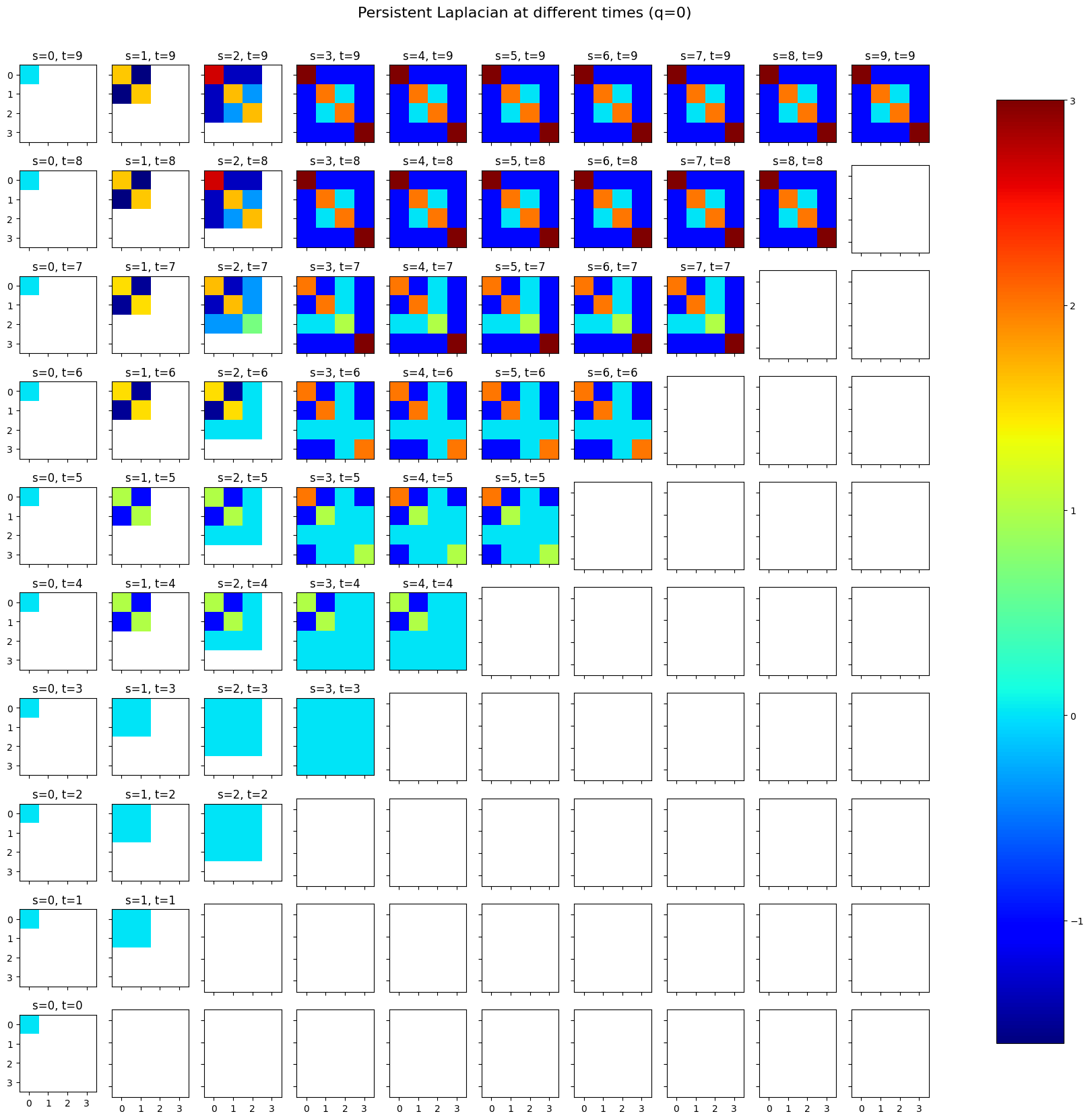}
    \caption{Heatmaps of the matrix representation $[\Delta^{s,t}_0]$ of the $0$-th persistent Laplacian for different start times $s$ and end times $t$, for the filtration in Figure~\ref{fig:filtration_example}.}
    \label{fig:pers_lap_heatmap_example_filtration}
\end{figure}

A direct way to visualize the persistent Laplacians of a filtration is to plot the matrix $[\Delta_q^{s,t}]$ as a heatmap for every pair $s\leq t$; see Figure~\ref{fig:pers_lap_heatmap_example_filtration} for dimension $0$. On the diagonal $s=t$, the persistent Laplacians coincide with the combinatorial Laplacians. Although this representation is not easy to read, it gives some intuition for how the operator behaves.

\begin{figure}
    \centering
    \includegraphics[width=\linewidth]{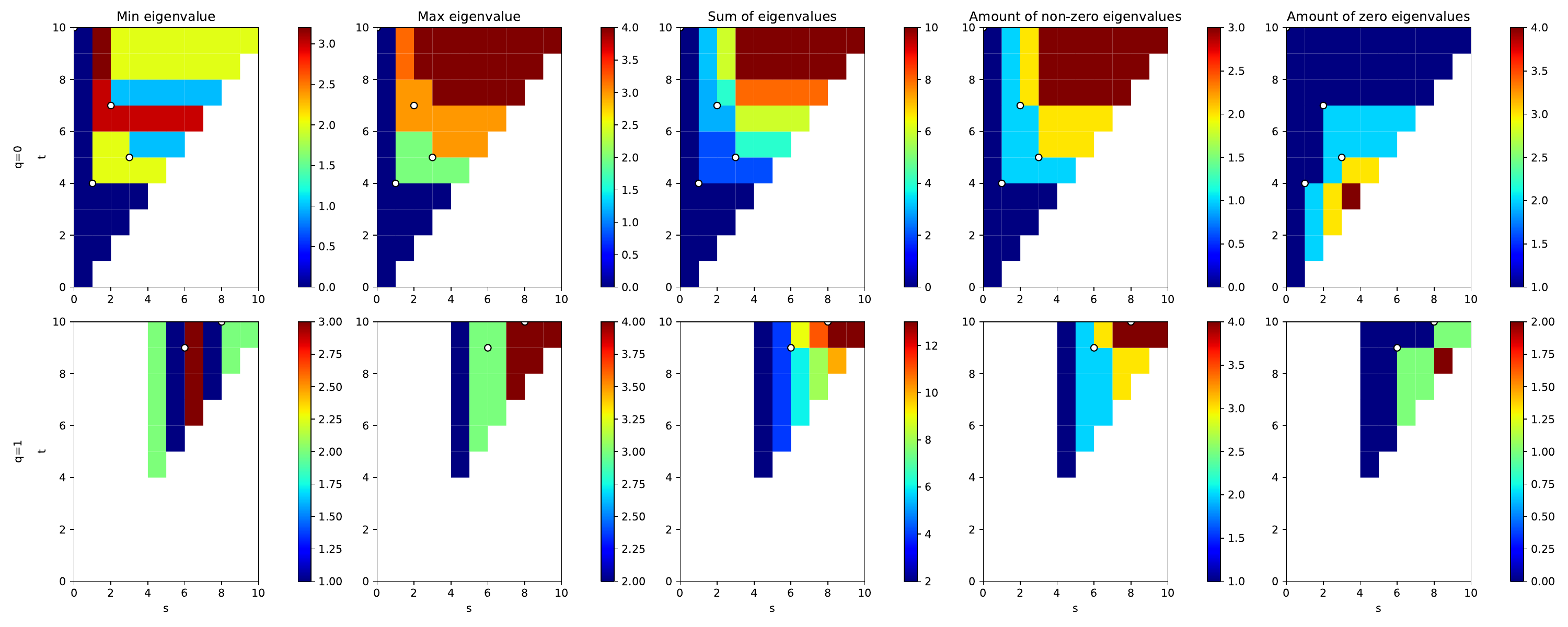}
    \caption{Aggregation functions applied to the eigenvalues of the persistent Laplacian $[\Delta^{s,t}_q]$, with the start time $s$ on the horizontal axis and the end time $t$ on the vertical axis, for the filtration in Figure~\ref{fig:filtration_example}. The Min eigenvalue column contains the minimum non-zero eigenvalue. White dots represent the bars of the persistence barcode, with the bar $[b,d)$ drawn at $(b,d)$.}
    % TODO: state which aggregation function each column shows and which dimension each row shows.
    \label{fig:pers_lap_evals_example_filtration}
\end{figure}

We are interested in the eigenvalues of these matrices rather than in their entries. For each matrix, we therefore compute the eigenvalues and condense them into a single value using an aggregation function. Figure~\ref{fig:pers_lap_evals_example_filtration} shows several such aggregation functions. In each plot, the matrix $[\Delta_q^{s,t}]$ is represented by a square whose bottom-left corner is $(s,t)$ and whose color is the value of the aggregation function applied to the eigenvalues of $\Delta_q^{s,t}$. The bars of the persistence barcode are shown as points, to highlight their effect on the persistent Laplacian.
% CHECK: the original said that the aggregation function is applied to the non-zero eigenvalues,
% but the rightmost plot shows the number of zero eigenvalues.

As an example, the point $(2,7)$ in dimension $0$ corresponds to the bar $[2,7)$ of the connected component created by vertex $2$ in the filtration of Figure~\ref{fig:filtration_example}. The square to the top right of this point shows that the persistent Laplacian $\Delta_0^{2,7}$ has two non-zero eigenvalues: the smallest is $1$ and the largest is $3$, so their sum is $4$. Moreover, the persistent Betti number is $\beta_0^{2,7}=1$.

To read Figure~\ref{fig:pers_lap_evals_example_filtration}, recall that the persistent Betti number counts the bars, with multiplicity, that are alive on the whole range from $s$ to $t$:
\[
    \beta_q^{s,t}=\#\bigl\{\text{bars } [b,d) \text{ in dimension } q:\ b\leq s \text{ and } t<d\bigr\}.
\]
Hence a bar $[b,d)$, drawn as the point $(b,d)$, contributes $1$ to $\beta_q^{s,t}$ exactly for the squares with $b\leq s$ and $t\leq d-1$, that is, the squares to the right of the point and below it, down to the diagonal $s=t$. By Theorem~\ref{thm:pers_lap_betti}, the rightmost plot, which shows the number of zero eigenvalues of $\Delta_q^{s,t}$, displays $\beta_q^{s,t}$. It is therefore a superposition of such regions, and the value of each square is the number of regions containing it. For example, in dimension $0$ the square $(3,3)$ lies in the regions of the bars $[0,\infty)$, $[1,4)$, $[2,7)$ and $[3,5)$, so $\beta_0^{3,3}=4$.

A complementary pattern appears in the plot of the number of non-zero eigenvalues. Since $\Delta_q^{s,t}$ acts on $C_q^s$, this number equals $n_q^s-\beta_q^{s,t}$, which can be written as
\[
    n_q^s-\beta_q^{s,t}
    =\big(n_q^s-\beta_q^{s,s}\big)
    +\#\bigl\{\text{bars } [b,d) \text{ in dimension } q:\ b\leq s<d\leq t\bigr\}.
\]
The first term is the number of non-zero eigenvalues of the combinatorial Laplacian $\Delta_q^{K_s}$, shown on the diagonal. For fixed $s$, the number of non-zero eigenvalues therefore increases by one for each bar born at or before $s$ whose death value $d$ is passed by $t$: the corresponding cycle becomes a boundary in $K_t$, and a zero eigenvalue becomes positive. Thus each finite bar $[b,d)$ now affects the squares above and to the right of the point $(b,d)$, namely those with $b\leq s<d\leq t$; for $s\geq d$, its effect is already contained in the diagonal value.
For instance, for $s=2$ in dimension $0$, the diagonal value is $0$, and the number of non-zero eigenvalues is $0$ for $t<4$, $1$ for $4\leq t<7$, after the death of $[1,4)$, and $2$ for $t\geq 7$, after the death of $[2,7)$. The same regions are visible in the plot of the sum of the eigenvalues, but there the changes are not simply additive, since the values of the non-zero eigenvalues also change with $s$ and $t$.

\subsubsection{Boundary operators of neighboring persistent Laplacians}\label{sec:bdry_subsequent_lap}
We first compare persistent Laplacians with consecutive end times, $\Delta_q^{s,t-1}$ and $\Delta_q^{s,t}$. By Definition~\ref{def:persistent_laplacian}, they have the same down-Laplacian $\Delta_{q,-}^{K_s}$, so they can only differ in their persistent up-Laplacians.

These up-Laplacians are defined using $\partial_{q+1}^{s,t-1}$ and $\partial_{q+1}^{s,t}$, which act on $C_{q+1}^{s,t-1}$ and $C_{q+1}^{s,t}$, respectively, and we have $C_{q+1}^{s,t-1}\subseteq C_{q+1}^{s,t}$. With respect to the inner product induced by the simplex basis, $C_{q+1}^{s,t}$ therefore decomposes orthogonally into $C_{q+1}^{s,t-1}$ and its orthogonal complement $C_{q+1}^{s,t}\cap\big(C_{q+1}^{s,t-1}\big)^\perp$ within $C_{q+1}^{s,t}$. Accordingly, every chain $c\in C_{q+1}^{s,t}$ can be written uniquely as $c=c_t+c_{t-1}$ with $c_{t-1}\in C_{q+1}^{s,t-1}$ and $c_t\in C_{q+1}^{s,t}\cap\big(C_{q+1}^{s,t-1}\big)^\perp$.

Now let $a\in C_q^s$ and $A:=\big(\partial_{q+1}^{s,t}\big)^*a\in C_{q+1}^{s,t}$. Since $\partial_{q+1}^{s,t}$ restricts to $\partial_{q+1}^{s,t-1}$ on $C_{q+1}^{s,t-1}$, the component $A_{t-1}$ equals $(\partial_{q+1}^{s,t-1})^*a$. Indeed, for every $b\in C_{q+1}^{s,t-1}$,
\[
    \langle A_{t-1},b\rangle
    =
    \langle A,b\rangle
    =
    \langle a,\partial_{q+1}^{s,t}b\rangle
    =
    \langle a,\partial_{q+1}^{s,t-1}b\rangle,
\]
and since $A_{t-1}\in C_{q+1}^{s,t-1}$, this characterizes $A_{t-1}$ as $(\partial_{q+1}^{s,t-1})^*a$. Consequently,
\begin{equation*}
    \Delta_{q,+}^{s,t}a=\partial_{q+1}^{s,t}(A_t+A_{t-1})=\partial_{q+1}^{s,t}A_t+\Delta_{q,+}^{s,t-1}a.
\end{equation*}
In other words, writing $\pi$ for the orthogonal projection of $C_{q+1}^{s,t}$ onto $C_{q+1}^{s,t}\cap\big(C_{q+1}^{s,t-1}\big)^\perp$, we have
\[
    \Delta_{q,+}^{s,t}=\Delta_{q,+}^{s,t-1}+\partial_{q+1}^{s,t}\,\pi\,\big(\partial_{q+1}^{s,t}\big)^*,
\]
so $\Delta_{q,+}^{s,t}$ only adds the contribution of the chains in the orthogonal complement of $C_{q+1}^{s,t-1}$. Since this additional term is positive semidefinite, this also explains the first inequality in Theorem~\ref{thm:incr_eig} and is made precise in Theorem \ref{thm:horizontal_operator}.

A similar argument applies to $\Delta_q^{s-1,t}$ and $\Delta_q^{s,t}$, whose up-Laplacians are defined using $\partial_{q+1}^{s-1,t}$ and $\partial_{q+1}^{s,t}$, respectively. Again, $C_{q+1}^{s-1,t}\subseteq C_{q+1}^{s,t}$: if a chain in $C_{q+1}^t$ has its boundary in $C_q^{s-1}$, then its boundary also lies in $C_q^s\supseteq C_q^{s-1}$. Moreover, viewed as a map into $C_q^s$, $\partial_{q+1}^{s-1,t}$ is the restriction of $\partial_{q+1}^{s,t}$ to $C_{q+1}^{s-1,t}$.
In this case, however, the two up-Laplacians act on different spaces, $C_q^{s-1}$ and $C_q^s$. Let $\iota:C_q^{s-1}\to C_q^s$ denote the inclusion, so that its adjoint $\iota^*$ is the orthogonal projection onto $C_q^{s-1}$. For $a\in C_q^s$, write again $A:=\big(\partial_{q+1}^{s,t}\big)^*a=A_s+A_{s-1}$, with $A_{s-1}\in C_{q+1}^{s-1,t}$ and $A_s$ in its orthogonal complement within $C_{q+1}^{s,t}$. For every $b\in C_{q+1}^{s-1,t}$, the boundary $\partial_{q+1}^{s,t}b$ lies in $C_q^{s-1}$, so
\[
    \langle A_{s-1},b\rangle
    =
    \langle a,\partial_{q+1}^{s,t}b\rangle
    =
    \langle \iota^*a,\partial_{q+1}^{s-1,t}b\rangle,
\]
which shows that $A_{s-1}=(\partial_{q+1}^{s-1,t})^*\iota^*a$. As before, we obtain
\[
    \Delta_{q,+}^{s,t}=\iota\,\Delta_{q,+}^{s-1,t}\,\iota^*+\partial_{q+1}^{s,t}\,\pi'\,\big(\partial_{q+1}^{s,t}\big)^*,
\]
where $\pi'$ is the orthogonal projection of $C_{q+1}^{s,t}$ onto $C_{q+1}^{s,t}\cap\big(C_{q+1}^{s-1,t}\big)^\perp$. Hence, also in this case, the persistent up-Laplacian $\Delta_{q,+}^{s,t}$ only adds information. The down-Laplacians $\Delta_{q,-}^{K_{s-1}}$ and $\Delta_{q,-}^{K_s}$, however, can differ. This is again made more precise in Theorem \ref{thm:vertical_operator}.

% Revised subsection: separating the two effects.
% TODO  = information that seems to be missing.
% CHECK = places where I had to guess the intended meaning, or changed the mathematics.
% NEW   = sentences I added; remove them if you prefer.
% The commented-out older versions (including the old diagram and the discussion of the points
% (1,6), (1,8), (2,8), (3,6), (3,7)) have been removed.

\subsection{Separating the two effects}\label{sec:effect_separation}
In Section~\ref{sec:duplicate_information}, we saw that the duplicate information in the persistent up-Laplacians of a filtration has two sources, corresponding to the two inequalities in Theorem~\ref{thm:incr_eig}: $\Delta_{q,+}^{s,t}$ shares information with $\Delta_{q,+}^{s,t-1}$, the square below it, and with $\Delta_{q,+}^{s-1,t}$, the square to its left. In this subsection, we remove these two sources separately. The \emph{horizontal operator} removes the information shared with the square below, which we call the \emph{vertical effect}. Similarly, the \emph{vertical operator} removes the information shared with the square to the left, which we call the \emph{horizontal effect}. This way the horizontal operator describes only the horizontal effect and the vertical operator describes only the vertical effect.

% In Section~\ref{sec:duplicate_information}, we saw that the duplicate information in the persistent up-Laplacians of a filtration has two sources, corresponding to the two inequalities in Theorem~\ref{thm:incr_eig}: $\Delta_{q,+}^{s,t}$ shares information with $\Delta_{q,+}^{s,t-1}$, the square below it, and with $\Delta_{q,+}^{s-1,t}$, the square to its left. We call these the \emph{vertical effect} and the \emph{horizontal effect}, respectively. In this subsection, we remove them separately. The \emph{horizontal operator} removes the vertical effect, so that only the horizontal effect remains; similarly, the \emph{vertical operator} removes the horizontal effect.

\subsubsection{Removing the vertical effect: the horizontal operator}\label{sec:horizontal_operator}
To understand the vertical effect, consider the diagram below, in which the maps defining the persistent Laplacians at $(s,t-1)$ and $(s,t)$ are shown in red and blue, respectively. Both persistent Laplacians have the same down-Laplacian, so the corresponding arrows belong to both, but they are only drawn in red. Recall from Section~\ref{sec:bdry_subsequent_lap} that $C_{q+1}^{s,t-1}\subseteq C_{q+1}^{s,t}$ and that $\partial_{q+1}^{s,t-1}$ is the restriction of $\partial_{q+1}^{s,t}$ to $C_{q+1}^{s,t-1}$. Note that the diagram is not commutative as a whole: $(\partial_{q+1}^{s,t})^*$ differs from $(\partial_{q+1}^{s,t-1})^*$ followed by the inclusion, and this difference is exactly what we study below.

\[
\begin{tikzcd}[column sep=large,row sep=large]
C_{q+1}^{s-1} \arrow[d, hook]
  &
  {}
  &
C_{q}^{s-1} \arrow[d, hook]
  &
  {}
  \\
C_{q+1}^{s} \arrow[d, hook]
&
  C_{q+1}^{s,t-1} \arrow[red, shift left=1ex]{r}{\partial_{q+1}^{s,t-1}} \arrow[d, hook] \arrow[hook]{dl}
  &
C_{q}^{s} \arrow[d, hook] \arrow[red, shift left=1ex]{r}{\partial_{q}^{s}} \arrow[red, shift left=1ex]{l}{(\partial_{q+1}^{s,t-1})^*} \arrow[blue, shift left=1ex, pos=0.8]{dl}{(\partial_{q+1}^{s,t})^*}
  &
C_{q-1}^{s} \arrow[red, shift left=1ex]{l}{(\partial_{q}^{s})^*}
\\
C_{q+1}^{t-1} \arrow[d, hook]
  &
C_{q+1}^{s,t} \arrow[blue, shift left=1ex, pos=0.2]{ur}{\partial_{q+1}^{s,t}} \arrow[hook]{dl}
  &
C_{q}^{t-1} \arrow[d, hook]
  & {}
  \\
C_{q+1}^{t} \arrow{rr}{\partial_{q+1}^t}
  &
{}
  &
C_{q}^{t}
  & {}
\end{tikzcd}
\]

Recall the intuition behind a bar $[s,t)$ in the persistence barcode. It corresponds to a $q$-cycle in $K_s$ that becomes a boundary in $K_t$, that is, the boundary of a $(q+1)$-chain of $K_t$. The bar ends exactly at $t$ if this $(q+1)$-chain requires a $(q+1)$-simplex that appears at $t$, so that the cycle is not yet a boundary in $K_{t-1}$. If the cycle is already a boundary in $K_{t-1}$, the bar is shorter.

We approach the geometric information in the same way, ignoring the down-Laplacian for now, and consider the $(q+1)$-chains in $C_{q+1}^{s,t}$, on which the persistent up-Laplacian is built. Since $C_{q+1}^{s,t-1}\subseteq C_{q+1}^{s,t}$, some of these chains were already captured by $\Delta_q^{s,t-1}$. To focus on the new information only, we instead consider
\[
    \widebar{C_{q+1}^{s,t}}:=C_{q+1}^{s,t}\cap\big(C_{q+1}^{s,t-1}\big)^\perp,
\]
the orthogonal complement of $C_{q+1}^{s,t-1}$ in $C_{q+1}^{s,t}$, and the restricted boundary map $\bar\partial_{q+1}^{s,t}:=\partial_{q+1}^{s,t}|_{\widebar{C_{q+1}^{s,t}}}:\widebar{C_{q+1}^{s,t}}\to C_q^s$. For $c\in C_q^s$, the chain $(\bar\partial_{q+1}^{s,t})^*c$ is the orthogonal projection of $(\partial_{q+1}^{s,t})^*c$ onto $\widebar{C_{q+1}^{s,t}}$, that is, the part of $(\partial_{q+1}^{s,t})^*c$ that is orthogonal to all $(q+1)$-chains already available at $t-1$.
Using this restricted boundary map in the same way as for the persistent up-Laplacian, we define the \emph{horizontal operator}
\begin{equation}
    \mathcal{H}_q^{s,t}:=\bar\partial_{q+1}^{s,t}(\bar\partial_{q+1}^{s,t})^*:C_q^s\to C_q^s.
\end{equation}
It can be computed directly from persistent up-Laplacians.

\begin{theorem}\label{thm:horizontal_operator}
    For $0\leq s<t\leq T$ and $q\geq 0$, the horizontal operator satisfies
    \begin{equation}
        \mathcal{H}_q^{s,t}=\Delta_{q,+}^{s,t}-\Delta_{q,+}^{s,t-1}.
    \end{equation}
    In particular, $[\mathcal{H}_q^{s,t}]=[\Delta_{q,+}^{s,t}]-[\Delta_{q,+}^{s,t-1}]$.
\end{theorem}
The proof of the theorem can be found in Appendix \ref{sec:proof_horizontal_vertical}.

\begin{figure}
    \centering
    \includegraphics[width=0.75\linewidth]{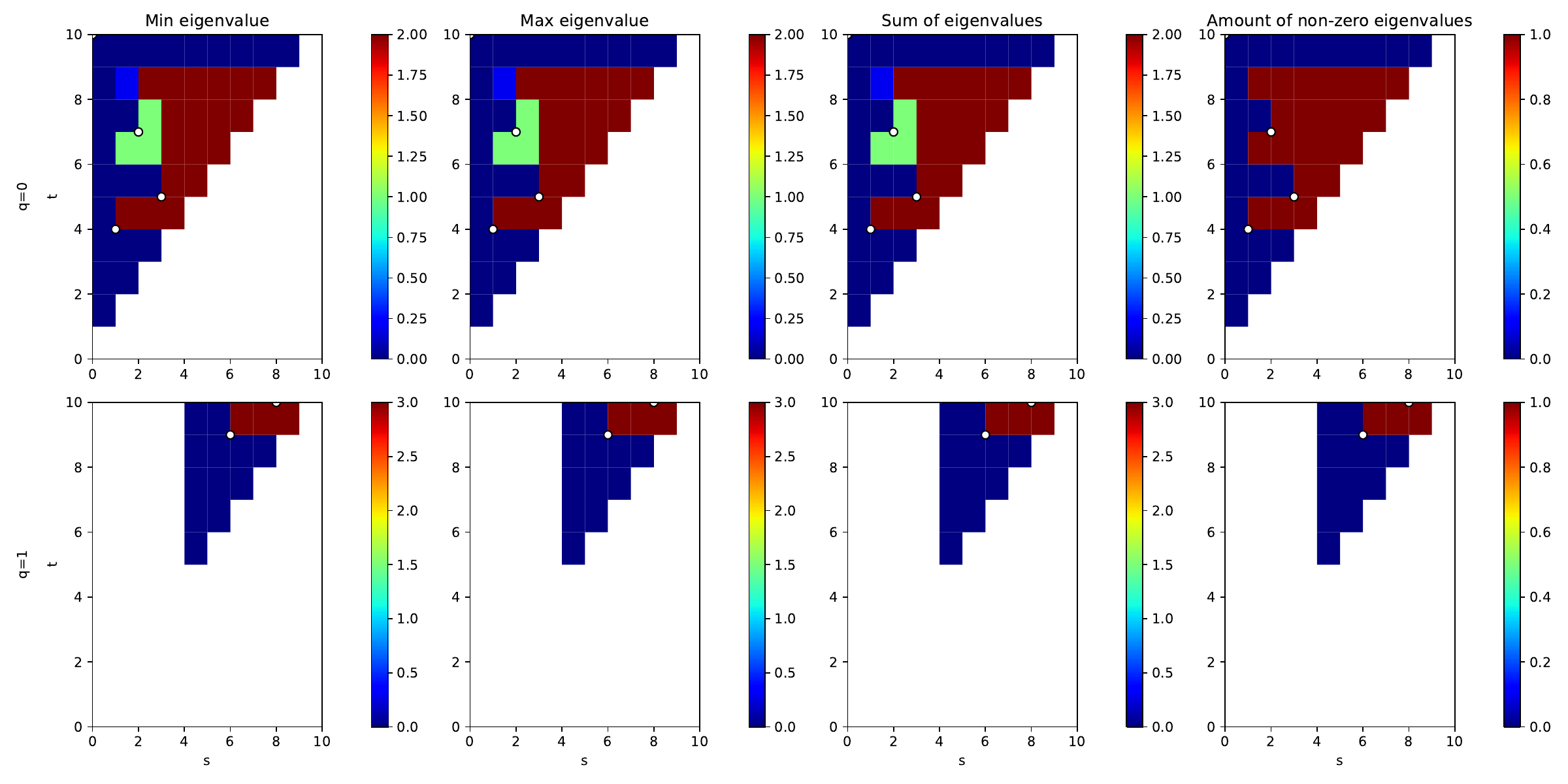}
    \caption{Aggregation functions applied to the eigenvalues of the horizontal operator $\mathcal{H}_q^{s,t}$, with the start time $s$ on the horizontal axis and the end time $t$ on the vertical axis, for the filtration in Figure~\ref{fig:filtration_example}. White dots represent the bars of the persistence barcode, with the bar $[b,d)$ drawn at $(b,d)$.}
    \label{fig:horizontal_laplacian}
\end{figure}

As in Section~\ref{sec:visual_pers_lap}, we can plot aggregation functions of the eigenvalues of $\mathcal{H}_q^{s,t}$; see Figure~\ref{fig:horizontal_laplacian}. Compared with Figure~\ref{fig:pers_lap_evals_example_filtration}, the number of non-zero eigenvalues is reduced: for every square with bottom-left corner $(s,t)$, the effect of the square below it, corresponding to $(s,t-1)$, has been removed. We can now clearly see which features are new at end time $t$. However, these features may already be present for earlier start times. For example, the square corresponding to $(2,4)$ still has a non-zero eigenvalue equal to $2$, although this effect originates at $(1,4)$.

\subsubsection{Removing the horizontal effect: the vertical operator}\label{sec:vertical_operator}
For the horizontal effect, we consider the diagram below. As before, the maps defining the persistent Laplacians at $(s-1,t)$ and $(s,t)$ are shown in red (dashed) and blue, respectively. In this case, the down-Laplacians at $(s-1,t)$ and $(s,t)$ differ, so both are drawn. Recall from Section~\ref{sec:bdry_subsequent_lap} that $C_{q+1}^{s-1,t}\subseteq C_{q+1}^{s,t}$. Since the codomains of the two boundary maps differ, we use the inclusion $\iota:C_q^{s-1}\hookrightarrow C_q^s$, so that $\iota\circ\partial_{q+1}^{s-1,t}$ is the restriction of $\partial_{q+1}^{s,t}$ to $C_{q+1}^{s-1,t}$.

\[
\begin{tikzcd}[column sep=large,row sep=large]
C_{q+1}^{s-1} \arrow[d, hook]
  &
 {}
  &
  {}
  &
C_{q}^{s-1} \arrow[d, hook] \arrow[red, shift left=1ex, dashed]{dll}{(\partial_{q+1}^{s-1,t})^*} \arrow[red, shift left=1ex, dashed]{r}{\partial_{q}^{s-1}}
  &
C_{q-1}^{s-1} \arrow[red, shift left=1ex, dashed]{l}{(\partial_{q}^{s-1})^*} \arrow[hook]{d}
  \\
C_{q+1}^{s} \arrow[d, hook]
  &
C_{q+1}^{s-1,t} \arrow[red, shift left=1ex, dashed]{urr}{\partial_{q+1}^{s-1,t}} \arrow[dr, hook, dashed] \arrow[hook]{ddl}
&
 {}
  &
C_{q}^{s} \arrow[d, hook] \arrow[blue, shift left=1ex]{r}{\partial_{q}^{s}} \arrow[blue, shift left=1ex, pos=0.8]{dl}{(\partial_{q+1}^{s,t})^*}
  &
C_{q-1}^{s} \arrow[blue, shift left=1ex]{l}{(\partial_{q}^{s})^*} \\
C_{q+1}^{t-1} \arrow[d, hook]
  &
{}
  &
C_{q+1}^{s,t} \arrow[blue, shift left=1ex]{ur}{\partial_{q+1}^{s,t}} \arrow[hook]{dll}
  &
C_{q}^{t-1} \arrow[d, hook]
  & {} \\
C_{q+1}^{t}
  &
  {}
  &
{}
  &
C_{q}^{t}
  & {}
\end{tikzcd}
\]

Besides requiring that the cycle becomes a boundary exactly at $t$, a bar $[s,t)$ also requires that the cycle appears at $s$ and not before; otherwise, the bar would start earlier and be longer. Translating this to the persistent Laplacian is harder, since both the persistent up-Laplacian and the down-Laplacian depend on $s$, so features that are new at $s$ can be captured by both.

We first focus on the persistent up-Laplacian, for which the discussion is very similar to the horizontal case, discussed in Section \ref{sec:horizontal_operator}. We consider
\[
    \widetilde{C_{q+1}^{s,t}}:=C_{q+1}^{s,t}\cap\big(C_{q+1}^{s-1,t}\big)^\perp,
\]
the orthogonal complement of $C_{q+1}^{s-1,t}$ in $C_{q+1}^{s,t}$, and the restricted boundary map $\tilde\partial_{q+1}^{s,t}:=\partial_{q+1}^{s,t}|_{\widetilde{C_{q+1}^{s,t}}}:\widetilde{C_{q+1}^{s,t}}\to C_q^s$. As before, we define the \emph{vertical operator}
\begin{equation}
    \mathcal{V}_q^{s,t}:=\tilde\partial_{q+1}^{s,t}(\tilde\partial_{q+1}^{s,t})^*:C_q^s\to C_q^s,
\end{equation}
which can again be computed from persistent up-Laplacians.

\begin{theorem}\label{thm:vertical_operator}
    For $1\leq s\leq t\leq T$ and $q\geq 0$, the vertical operator satisfies
    \begin{equation}
        \mathcal{V}_q^{s,t}=\Delta_{q,+}^{s,t}-\iota\,\Delta_{q,+}^{s-1,t}\,\iota^*,
    \end{equation}
    where $\iota:C_q^{s-1}\hookrightarrow C_q^s$ denotes the inclusion. In particular, $[\mathcal{V}_q^{s,t}]=[\Delta_{q,+}^{s,t}]-[\iota\,\Delta_{q,+}^{s-1,t}\,\iota^*]$.
\end{theorem}
The proof can be found in Appendix \ref{sec:proof_horizontal_vertical}. Since the simplices of $K_{s-1}$ come first in the basis of $C_q^s$, the matrix $[\iota\,\Delta_{q,+}^{s-1,t}\,\iota^*]$ is obtained from $[\Delta_{q,+}^{s-1,t}]$ by adding zero rows and columns for the simplices of $K_s\setminus K_{s-1}$.

For the down-Laplacian, the same reasoning is less clear. In homology, the down-Laplacian is used to check whether a chain $c$ is a cycle, that is, whether $\partial_q^s c=0$. For the persistence barcode, this only matters if the cycle appears at $s$ and not before, so we would require $\partial_q^s c=0$ and $c\notin\iota\big(C_q^{s-1}\big)$.
% CHECK: the original condition was \partial_q^{s-1}\circ\iota(c)\neq 0, which does not type-check for c in C_q^s (see notes).
On the other hand, the geometric information in the down-Laplacian comes from chains $c$ with non-zero boundary. Since $\partial_q^s$ restricts to $\partial_q^{s-1}$ on $C_q^{s-1}$, some duplicate information is again present. However, restricting the down-Laplacian to the orthogonal complement of $C_q^{s-1}$ would change the space on which it acts. Instead, one could use $\hat\partial_q^s:=\partial_q^s\circ P$, where $P$ is the orthogonal projection of $C_q^s$ onto the orthogonal complement of $C_q^{s-1}$; this map is still defined on all of $C_q^s$. However, $\hat\partial_q^s\circ\partial_{q+1}^s$ need not be zero. As a consequence, the images of the up- and down-parts are no longer orthogonal, and the eigenvalues no longer split into an up-part and a down-part.

Nevertheless, the spectral information of the down-Laplacian is already captured by the up-Laplacian of one dimension lower on the diagonal. Indeed, $\Delta_{q+1,-}^{K_s}=(\partial_{q+1}^s)^*\partial_{q+1}^s$ and $\Delta_{q,+}^{K_s}=\partial_{q+1}^s(\partial_{q+1}^s)^*=\Delta_{q,+}^{s,s}$ have the same non-zero eigenvalues, counted with multiplicity, and $\partial_{q+1}^s$ maps eigenvectors of the former with non-zero eigenvalues to eigenvectors of the latter.
Therefore, in the rest of the paper we only consider the persistent up-Laplacian, unless stated otherwise.

\begin{figure}
    \centering
    \includegraphics[width=0.75\linewidth]{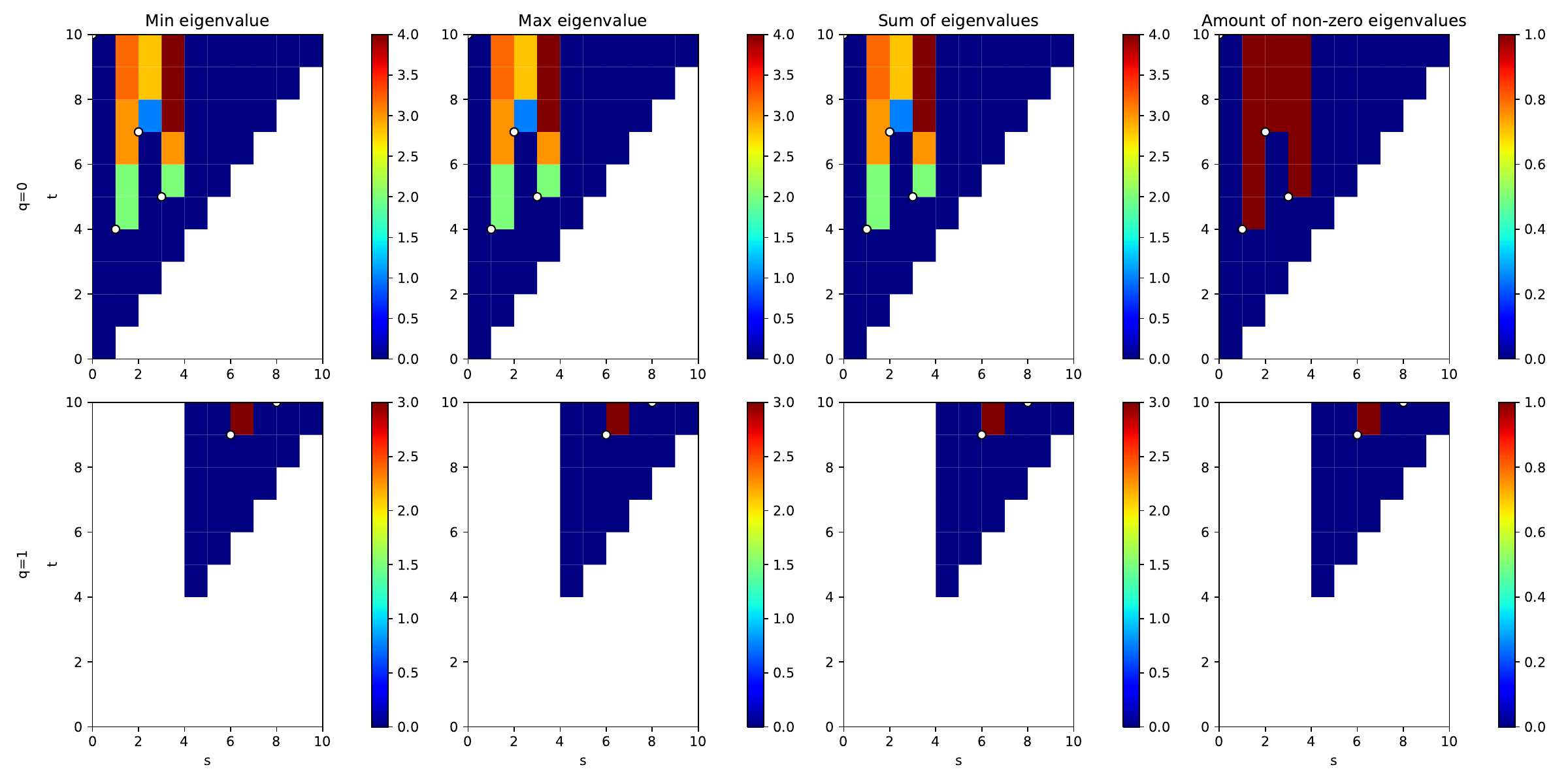}
    \caption{Aggregation functions applied to the eigenvalues of the vertical operator $\mathcal{V}_q^{s,t}$, with the start time $s$ on the horizontal axis and the end time $t$ on the vertical axis, for the filtration in Figure~\ref{fig:filtration_example}. White dots represent the bars of the persistence barcode, with the bar $[b,d)$ drawn at $(b,d)$.}
    \label{fig:vertical_laplacian}
\end{figure}

By definition, both $\mathcal{H}_q^{s,t}$ and $\mathcal{V}_q^{s,t}$ are positive semidefinite, which is not immediately obvious from the differences in Theorems~\ref{thm:horizontal_operator} and~\ref{thm:vertical_operator}. We plot the eigenvalues of the vertical operator in the same way as before; see Figure~\ref{fig:vertical_laplacian}. Now, for every square, the effect of the square to its left, corresponding to $(s-1,t)$, has been removed, instead of the effect of the square below it. All the previously discussed points at which the persistent up-Laplacian changes are also visible here.
The goal is now to combine both operations and remove both the horizontal and the vertical effect, to obtain an operator that is non-zero only at the pairs where the persistent up-Laplacian changes.

% Revised subsection: new persistent chain spaces.
% TODO  = information that seems to be missing.
% CHECK = places where I had to guess the intended meaning, or changed the mathematics.
% NEW   = sentences I added; remove them if you prefer.
% The commented-out older versions have been removed.

\subsection{New persistent chain spaces}
So far, we have discussed how to remove the vertical effect, coming from $t-1$, and the horizontal effect, coming from $s-1$. The next step is to remove both at the same time. This is analogous to the persistence barcode, where a bar $[s,t)$ requires both that a cycle appears at $s$ and that it becomes a boundary at $t$. To this end, we decompose $C_{q+1}^{s,t}$ into an old vertical, an old horizontal, an old ``diagonal'' part, and a remaining new part. In this subsection and for the rest of the paper, we use the conventions $C_{q}^{-1,t}:=\{0\}$ and $C_q^{s,s-1}:=\{0\}$.

We first show that the persistent chains at the two neighboring pairs $(s-1,t)$ and $(s,t-1)$ only overlap in the persistent chains at $(s-1,t-1)$. As a consequence, the horizontal effect coming from the left and the vertical effect coming from below do not overlap.
 
\begin{lemma}\label{lemma:horizontal_vertical_distinct}
    For $q\geq 0$ and $1\leq s\leq t$, we have
    \begin{enumerate}
        \item[(i)] $C_{q+1}^{s-1,t}\cap C_{q+1}^{s,t-1}\subseteq C_{q+1}^{s-1,t-1}$, with equality if $s\leq t-1$.
        \item[(ii)] $\widebar{C_{q+1}^{s-1,t}}\cap\widetilde{C_{q+1}^{s,t-1}}=\{0\}$.
    \end{enumerate}
\end{lemma}

Although Lemma~\ref{lemma:horizontal_vertical_distinct} shows that these subspaces intersect trivially, they are in general not orthogonal, which is why, in Section~\ref{ch:algorithm}, the new persistent chain spaces have to be orthogonalized with respect to each other whenever a persistent Laplacian is computed. Nevertheless, the lemma yields a decomposition of $C_{q+1}^{s,t}$. By definition, both $\widebar{C_{q+1}^{s-1,t}}$ and $\widetilde{C_{q+1}^{s,t-1}}$ are orthogonal to $C_{q+1}^{s-1,t-1}$, and all three spaces are subspaces of $C_{q+1}^{s,t}$. Together with Lemma~\ref{lemma:horizontal_vertical_distinct}, this implies that their sum is direct: if $a+b+c=0$ with $a$, $b$ and $c$ in these three spaces, respectively, then projecting onto $C_{q+1}^{s-1,t-1}$ gives $c=0$, and hence $a=-b=0$ by the lemma. % NEW: argument for directness
Decomposing $C_{q+1}^{s,t}$ into this direct sum and its orthogonal complement gives
\begin{equation}\label{eq:new_chain_decomposition}
    C_{q+1}^{s,t}=\widebar{C_{q+1}^{s-1,t}}\oplus \widetilde{C_{q+1}^{s,t-1}}\oplus C_{q+1}^{s-1,t-1}\oplus \Big(C_{q+1}^{s,t}\cap\big(\widebar{C_{q+1}^{s-1,t}}\oplus \widetilde{C_{q+1}^{s,t-1}}\oplus C_{q+1}^{s-1,t-1}\big)^\perp\Big).
\end{equation}
Here and in the rest of the paper, $\oplus$ denotes an internal direct sum, which is not necessarily orthogonal. % NEW

\begin{figure}
    \centering
    \includegraphics[width=0.3\linewidth]{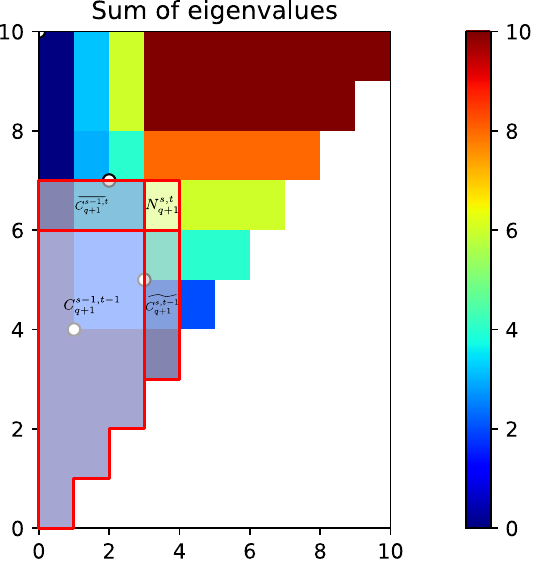}
    \caption{Regions of the plot whose information is captured by the different parts of the decomposition~\eqref{eq:new_chain_decomposition}, for $(s,t)=(3,6)$. The background shows the trace of the persistent Laplacian.}
    \label{fig:decomposition_example}
\end{figure}

Each part of this decomposition corresponds to a region of squares in the plot of the persistent Laplacian; see Figure~\ref{fig:decomposition_example}. The space $\widebar{C_{q+1}^{s-1,t}}$ captures the effect coming from the left of the square $(s,t)$, the space $\widetilde{C_{q+1}^{s,t-1}}$ the effect coming from below, and $C_{q+1}^{s-1,t-1}$ the effect of all squares to the bottom left. The remaining orthogonal complement can therefore be regarded as the ``new'' effect at $(s,t)$.

\begin{definition}[New persistent chain space]\label{def:new_chain_space}
    The new persistent chain space $N_{q+1}^{s,t}\subseteq C_{q+1}^{s,t}$ is defined by
    \begin{equation}
        N_{q+1}^{s,t}:=C_{q+1}^{s,t}\cap\big(\widebar{C_{q+1}^{s-1,t}}\oplus \widetilde{C_{q+1}^{s,t-1}}\oplus C_{q+1}^{s-1,t-1}\big)^\perp.
    \end{equation}
\end{definition}
% NEW, CHECK: simpler equivalent description (see notes).
Since $\widebar{C_{q+1}^{s-1,t}}\oplus C_{q+1}^{s-1,t-1}=C_{q+1}^{s-1,t}$ and $\widetilde{C_{q+1}^{s,t-1}}\oplus C_{q+1}^{s-1,t-1}=C_{q+1}^{s,t-1}$ for $s<t$, the new persistent chain space is the orthogonal complement in $C_{q+1}^{s,t}$ of everything that is already present at the two neighboring pairs:
\[
    N_{q+1}^{s,t}=C_{q+1}^{s,t}\cap\big(C_{q+1}^{s-1,t}+C_{q+1}^{s,t-1}\big)^\perp.
\]
This also holds for $s=t$, since then $C_{q+1}^{s,s-1}=\{0\}$ and $C_{q+1}^{s-1,s-1}\subseteq C_{q+1}^{s-1,s}$.

\begin{figure}
    \centering
    \includegraphics[width=\linewidth]{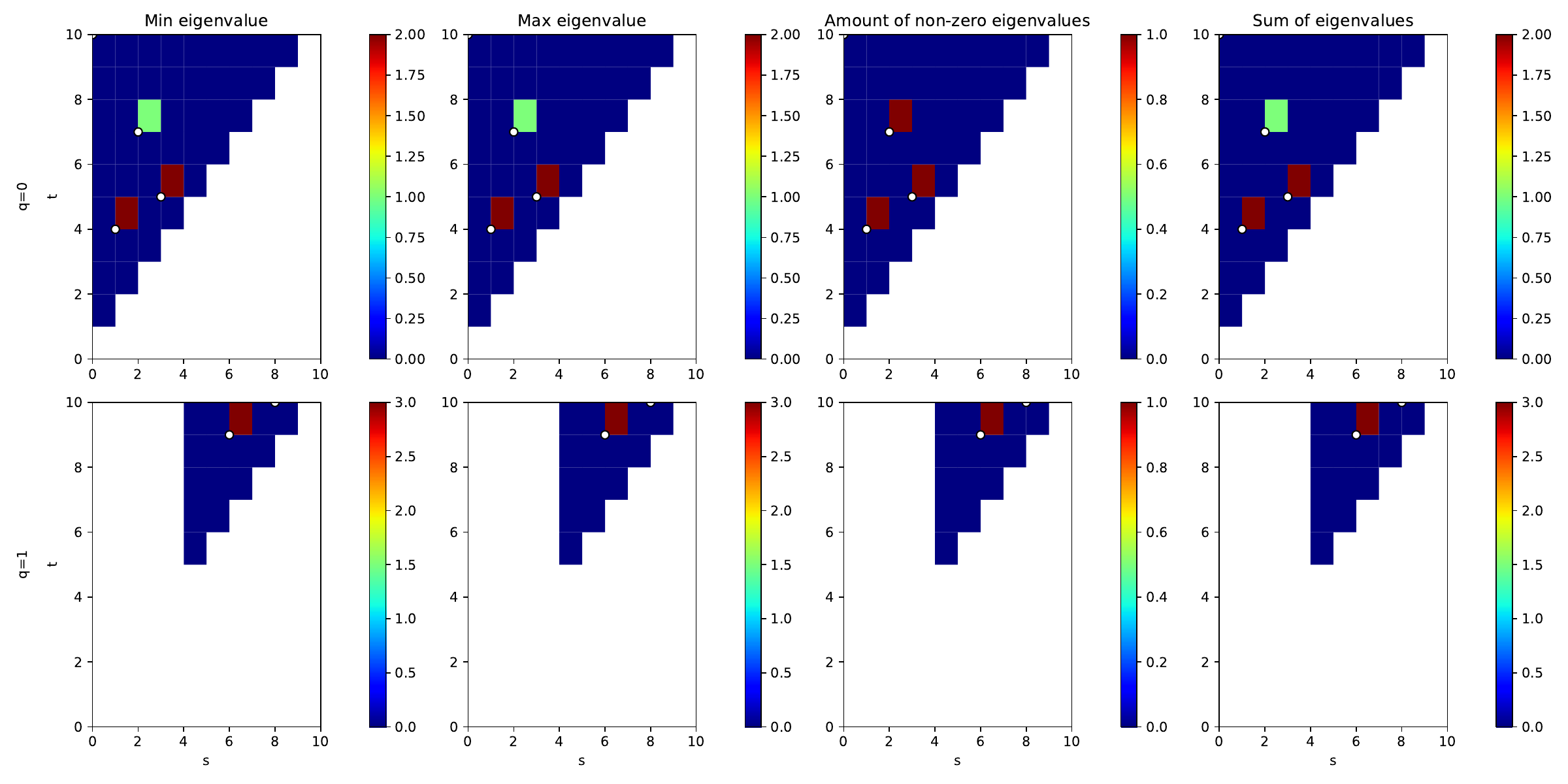}
    \caption{Aggregation functions applied to the eigenvalues of $\mathcal{R}_q^{s,t}$, with the start time $s$ on the horizontal axis and the end time $t$ on the vertical axis, for the filtration in Figure~\ref{fig:filtration_example}. White dots represent the bars of the persistence barcode, with the bar $[b,d)$ drawn at $(b,d)$.}
    % TODO: state which aggregation function each column shows and which dimension each row shows.
    \label{fig:full_multi_example}
\end{figure}

Let $\partial_{N_{q+1}^{s,t}}:=\partial_{q+1}^{s,t}|_{N_{q+1}^{s,t}}:N_{q+1}^{s,t}\to C_q^s$ and $\mathcal{R}_q^{s,t}:=\partial_{N_{q+1}^{s,t}}\big(\partial_{N_{q+1}^{s,t}}\big)^*:C_q^s\to C_q^s$.
% CHECK: the original restricted \partial_{q+1}^t, whose codomain is C_q^t (see notes).
Figure~\ref{fig:full_multi_example} shows aggregation functions of the eigenvalues of $\mathcal{R}_q^{s,t}$ for the example filtration. Clearly, a lot of information is lost this way: only the squares at the points of the persistence barcode remain non-zero. These are exactly the pairs at which $N_{q+1}^{s,t}$ contains chains with non-zero boundary.
% CHECK: is this an observation for this example, or a general result? If general, refer to where it is proven.
At other pairs, $N_{q+1}^{s,t}$ may be non-zero but contain only chains with zero boundary. For example, a $1$-cycle is born at time $6$. Its boundary is zero, so it lies in $C_1^{0,6}$, and since $C_1^{-1,6}=\{0\}$ by our convention, it contributes to $N_1^{0,6}$.

\subsubsection{Recovering all information of the persistent Laplacian}
The space $N_{q+1}^{s,t}$ contains the chains that are new at $(s,t)$; one can view them as being born at $(s,t)$. These chains remain part of $C_{q+1}^{s',t'}$ for all $s'\geq s$ and $t'\geq t$. Hence, to compute the persistent up-Laplacian at such a pair, they have to be combined with the chains born at the other pairs to the bottom left of $(s',t')$. The following theorem makes this precise.

\begin{theorem}\label{thm:new_chain_decomposition}
    For $0\leq s\leq t\leq T$,
    \[
        C_{q+1}^{s,t}=\bigoplus_{\substack{0\leq s^-\leq s,\ t^-\leq t,\\ s^-\leq t^-}} N_{q+1}^{s^-,t^-}.
    \]
\end{theorem}
% CHECK: I stated the theorem directly for C^{s,t}, dropping S^{s,t} (see notes), and added the index condition s^- <= t^-.

Consequently, the persistent up-Laplacian is completely determined by the new persistent chain spaces. To compute $\Delta_{q,+}^{s,t}$, we gather bases of all spaces $N_{q+1}^{s^-,t^-}$ with $s^-\leq s$ and $t^-\leq t$, and orthogonalize them, since these spaces are independent, but not mutually orthogonal. If the columns of $Q$ form the resulting orthonormal basis of $C_{q+1}^{s,t}$, then
\[
    [\Delta_{q,+}^{s,t}]=B_s\,QQ^T B_s^T,
\]
where $B_s$ consists of the rows of $[\partial_{q+1}^t]$ corresponding to the $q$-simplices of $K_s$. % NEW, CHECK

Moreover, the persistent up-Laplacian only needs to be reevaluated at a few pairs. If $C_{q+1}^{s,t}=C_{q+1}^{s,t-1}$, then $\Delta_{q,+}^{s,t}=\Delta_{q,+}^{s,t-1}$. If $C_{q+1}^{s,t}=C_{q+1}^{s-1,t}$, then $\Delta_{q,+}^{s,t}=\iota\,\Delta_{q,+}^{s-1,t}\,\iota^*$ by Theorem~\ref{thm:vertical_operator}, which has the same non-zero eigenvalues as $\Delta_{q,+}^{s-1,t}$. By Theorem~\ref{thm:new_chain_decomposition}, $C_{q+1}^{s,t}$ differs from both $C_{q+1}^{s,t-1}$ and $C_{q+1}^{s-1,t}$ in only two cases:
\begin{enumerate}
    \item $N_{q+1}^{s,t}\neq\{0\}$. New chains are born at $(s,t)$, so $C_{q+1}^{s,t}$ changes, and the eigenvalues of the persistent up-Laplacian may change as well.
    \item There exist $s'<s$ and $t'<t$ with $s\leq t'$ such that $N_{q+1}^{s',t}$ and $N_{q+1}^{s,t'}$ are both non-zero. Otherwise, $C_{q+1}^{s,t}$ would equal $C_{q+1}^{s,t-1}$ or $C_{q+1}^{s-1,t}$, respectively.
\end{enumerate}
These cases can be combined by allowing $s'\leq s$ and $t'\leq t$, but we state them separately for clarity.

\begin{figure}
    \centering
    \includegraphics[width=\linewidth]{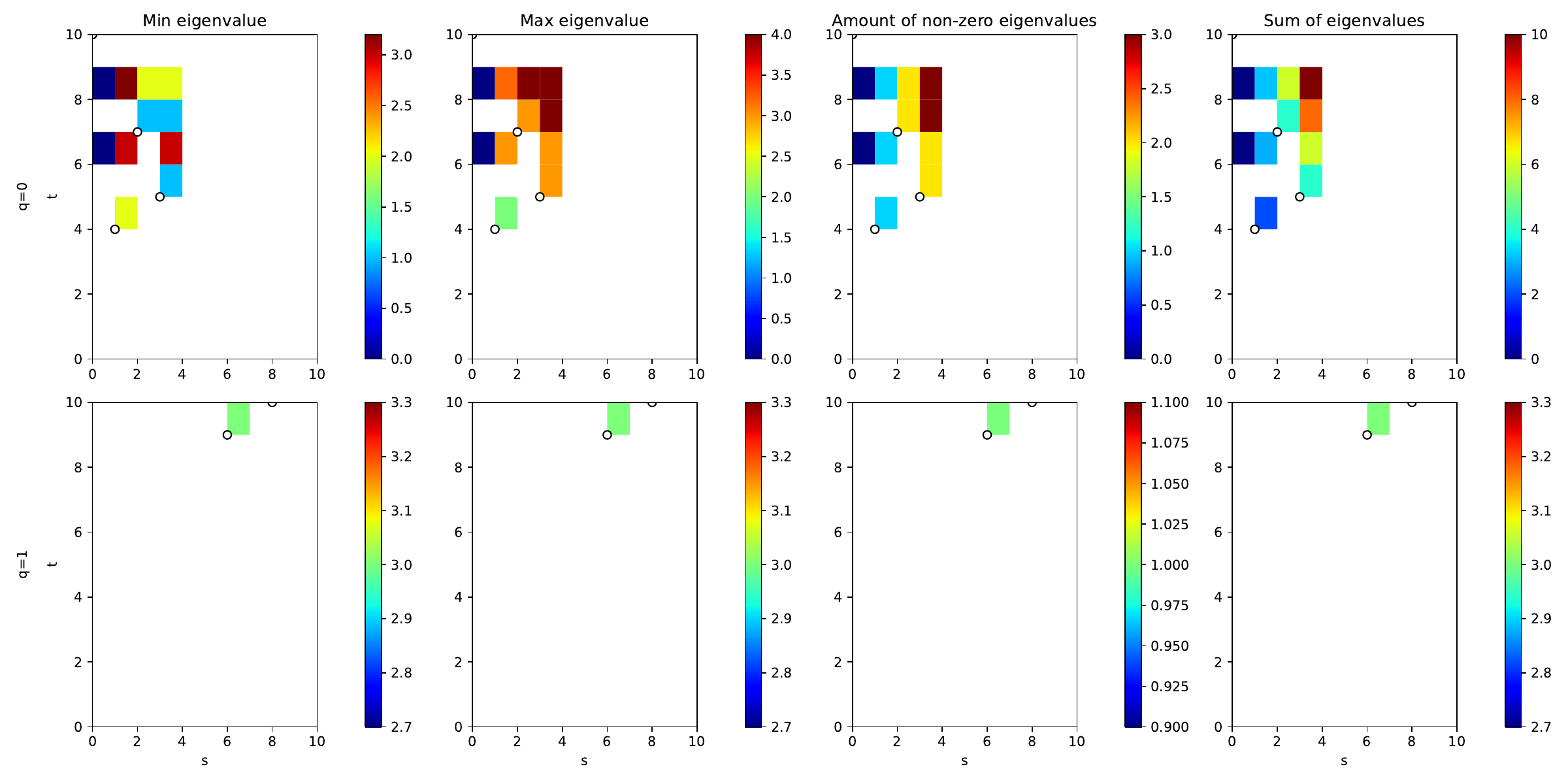}
    \caption{Aggregation functions applied to the eigenvalues of the persistent up-Laplacian $\Delta_{q,+}^{s,t}$, computed from the new persistent chain spaces, for the filtration in Figure~\ref{fig:filtration_example}. Only the squares at the pairs where $\Delta_{q,+}^{s,t}$ needs to be reevaluated are shown. White dots represent the bars of the persistence barcode, with the bar $[b,d)$ drawn at $(b,d)$.}
    \label{fig:reduced_lap_example}
\end{figure}

Figure~\ref{fig:reduced_lap_example} shows the eigenvalues of $\Delta_{q,+}^{s,t}$, computed from the new persistent chain spaces, only at the pairs that satisfy one of these two conditions. This figure contains exactly the same information as the up-Laplacian information in Figure~\ref{fig:pers_lap_evals_example_filtration}: every missing square can be recovered from a neighboring square, namely from the square below it if $C_{q+1}^{s,t}=C_{q+1}^{s,t-1}$, and from the square to its left, up to additional zero eigenvalues, if $C_{q+1}^{s,t}=C_{q+1}^{s-1,t}$. 
% CHECK: if Figure~\ref{fig:pers_lap_evals_example_filtration} shows the full persistent Laplacian, the panels for q>=1
% also contain the down-Laplacian, so "exactly the same information" only holds for the up part.
The advantage is that we know which persistent Laplacians are relevant and should be used in further analysis. Moreover, this leads to a faster algorithm, as we show in Section~\ref{ch:algorithm}.

%% file: shortened_sections/algorithm_v2.tex
% Revised algorithm section (version without the lemma).
% TODO  = information that seems to be missing.
% CHECK = places where I had to guess the intended meaning, or changed the mathematics.
% NEW   = sentences I added; remove them if you prefer.
% The commented-out paragraphs have been removed. The paragraph on the down-Laplacian has been moved
% to after the description of the algorithm, since it refers to the evaluated pairs.

In this section, we present an algorithm, based on the results of Section~\ref{ch:new_chain_spaces}, that computes the persistent up-Laplacians of a filtration at all pairs where they need to be evaluated. It consists of two steps. First, the standard algorithm for persistent homology determines the new persistent chain spaces, and thereby the pairs at which the persistent up-Laplacian has to be evaluated. Second, the persistent up-Laplacian is computed at these pairs. 
We illustrate the algorithm with an example.

\subsection{Obtaining the new persistent chain spaces}
Throughout this section, $[\partial_{q+1}]$ denotes the matrix of the boundary operator $\partial_{q+1}$ of the final complex $K_T$, with the simplices ordered compatibly with the filtration as in Section~\ref{ch:mathematical_background}, and we write $n_q:=n_q^T$ and $n_{q+1}:=n_{q+1}^T$. For a column $R_j$ of a matrix $R$, let $\operatorname{low}(R_j)$ denote the row index of its lowest non-zero entry, with $\operatorname{low}(R_j):=0$ if $R_j=0$. A matrix $R$ is \emph{reduced} if $\operatorname{low}(R_i)\neq\operatorname{low}(R_j)$ for all $i\neq j$ with $R_i\neq 0$ and $R_j\neq 0$. The standard algorithm for persistent homology (Algorithm~\ref{alg:standard}) reduces $[\partial_{q+1}]$ by subtracting multiples of columns from columns to their right. This yields a reduced matrix $R=[\partial_{q+1}]V$, where $V$ is an invertible upper-triangular matrix, in $\mathcal{O}(n_qn_{q+1}^2)$ time. The columns of $V$ corresponding to the zero columns of $R$ form a basis of the kernel of $[\partial_{q+1}]$, that is, of the space of $(q+1)$-cycles of $K_T$. This holds more generally: if $R=MV$ is reduced for some matrix $M$ and an invertible matrix $V$, then the columns of $V$ corresponding to the zero columns of $R$ form a basis of the kernel of $M$, since the non-zero columns of a reduced matrix are linearly independent. %

For the persistent Laplacian, we instead need a basis of $C_{q+1}^{s,t}$. Since $C_{q+1}^{s,t}$ consists of the chains in $C_{q+1}^t$ whose boundary lies in $C_q^s$, it is the kernel of the submatrix $M^{s,t}$ of $[\partial_{q+1}]$ obtained by keeping only the first $n_{q+1}^t$ columns and removing the first $n_q^s$ rows. The reduction of the full boundary matrix also reduces this submatrix. Removing the last $n_{q+1}-n_{q+1}^t$ columns does not affect the reduction, since the standard algorithm only works from left to right. Since $V$ is upper-triangular, we obtain $M^{s,t}V^t=R^{s,t}$, where $V^t$ consists of the first $n_{q+1}^t$ rows and columns of $V$, which is again invertible and upper-triangular, and $R^{s,t}$ is obtained from $R$ by keeping the first $n_{q+1}^t$ columns and removing the first $n_q^s$ rows. Removing these rows keeps the matrix reduced: a column with $\operatorname{low}(R_j)>n_q^s$ keeps its lowest non-zero entry, and a column with $\operatorname{low}(R_j)\leq n_q^s$ becomes zero. Hence, increasing $s$ creates new zero columns, which indicate new features. By the observation above, the columns $V_j$ with $j\leq n_{q+1}^t$ and $\operatorname{low}(R_j)\leq n_q^s$ form a basis of $C_{q+1}^{s,t}$.
% CHECK: this is the step that justifies the orthogonalization below (see notes).
% TODO (optional): cite a standard reference for the fact that left-to-right column operations also reduce
% the lower-left submatrices of the boundary matrix (e.g. the proof of the Pairing Lemma in Edelsbrunner and Harer, Computational Topology).

For each column $j$, we define the start and end time
\[
    s_j:=\min\{s:n_q^s\geq\operatorname{low}(R_j)\},
    \qquad
    t_j:=\min\{t:n_{q+1}^t\geq j\}.
\]
The end time $t_j$ is the time at which the $(q+1)$-simplex $j$ appears, and for $R_j\neq 0$, the start time $s_j$ is the time at which the $q$-simplex $\operatorname{low}(R_j)$ appears. Zero columns of $R$ have $s_j=0$, in agreement with the convention $C_{q+1}^{-1,t}=\{0\}$. Since $s_j\leq s$ if and only if $\operatorname{low}(R_j)\leq n_q^s$, and $t_j\leq t$ if and only if $j\leq n_{q+1}^t$, the columns $V_j$ with $s_j\leq s$ and $t_j\leq t$ form a basis of $C_{q+1}^{s,t}$.

Consequently, the new persistent chain spaces can be read off from the reduction. Let $W^{s,t}:=\{V_j:(s_j,t_j)=(s,t)\}$. Applying the previous statement to $(s-1,t)$ and $(s,t-1)$ shows that the columns $V_j$ with $s_j\leq s$, $t_j\leq t$ and $(s_j,t_j)\neq(s,t)$ form a basis of $C_{q+1}^{s-1,t}+C_{q+1}^{s,t-1}$, and the vectors in $W^{s,t}$ are linearly independent of them, since $V$ is invertible. Orthogonalizing the vectors in $W^{s,t}$ with respect to this basis, for example with the Gram--Schmidt process, therefore yields a basis of $N_{q+1}^{s,t}$. In particular, $\dim N_{q+1}^{s,t}=|W^{s,t}|$, so $N_{q+1}^{s,t}\neq\{0\}$ only at the pairs $(s_j,t_j)$, which are determined by the persistence pairing.

\begin{example}\label{ex:algorithm}
    Consider the persistent up-Laplacian in dimension $0$ for the filtration of Figure~\ref{fig:filtration_example}. Here, vertex $i$ appears at time $i$ for $i=0,\dots,3$, and the five edges appear at the times $4,\dots,8$.
    Ordering the vertices and edges by their appearance, we obtain the boundary matrix
    \[[\partial_{1}]=
    \begin{bmatrix}
        -1 & -1 &  0 &  0 & -1\\
         1 &  0 & -1 &  0 &  0\\
         0 &  0 &  0 & -1 &  1\\
         0 &  1 &  1 &  1 &  0
    \end{bmatrix}.
    \]
    The standard algorithm yields
    \[
    R = \begin{bmatrix}
        -1 & -1 &  0 &  1 &  0\\
         1 &  0 &  0 &  0 &  0\\
         0 &  0 &  0 & -1 &  0\\
         0 &  1 &  0 &  0 &  0
    \end{bmatrix}, \qquad V = \begin{bmatrix}
         1 &  0 &  1 &  0 &  0\\
         0 &  1 & -1 & -1 & -1\\
         0 &  0 &  1 &  0 &  0\\
         0 &  0 &  0 &  1 &  1\\
         0 &  0 &  0 &  0 &  1
    \end{bmatrix}.
    \]
    Columns $3$ and $5$ of $R$ are zero; they correspond to the $1$-cycles born at times $6$ and $8$, and have start time $s_3=s_5=0$. For column $3$, the end time is $t_3=\min\{t:n_1^t\geq 3\}=6$. No other column $j$ has $s_j\leq 0$ and $t_j\leq 6$, so
    \[
    N_{1}^{0,6}=\operatorname{span}\left\{\frac{1}{\sqrt{3}}\begin{bmatrix}
        1 \\ -1 \\ 1 \\ 0 \\ 0
    \end{bmatrix}\right\}.
    \]
    Column $5$ has $(s_5,t_5)=(0,8)$. Since $(s_3,t_3)=(0,6)$ satisfies $s_3\leq 0$ and $t_3\leq 8$, we orthogonalize $V_5$ with respect to $V_3$, which gives
    \[
    N_1^{0,8}=\operatorname{span}\left\{\frac{1}{\sqrt{24}}\begin{bmatrix}
        -1 \\ -2 \\ -1 \\ 3 \\ 3
    \end{bmatrix}\right\}.
    \]
    % CHECK: the original had +2 as second entry; that vector is neither orthogonal to V_3 nor a cycle.
    Finally, the first column has $\operatorname{low}(R_1)=2$, so its start time is the first $s$ with $n_0^s\geq 2$, that is, $s_1=1$, and its end time is the first $t$ with $n_1^t\geq 1$, that is, $t_1=4$. No other column $j$ has $s_j\leq 1$ and $t_j\leq 4$, hence
    \[
    N_{1}^{1,4}=\operatorname{span}\left\{\begin{bmatrix}
        1 \\ 0 \\ 0 \\ 0 \\ 0
    \end{bmatrix}\right\}.
    \]
    % TODO (optional): also give N_1^{3,5} and N_1^{2,7} for completeness.
\end{example}

\begin{remark}\label{rem:clearing}
For standard persistent homology, the standard algorithm can be accelerated with clearing~\cite{chen_persistent_2011,bauer_clear_2014}. In this approach, the boundary matrices are reduced in decreasing dimension: every non-zero column $R'_i$ of the reduced matrix of $[\partial_{q+2}]$ with $\operatorname{low}(R'_i)=j$ shows that column $j$ of $[\partial_{q+1}]$ reduces to zero, so this column can be skipped. The same technique can be applied to the reduction step of our method. In that case, the column $V_j$ is replaced by $R'_i$, which is a $(q+1)$-cycle with lowest non-zero entry $j$, so that $V$ remains invertible and upper-triangular and the arguments of this section still apply.
\end{remark}

\subsection{Computing the persistent Laplacians}
In practice, we do not orthogonalize the new persistent chain spaces with respect to each other when they are found. As shown above, the vectors $V_j$ with $s_j\leq s$ and $t_j\leq t$ form a basis of $C_{q+1}^{s,t}$, regardless of whether they are orthogonal. In particular, they are linearly independent, which reflects the direct sum in Theorem~\ref{thm:new_chain_decomposition}. Hence, the Gram--Schmidt process can be applied to them directly: every step yields a non-zero vector, and no dependent vectors have to be detected or discarded. Since the spaces $N_{q+1}^{s,t}$ are not mutually orthogonal anyway, these vectors have to be orthonormalized whenever a persistent up-Laplacian is computed. It therefore suffices to store the columns $V_j$ and to orthonormalize them only when a persistent up-Laplacian is evaluated. This basic version of the method is given in Algorithm~\ref{alg:new_chains}.

This step can be optimized further, since many computations are repeated at different pairs. When $\Delta_{q,+}^{s,t}$ is evaluated and a pair $(s^-,t^-)$ with $s^-\leq s$ and $t^-\leq t$ has already been evaluated, we can reuse its orthonormal basis and only add the new vectors. Adding $m$ vectors to an orthonormal basis of $k$ vectors in $\mathbb{R}^{n_{q+1}}$ with the Gram--Schmidt process takes $\mathcal{O}\big(n_{q+1}m(k+m)\big)$ time, which is $\mathcal{O}(n_{q+1}mk)$ if $m\leq k$. Since there are $n_{q+1}$ columns in total and the added vectors are always linearly independent of the existing basis, we have $m+k\leq n_{q+1}$, so a single update takes $\mathcal{O}\big(n_{q+1}k(n_{q+1}-k)\big)$ time.

Finally, at each evaluated pair we compute $[\Delta_{q,+}^{s,t}]=B_sQQ^TB_s^T$ as in Section~\ref{ch:new_chain_spaces}, where the columns of $Q$ form an orthonormal basis of $C_{q+1}^{s,t}$ and $B_s$ consists of the first $n_q^s$ rows of $[\partial_{q+1}]$. This is done most efficiently by first computing $U:=B_sQ$ and then $[\Delta_{q,+}^{s,t}]=UU^T$. Previously computed results can be reused here as well. If the orthonormal basis of an evaluated pair $(s^-,t^-)$ is extended by $m$ new vectors, which form the columns of $Q'$, then
\[
    [\Delta_{q,+}^{s,t}]=\big[\iota\,\Delta_{q,+}^{s^-,t^-}\iota^*\big]+U'U'^T,
    \qquad U':=B_sQ',
\]
where $\big[\iota\,\Delta_{q,+}^{s^-,t^-}\iota^*\big]$ is obtained from $[\Delta_{q,+}^{s^-,t^-}]$ by adding zero rows and columns for the $q$-simplices of $K_s\setminus K_{s^-}$. % CHECK: zero padding added
Computing $U'$ takes $\mathcal{O}(n_qn_{q+1}m)$ time, and computing $U'U'^T$ takes $\mathcal{O}(n_q^2m)$ time. Computing the eigenvalues of the resulting symmetric matrix takes $\mathcal{O}(n_q^3)$ time.

The persistent up-Laplacians at all other pairs follow from the evaluated ones. For a pair $(s,t)$, let $(s',t')$ be the evaluated pair with $s'\leq s$ and $t'\leq t$ for which $t'$ is maximal and, given $t'$, $s'$ is maximal. Then $\Delta_{q,+}^{s,t}=\iota\,\Delta_{q,+}^{s',t'}\iota^*$, so both have the same non-zero eigenvalues. If no such pair exists, then $\Delta_{q,+}^{s,t}=0$. 

In particular, the algorithm also yields the down-Laplacian. As discussed in Section~\ref{sec:effect_separation}, the non-zero eigenvalues of $\Delta_{q,-}^{K_s}$ are those of $\Delta_{q-1,+}^{K_s}=\Delta_{q-1,+}^{s,s}$. Given the persistent up-Laplacians in dimension $q-1$ at the evaluated pairs, either $\Delta_{q-1,+}^{s,s}$ has been evaluated, or its non-zero eigenvalues are those of $\Delta_{q-1,+}^{s',s}$, with $s'<s$ maximal among the evaluated pairs. If no such pair exists, we decrease the end time by one and repeat, until we find an evaluated pair or the end time reaches $0$. If no evaluated pair is found, the down-Laplacian is zero.

\begin{algorithm}
\caption{The standard algorithm.}\label{alg:standard}
    \KwData{Boundary matrix $[\partial_{q+1}]\in\mathbb{R}^{n_q\times n_{q+1}}$}
    \KwResult{Reduced matrix $R\in\mathbb{R}^{n_q\times n_{q+1}}$ and invertible upper-triangular matrix $V\in\mathbb{R}^{n_{q+1}\times n_{q+1}}$ such that $R=[\partial_{q+1}]V$}
    $R\gets[\partial_{q+1}]$\;
    $V\gets I_{n_{q+1}}$ (the identity matrix)\;
    \For{$j\gets 1$ \KwTo $n_{q+1}$}{
        \While{there exists $i<j$ such that $0\neq\operatorname{low}(R_i)=\operatorname{low}(R_j)$}{
            $k\gets\operatorname{low}(R_j)$\;
            $\alpha\gets R(k,j)/R(k,i)$\;
            subtract $\alpha$ times column $i$ of $V$ from column $j$ of $V$\;
            subtract $\alpha$ times column $i$ of $R$ from column $j$ of $R$\;
        }
    }
\end{algorithm}

\begin{algorithm}
\caption{Computation of the persistent up-Laplacians at all pairs where they need to be evaluated.}\label{alg:new_chains}
    \KwData{Boundary matrix $[\partial_{q+1}]$}
    \KwResult{$[\Delta_{q,+}^{s,t}]$ for all pairs $(s,t)$ where it needs to be evaluated}
    $R,V\gets\operatorname{Standard\_Algorithm}([\partial_{q+1}])$\;
    all sets $W^{s,t}$ are initially empty\;
    \For{$j\gets 1$ \KwTo $n_{q+1}$}{
        $s_j\gets\min\{s:n_q^s\geq\operatorname{low}(R_j)\}$\;
        $t_j\gets\min\{t:n_{q+1}^t\geq j\}$\;
        add $V_j$ to $W^{s_j,t_j}$\;
    }
    \For{$t\gets 0$ \KwTo $T$}{
        \For{$s\gets 0$ \KwTo $t$}{
            \If{there exist $s^-\leq s$ and $t^-\leq t$ such that $W^{s^-,t}$ and $W^{s,t^-}$ are non-empty}{
                $Q\gets\operatorname{Orthonormalize\_Columns}\big(\bigcup_{s^-\leq s,\,t^-\leq t}W^{s^-,t^-}\big)$\;
                $B_s\gets[\partial_{q+1}]([n_q^s],:)$\;
                $[\Delta_{q,+}^{s,t}]\gets B_sQQ^TB_s^T$\;
            }
        }
    }
\end{algorithm}

%% file: shortened_sections/comparison_v2.tex
% Revised benchmark section.
% TODO  = information that seems to be missing.
% CHECK = places where I had to guess the intended meaning.

In this section, we compare our method with two existing methods for computing persistent Laplacians: the Schur complement method~\cite{laplacians_algorithms} and PETLS~\cite{computationalefficiency}. For each method, we measure the time needed to compute all persistent Laplacians up to dimension 3 of a filtration built on randomly sampled points, for several types of point clouds. We measure the computation time of the matrices and their spectra; the construction of the filtration is excluded from the timing.
All experiments were run in WSL, started from the Windows console with no other applications running. The machine has an AMD Ryzen~7 5700X3D CPU and 32\,GB of RAM, of which 24\,GB is available to WSL. None of the tested methods uses the GPU.

We consider four types of point clouds. First, to study low-dimensional simplicial complexes, we sample points from the unit square in $\mathbb{R}^2$ and build the alpha filtration. This filtration contains simplices of dimension at most two, so the persistent Laplacian is fully described by $\Delta_0$ and $\Delta_1$, dimension 2 only contains the down-Laplacian and therefore no new information. Second, we sample points from the unit cube in $\mathbb{R}^3$ and build the Rips filtration up to radius $r = 0.3$. Third, to obtain a point cloud with more structure, we sample points on the unit sphere $S^2 \subset \mathbb{R}^3$ and build the alpha filtration. These points are generated by sampling from the standard normal distribution in $\mathbb{R}^3$ and then rescaling each point to have norm one.
Finally, to simulate high-dimensional data, we sample points from the unit cube in $\mathbb{R}^{100}$ and build the Rips filtration up to radius $r = 5$. This filtration contains many high-dimensional simplices, but we restrict our analysis to dimension at most four.

Our method consists of two steps: we first determine the locations at which the persistent Laplacian needs to be computed, and then compute it efficiently at these locations. In the figures, our method is labeled \emph{New Chains}. Obtaining the same information with PETLS alone requires computing the persistent Laplacian for every pair of start and end values of the filtration; we refer to this as PETLS. The two approaches can also be combined: we use the first step of our method to find the relevant locations, and then compute the persistent Laplacian only at these locations using PETLS. We include this combination, referred to as PETLS (Relevant Locations), in the comparison as well. Furthermore, PETLS contains a method for alpha complexes using the fact that they are non-branching complexes \cite{dong_faster_2024, botnan_nonbranching_2026}, which has been used for the 3D unit sphere point cloud.

\begin{figure}
    \centering
    \includegraphics[width=\linewidth]{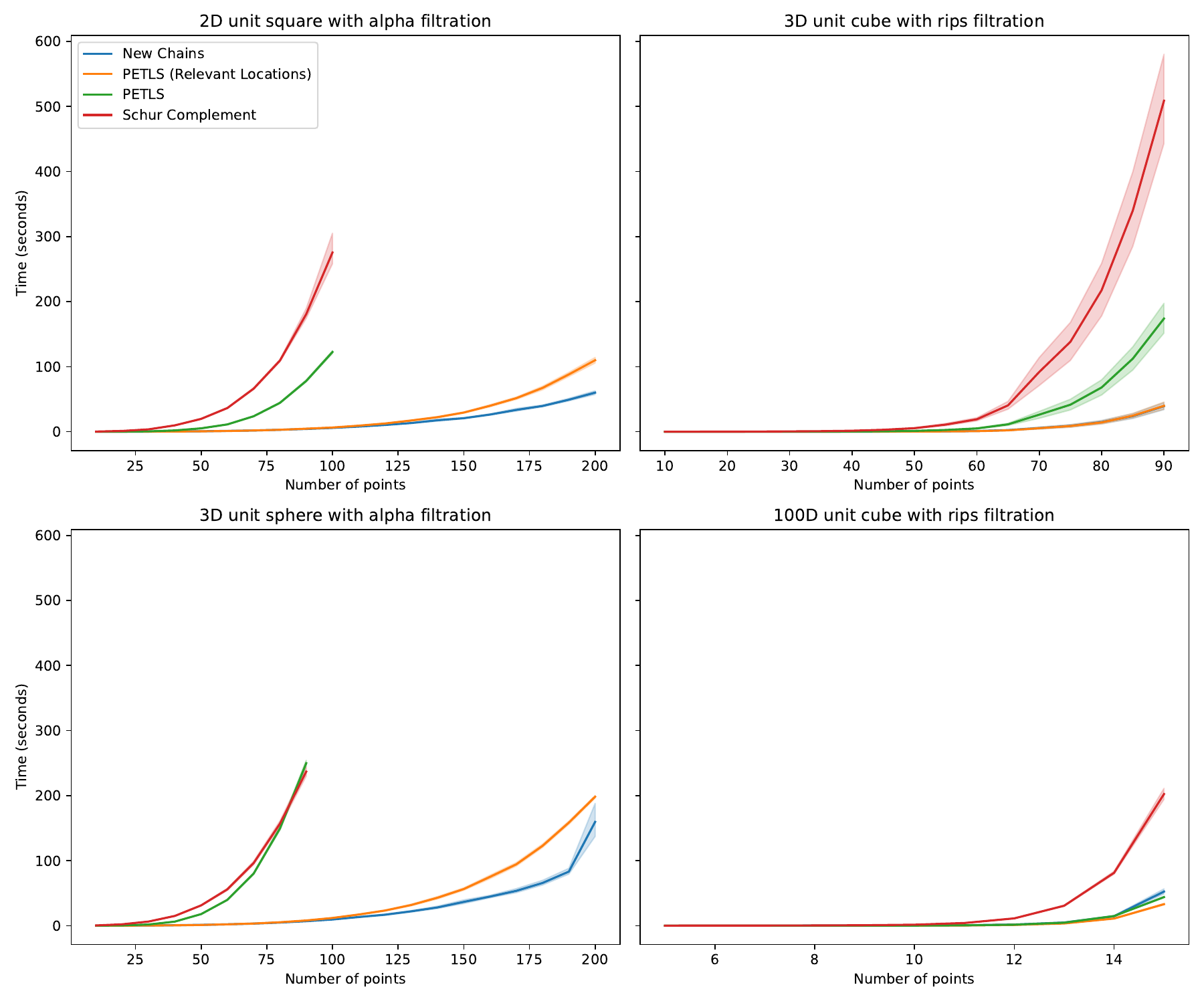}
    \caption{Time needed to compute all persistent Laplacians of a filtration built on randomly sampled points. For each number of points, five seeds are used and each seed is run five times; we show the mean together with a 95\% confidence interval.
    (Top left) Unit square in $\mathbb{R}^2$ with the alpha filtration, for $n = 10, 20, \dots, 200$ points; PETLS and the Schur complement method are only run up to $n = 100$.
    (Top right) Unit cube in $\mathbb{R}^3$ with the Rips filtration up to $r = 0.3$, in steps of 5 points.
    (Bottom left) Unit sphere with the alpha filtration, in steps of 10 points.
    (Bottom right) Unit cube in $\mathbb{R}^{100}$ with the Rips filtration up to $r = 5$, for $n = 5, 6, \dots, 15$ points.
    Our method (New Chains) is the fastest in most settings, in particular on the alpha filtrations. In the high-dimensional case, PETLS (relevant locations) is the fastest.}
    \label{fig:timings_compared}
\end{figure}

Each of the point clouds is randomly sampled with five distinct seeds and each method is run five times on each seed. The results are shown in Figure~\ref{fig:timings_compared}. Our method performs well in most settings. On both alpha filtrations it is faster than all other methods, and on the three-dimensional Rips filtration its running time is nearly identical to that of PETLS (Relevant Locations). In the high-dimensional case, however, our method is slower than both PETLS variants. We expect this to be caused by the higher-dimensional simplices, whose associated matrices are larger but also sparser: PETLS uses algorithms that exploit sparsity, whereas our implementation does not distinguish between sparse and dense matrices.

As expected, the Schur complement method is the slowest in all cases except the unit sphere, where regular PETLS takes the longest for larger amounts of points. The difference on this point set is substantial, and we do not yet have an explanation for it. Finally, PETLS is almost always faster when it is evaluated only at the relevant locations, as one would expect. For small numbers of points it can occasionally be slower, since the first step of our method then has to be run as well, although this is hard to see in the figure.

\begin{figure}
    \centering
    \includegraphics[width=\linewidth]{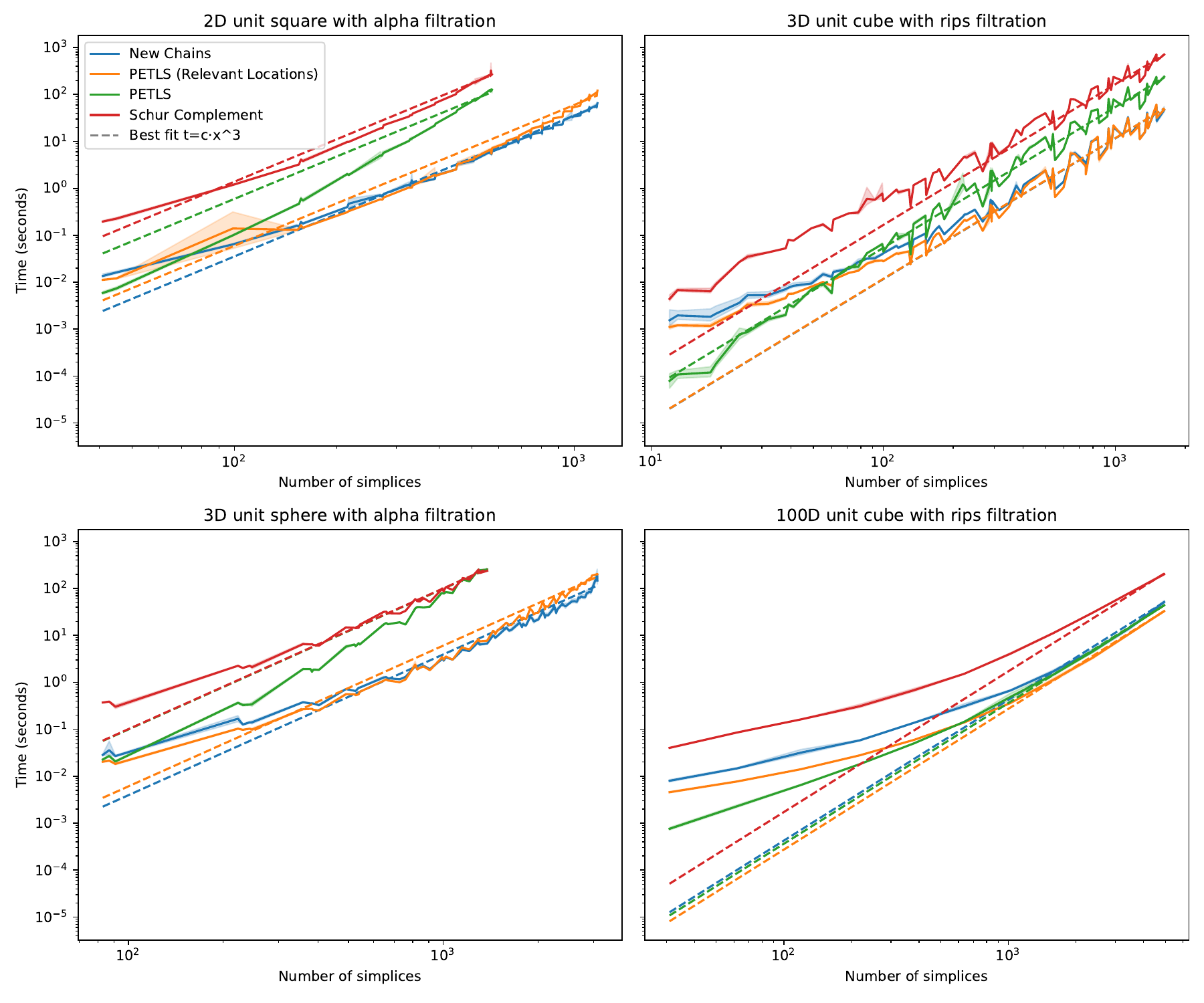}
    \caption{The timings of Figure~\ref{fig:timings_compared}, plotted against the total number of simplices~$N$ in the filtration, with both axes on a logarithmic scale. The point sets, filtrations and numbers of points are the same as in Figure~\ref{fig:timings_compared}: (top left) unit square in $\mathbb{R}^2$ with the alpha filtration, (top right) unit cube in $\mathbb{R}^3$ with the Rips filtration up to $r = 0.3$, (bottom left) unit sphere with the alpha filtration, and (bottom right) unit cube in $\mathbb{R}^{100}$ with the Rips filtration up to $r = 5$. Each configuration is run five times per seed. If multiple seeds give the same value of $N$, their runs are combined into one mean with a 95\% confidence interval; otherwise they are shown as separate points. The dashed lines show the least-squares fit of $t = c\,N^3$ for each method, in the same color.
    Methods with polynomial running time should approach a straight line in a log-log plot. The cubic fits describe the running times well, except for PETLS on the alpha filtrations, whose slope is closer to that of a quartic dependence.}
    \label{fig:timings_n_simplices}
\end{figure}

The number of points is only a rough measure of the size of the problem. As Figure~\ref{fig:timings_compared} shows, the computation time also depends strongly on the type of filtration and on the dimension of the data. A more informative measure is the total number of simplices~$N$ in the filtration. In Figure~\ref{fig:timings_n_simplices}, we therefore show the same results with $N$ on the horizontal axis. Different seeds may give slightly different values of $N$ for the same number of points. Runs with the same value of $N$ are combined into a single mean and confidence interval, while runs with different values of $N$ are shown as separate points. As a result, the curves may jump up and down between seeds.

Both axes are logarithmically scaled, so that small and large problems are visible in the same plot. Moreover, in such a plot a power law $t = c\,N^k$ appears as a straight line with slope $k$, and a polynomial of degree $k$ approaches such a line as $N$ grows. To compare each method with cubic growth, we fit the model $t = c\,N^3$ by least squares. In the log-log plot, this amounts to fitting only the intercept $\log c$, since the slope is fixed at three. The fitted value of $c$ differs between point sets, which indicates that the running time depends not only on the number of simplices, but also on other properties of the filtration.

Some methods only follow a straight line once $N$ is sufficiently large, but in almost all cases the timings eventually follow the cubic fit. %Based on the analysis in the previous section, we expected the running time to grow as $N^4$, but for most methods we do not observe this. The tested cases appear to be closer to the second extremal case discussed in the previous section, which gives a cubic relation. 
A quartic relation, which corresponds to the worst case discussed in Section~\ref{ch:conclusion}, is only observed for PETLS on the alpha filtrations.

Note that our experiments may not reach values of $N$ large enough to clearly reveal super-cubic behavior. However, computation times are already long at these sizes, so the regime in which a higher-order term would dominate may lie beyond any practical time budget. This makes the distinction less relevant in practice.

%% file: shortened_sections/conclusion_v3.tex
% The field of topological data analysis has received a lot of attention in recent years. In particular, persistent homology has been studied extensively and has proven useful in many applications. However, persistent homology only captures topological features, which are invariant under deformations such as stretching, twisting and bending, and this may not be sufficient for every application. For example, the persistent Laplacian, which also captures geometric features, has been shown to improve results in handwritten digit recognition~\cite{pers_lap_mnist} over using standard persistent homology.

In this paper, we studied how the persistent Laplacians of an entire filtration can be obtained without evaluating them independently at every possible location, which is expensive and largely redundant. In Section~\ref{sec:duplicate_information}, we showed that the persistent Laplacians at different locations largely contain the same information.
% Computing the persistent Laplacian at every possible location is, however, expensive. In Section~\ref{sec:duplicate_information}, we showed that the persistent Laplacians at different locations largely contain the same information.
We introduced a new visualization technique, which shows that only certain locations influence all subsequent persistent Laplacians, and we reached the same conclusion using the structure of the boundary matrix.

The influence of these locations on other persistent Laplacians can be split into two effects: a horizontal effect, which influences locations with later start times, and a vertical effect, which influences locations with later end times. As shown by Theorems~\ref{thm:horizontal_operator} and~\ref{thm:vertical_operator}, these effects can be removed independently of each other.
This leads to a decomposition of the chain space $C_{q+1}^{s,t}$ into a horizontal, a vertical, a diagonal and a new component. We showed that at every location, only this new component, which we call the new persistent chain space, is needed to obtain all persistent Laplacians, see Theorem \ref{thm:new_chain_decomposition}.

The locations at which the new persistent chain space is non-trivial can be found with the standard algorithm for persistent homology, as shown in Section~\ref{ch:algorithm}, and a basis for each new persistent chain space can then be obtained with the Gram--Schmidt process. To reconstruct the persistent Laplacian at an arbitrary location, however, all new persistent chain spaces that appear before this location have to be found and orthogonalized with respect to each other. 
When the persistent Laplacian is computed at all relevant locations of a filtration, the orthonormal bases can be built incrementally, so that previously computed bases are reused (Section~\ref{ch:algorithm}). 
% A crude bound on the cost follows directly from the construction: the persistent up-Laplacian has to be evaluated at no more than $n_{q+1}^2$ pairs, since at every such pair $s$ is a start value $s_j$ and $t$ is an end value $t_j$. At each of them, the orthonormalization, the computation of the matrix and its eigendecomposition take $\mathcal{O}(n_{q+1}^3+n_qn_{q+1}^2+n_q^2n_{q+1}+n_q^3)$ time.%When the persistent Laplacian is computed at all locations of a filtration, this can be organized efficiently, which reduces the time complexity compared to computing each orthogonalization separately.
%We did not derive the time complexity of this algorithm in general, but discussed two extremal cases in which it is $\mathcal{O}(N^4)$ and $\mathcal{O}(N^3)$, respectively, where $N$ is the number of simplices. Worst case in a filtration with a single simplex added at each step, this is not faster than the Schur complement method~\cite{pers_lap_mnist}. However, our method is less reliant on the amount of timesteps.
% and will also work on filtration with multiple simplices added in each step.

In Section~\ref{ch:comparison}, we compared the running time of our method with that of existing methods. On filtrations with low-dimensional simplices, our method is faster than or on par with all other methods. Moreover, the relevant locations found by our method can also be used to speed up existing methods: PETLS~\cite{computationalefficiency} evaluated only at these locations comes close to our method, and is the fastest method in the high-dimensional case. There, our method is slower, which we attribute to our implementation not exploiting the sparsity of the matrices.

Besides the algorithm, we believe that the relevant locations themselves are an important result of this paper. Many previous works evaluate the persistent Laplacian only on a preselected grid of locations \cite{pers_lap_mnist}.
% TODO: add citations for previous works that use a fixed set of locations.
Such a grid has to be chosen before computation, so the locations themselves carry no information. In contrast, the relevant locations we introduce already encode information about the filtration. Moreover, persistent Laplacians that only ``last'' for a short time are captured as well, without having to choose a sufficiently fine grid. This is one of the strengths of standard persistent homology, and our approach makes it available for persistent Laplacians as well.

\subsection*{Future research}
\begin{figure}
   \centering
    \includegraphics[width=0.5\linewidth]{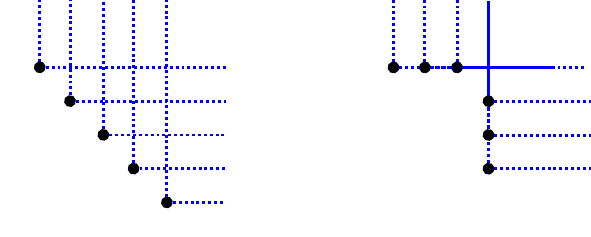}
    \caption{Two extreme configurations of the pairs with non-trivial new persistent chain spaces. Left: their regions of influence overlap as often as possible, but at each evaluated pair only one vector is added. Right: they overlap less often, but at one pair many vectors have to be combined}
    \label{fig:extremal_cases}
\end{figure}

The crude bound previously discussed does not take into account that previously computed orthonormal bases are reused. Two extreme configurations of the relevant pairs suggest a better bound (Figure \ref{fig:extremal_cases}).
In the first, the regions of influence of the new persistent chain spaces overlap as often as possible, so that the persistent up-Laplacian has to be evaluated at $\mathcal{O}(n_{q+1}^2)$ pairs, but at each of them only a single vector is added to a previously computed basis; the orthogonalization then takes $\mathcal{O}(n_{q+1}^4)$ time in total. In the second, there are only $n_{q+1}+1$ relevant pairs, but at one of them $n_{q+1}/2$ vectors have to be added at once; the orthogonalization then takes $\mathcal{O}(n_{q+1}^3)$ time. We therefore conjecture that the total cost, excluding the eigendecompositions, is $\mathcal{O}(n_{q+1}^4+n_qn_{q+1}^3+n_q^2n_{q+1}^2)$. For a filtration in which a single simplex is added at each step, this is of the same order as the cost of the Schur complement approach over all pairs~\cite{laplacians_algorithms}. In our experiments (Section~\ref{ch:comparison}), the running times grow approximately cubically in the number of simplices, which suggests that typical filtrations are closer to the second configuration. Proving the conjectured bound, or a sharper one in terms of the persistence pairing, is left for future work.

PETLS uses several techniques to speed up its computations, such as sparse matrix multiplication. Our implementation does not use these, and it is slower in settings where the matrices are sparse; we believe this is the main reason. Using sparse matrix computations in our implementation could therefore yield substantial speedups. In addition, our implementation is written in Python, and rewriting it in a compiled language such as C++ may further improve its performance.

The standard algorithm itself can also be accelerated, for example using clearing~\cite{chen_persistent_2011, bauer_clear_2014}.
In our experiments, the standard algorithm accounted for only a small part of the total running time, but implementing clearing may still be worthwhile, especially in combination with the approach described below.

Finally, although our method yields the full spectrum of every persistent Laplacian in the filtration, computing all of these spectra remains time-consuming. However, all persistent Laplacians are determined by a small set of vectors that is available once the standard algorithm part of our method has been run. For deep learning applications, it may therefore suffice to compute only these vectors and let the model learn the relevant information from them, instead of computing every persistent Laplacian explicitly. This would greatly reduce the computation time, since only the standard algorithm needs to be run, and it can be accelerated further with clearing.

%% file: shortened_sections/supplementary_v3.tex
% Appendix with light text edits. The mathematics is unchanged.
% CHECK: = small mistake or notation issue, explained in the comment.
% MAJOR: = a problem in the argument itself, explained in more detail in the accompanying notes.

\section{Properties of the Moore--Penrose pseudo inverse} \label{sec:preliminaries}
\begin{definition}[Moore--Penrose pseudo inverse, \cite{pinv_properties}]\label{def:pinv}
    The Moore--Penrose pseudo inverse of a matrix $A\in\RR^{m\times n}$ is the $n\times m$ matrix $A^\dagger$, which satisfies:
    \begin{align}
        AA^\dagger A&=A \label{eq:pinv_def_A}\\
        A^\dagger AA^\dagger &= A^\dagger \label{eq:pinv_def_A_inv}\\
        (AA^\dagger)^T&=AA^\dagger\\
        (A^\dagger A)^T&=A^\dagger A.\label{eq:Q_symmetric}
    \end{align}
\end{definition}

We require three properties of the Moore--Penrose pseudo inverse. They are formulated in the following theorem.

\begin{theorem}[\cite{pinv_properties}]
    For a matrix $A\in\mathbb{R}^{m\times n}$, with Moore-Penrose pseudo inverse $A^\dagger$, we have,
    \begin{align}
        A^\dagger=A^T(AA^T)^\dagger&=(A^TA)^\dagger A^T \label{eq:pinv_A_to_pinv_AAT}\\
        (AA^T)^\dagger&=(A^T)^\dagger A^\dagger \label{eq:pinv_AAT}\\
        \text{For }a\in\mathbb{Z}^+, \text{let }B\in\mathbb{R}^{a\times n} \text{, then } \ker A\subseteq\ker B &\iff B=BA^\dagger A. \label{eq:pinv_subset_ker}
    \end{align}
\end{theorem}

Furthermore, the matrix $Q=A^\dagger A$ has some interesting properties. $Q$ is a projection matrix as, by Equation \eqref{eq:Q_symmetric}, it is symmetric and using Equation \eqref{eq:pinv_def_A_inv}, one can show it is idempotent. Indeed, 
\begin{equation}\label{eq:Q_idempotent}
    Q^2=A^\dagger AA^\dagger A=A^\dagger A=Q.
\end{equation}
It is known that this matrix is the projection matrix that projects onto the orthogonal complement of the kernel of $A$ \cite{projection_matrices}.

\section{Proof of Theorems \ref{thm:horizontal_operator} and \ref{thm:vertical_operator}}\label{sec:proof_horizontal_vertical}
To prove these theorems, we use the following lemmas, whose proofs are given in Appendix~\ref{ch:lemma_proofs}. For all proofs in the appendix, we use for matrices $A\in\mathbb{R}^{0\times n}$ and $B\in\mathbb{R}^{m\times 0}$ with $m,n\in\mathbb{Z}^+$, the convention that $BA=0_{m\times n}$, the $m\times n$ zero matrix, and when $m=n$, $AB=0$. Besides this, we use $\ker A=\mathbb{R}^n$ and $\ker B=\{0\}$. Furthermore, we use $A^\dagger A=A^TA=0_{n\times n}$ and $B^\dagger B=B^TB=0$. Finally, in general we define $C\in\mathbb{R}^{0\times 0}$ to be $C=0$.
\begin{lemma}\label{lemma:up_laplacian_boundary}
    Let $K$ and $L$ be simplicial complexes with $K\subseteq L$. Then the matrix representation of the persistent up-Laplacian is given by
    \begin{equation*}\label{eq:up_laplacian_boundary}
        [\Delta_{q,+}^{K,L}]=[\Delta_{q,+}^K]+B_{q,1}^{K,L}(I-(B_{q,2}^{K,L})^\dagger B_{q,2}^{K,L})(B_{q,1}^{K,L})^T,
    \end{equation*}
    where $B_{q,1}^{K,L}:=[\partial^L_{q+1}]([n_{q}^K],I_{q+1}^{K,L})$, 
    $B_{q,2}^{K,L}:=[\partial^L_{q+1}](I_{q}^{K,L},I_{q+1}^{K,L})$, 
    and $^\dagger$ denotes the Moore--Penrose pseudoinverse.
\end{lemma}

\begin{lemma}\label{lemma:B22_direct_sum}
    Let $K\hookrightarrow L$ be a simplicial pair and let $B_{q,2}^{K,L}:=[\partial^L_{q+1}](I_{q}^{K,L},I_{q+1}^{K,L})$. Then
    \begin{equation*}
        C_{q+1}^{K,L}\cong C^K_{q+1}\oplus\ker B_{q,2}^{K,L}.
    \end{equation*}
    We identify $\ker B_{q,2}^{K,L}\subseteq \mathbb{R}^{n_{q+1}^L-n_{q+1}^K}$ with the chains in $C_{q+1}^L$ supported on the $(q+1)$-simplices of $L\setminus K$. Using the identification, we have
    \begin{equation*}
        C_{q+1}^{K,L}= C^K_{q+1}\oplus\ker B_{q,2}^{K,L}.
    \end{equation*}
    % CHECK: state the identification once: ker B_{q,2}^{K,L} is a subspace of R^{n_{q+1}^L-n_{q+1}^K}, identified with the chains
    % supported on the (q+1)-simplices of L\K (padding with zeros). With this identification the sum is an orthogonal internal
    % direct sum in C_{q+1}^L, and "\cong" can be "=". All later "\cong" and "equal to" statements use this identification implicitly.
\end{lemma}

\subsection{Proof of Theorem \ref{thm:horizontal_operator}}
For this theorem, we use the following lemma.
\begin{lemma}\label{lemma:B22_horizontal}
    In a simplicial filtration $\{K_t\}_{0\leq t\leq T}$, for $0\leq s< t\leq T$, we have
    \begin{equation}\label{eq:B22_horizontal}
        (B_{q,2}^{s,t})^\dagger B_{q,2}^{s,t}=(B_{q,2}^{s,t})^\dagger B_{q,2}^{s,t}\begin{bmatrix}
        (B_{q,2}^{s,t-1})^\dagger B_{q,2}^{s,t-1} & 0\\
        0 & I_{n_{q+1}^t-n_{q+1}^{t-1}}
    \end{bmatrix}=\begin{bmatrix}
        (B_{q,2}^{s,t-1})^\dagger B_{q,2}^{s,t-1} & 0\\
        0 & I_{n_{q+1}^t-n_{q+1}^{t-1}}
    \end{bmatrix}(B_{q,2}^{s,t})^\dagger B_{q,2}^{s,t}.
    \end{equation}
\end{lemma}
We now prove Theorem~\ref{thm:horizontal_operator}.
\begin{proof}
    We prove the identity in the form $[\Delta_{q,+}^{s,t}]-[\Delta_{q,+}^{s,t-1}]=[\bar\partial_{q+1}^{s,t}(\bar\partial_{q+1}^{s,t})^*]$, which is more natural here.
    By Lemma~\ref{lemma:up_laplacian_boundary}, the difference equals
    \begin{align*}
        [\Delta_{q,+}^{s,t}]-[\Delta_{q,+}^{s,t-1}]=& B_{q,1}^{s,t}\left(I-(B_{q,2}^{s,t})^\dagger B_{q,2}^{s,t}\right) (B_{q,1}^{s,t})^T-B_{q,1}^{s,t-1}\left(I-(B_{q,2}^{s,t-1})^\dagger B_{q,2}^{s,t-1}\right)(B_{q,1}^{s,t-1})^T.
    \end{align*}
    By definition, $B_{q,1}^{s,t}$ consists of the columns of $B_{q,1}^{s,t-1}$, followed by the columns of the $(q+1)$-simplices of $K_t\setminus K_{t-1}$. Writing the second term in the same block form, the two terms can be combined:
    \begin{align*}
        [\Delta_{q, +}^{s,t}]-[\Delta_{q,+}^{s,t-1}]=&\begin{bmatrix}
            B_{q,1}^{s,t-1}&[\partial_{q+1}^t]([n_q^s],I_{q+1}^{t-1,t})
        \end{bmatrix}\left(I-(B_{q,2}^{s,t})^\dagger B_{q,2}^{s,t}\right)\begin{bmatrix}
            (B_{q,1}^{s,t-1})^T\\
            ([\partial_{q+1}^t]([n_q^s],I_{q+1}^{t-1,t}))^T
        \end{bmatrix}\\
        &-\begin{bmatrix}
            B_{q,1}^{s,t-1}&[\partial_{q+1}^t]([n_q^s],I_{q+1}^{t-1,t})
        \end{bmatrix}\left(I-\begin{bmatrix}
        (B_{q,2}^{s,t-1})^\dagger B_{q,2}^{s,t-1} & 0\\
        0 & I_{n_{q+1}^t-n_{q+1}^{t-1}}
    \end{bmatrix}\right)\begin{bmatrix}
            (B_{q,1}^{s,t-1})^T\\
            ([\partial_{q+1}^t]([n_q^s],I_{q+1}^{t-1,t}))^T
        \end{bmatrix}\nonumber\\
        =&B_{q,1}^{s,t}\left(\begin{bmatrix}
        (B_{q,2}^{s,t-1})^\dagger B_{q,2}^{s,t-1} & 0\\
        0 & I_{n_{q+1}^t-n_{q+1}^{t-1}}
    \end{bmatrix} -(B_{q,2}^{s,t})^\dagger B_{q,2}^{s,t}\right)(B_{q,1}^{s,t})^T.
    \end{align*}
    By Lemma~\ref{lemma:B22_horizontal}, we obtain
    \begin{align*}
        [\Delta_{q,+}^{s,t}]-[\Delta_{q,+}^{s,t-1}]&=B_{q,1}^{s,t}\left(\begin{bmatrix}
        (B_{q,2}^{s,t-1})^\dagger B_{q,2}^{s,t-1} & 0\\
        0 & I_{n_{q+1}^t-n_{q+1}^{t-1}}
    \end{bmatrix}-\begin{bmatrix}
        (B_{q,2}^{s,t-1})^\dagger B_{q,2}^{s,t-1} & 0\\
        0 & I_{n_{q+1}^t-n_{q+1}^{t-1}}
    \end{bmatrix}(B_{q,2}^{s,t})^\dagger B_{q,2}^{s,t}\right)(B_{q,1}^{s,t})^T\\
    &=B_{q,1}^{s,t}\begin{bmatrix}
        (B_{q,2}^{s,t-1})^\dagger B_{q,2}^{s,t-1} & 0\\
        0 & I_{n_{q+1}^t-n_{q+1}^{t-1}}
    \end{bmatrix}\left(I-(B_{q,2}^{s,t})^\dagger B_{q,2}^{s,t}\right)(B_{q,1}^{s,t})^T.
    \end{align*}
    Note that
    \begin{equation*}
        \begin{bmatrix}
        (B_{q,2}^{s,t-1})^\dagger B_{q,2}^{s,t-1} & 0\\
        0 & I_{n_{q+1}^t-n_{q+1}^{t-1}}
    \end{bmatrix}=\begin{bmatrix}
        B_{q,2}^{s,t-1} & 0\\
        0 & I_{n_{q+1}^t-n_{q+1}^{t-1}}
    \end{bmatrix}^\dagger \begin{bmatrix}
        B_{q,2}^{s,t-1} & 0\\
        0 & I_{n_{q+1}^t-n_{q+1}^{t-1}}
    \end{bmatrix}.
    \end{equation*}
    Therefore, it is the orthogonal projection onto $\left(\ker \begin{bmatrix}
        B_{q,2}^{s,t-1} & 0\\
        0 & I_{n_{q+1}^t-n_{q+1}^{t-1}}
    \end{bmatrix}\right)^{\mathrel\bot}$. Similarly, $I-(B_{q,2}^{s,t})^\dagger B_{q,2}^{s,t}$ is the orthogonal projection onto $\left(\left(\ker B_{q,2}^{s,t}\right)^{\mathrel\bot}\right)^{\mathrel\bot}=\ker B_{q,2}^{s,t}$, see Section \ref{sec:preliminaries}. By Lemma~\ref{lemma:B22_horizontal}, these two projections commute. Therefore, their product is the orthogonal projection onto the intersection of the two spaces~\cite{projection_matrices}, $\Psi_h:=\ker B_{q,2}^{s,t}\cap\left(\ker \begin{bmatrix}
        B_{q,2}^{s,t-1} & 0\\
        0 & I_{n_{q+1}^t-n_{q+1}^{t-1}}
    \end{bmatrix}\right)^{\mathrel\bot}$.

    Let $[\Psi_h]$ be an orthonormal basis of $\Psi_h$, and let $V$ be the matrix whose columns are the vectors in $[\Psi_h]$. Then this projection equals $VV^T$, and we obtain
    \[
      [\Delta_{q,+}^{s,t}]-[\Delta_{q,+}^{s,t-1}]=B_{q,1}^{s,t}VV^T(B_{q,1}^{s,t})^T=B_{q,1}^{s,t}V\left(B_{q,1}^{s,t}V\right)^T.
    \]
    We claim that $B_{q,1}^{s,t}V=[\bar\partial_{q+1}^{s,t}]$, which proves the theorem. To show this, we describe $\widebar{C_{q+1}^{s,t}}$ using Lemma~\ref{lemma:B22_direct_sum}:
    \[
        \widebar{C_{q+1}^{s,t}}\cong \left(C^s_{q+1}\oplus\ker B_{q,2}^{s,t}\right)\cap\left(C^s_{q+1}\oplus\ker B_{q,2}^{s,t-1}\oplus \ker I_{n_{q+1}^t-n_{q+1}^{t-1}}\right)^\perp,
    \]
    with the complement taken in $C_{q+1}^t$. We show that the following holds,
    \begin{equation}\label{eq:intersect_sum_horizontal}
    \left(C^s_{q+1}\oplus\ker B_{q,2}^{s,t}\right)\cap\left(C^s_{q+1}\oplus\ker B_{q,2}^{s,t-1}\oplus \ker I_{n_{q+1}^t-n_{q+1}^{t-1}}\right)^\perp\cong\ker B_{q,2}^{s,t}\cap\left(\ker \begin{bmatrix}
        B_{q,2}^{s,t-1} & 0\\
        0 & I_{n_{q+1}^t-n_{q+1}^{t-1}}
    \end{bmatrix}\right)^\perp.
    \end{equation} 
    We have $v\in \left(C^s_{q+1}\oplus\ker B_{q,2}^{s,t}\right)\cap\left(C^s_{q+1}\oplus\ker \begin{bsmallmatrix}
        B_{q,2}^{s,t-1} & 0\\
        0 & I_{n_{q+1}^t-n_{q+1}^{t-1}}
    \end{bsmallmatrix}\right)^\perp$ if and only if $v=v_\ell+v_r$ with $v_\ell\in C_{q+1}^s, v_r\in\ker B_{q,2}^{s,t}$ and $\langle v, w\rangle=0$ for all $w\in C^s_{q+1}\oplus\ker \begin{bsmallmatrix}
        B_{q,2}^{s,t-1} & 0\\
        0 & I_{n_{q+1}^t-n_{q+1}^{t-1}}
    \end{bsmallmatrix}$. Decomposing $w$ in the same way, $\langle v, w\rangle=0$ holds if and only if
    \[
    \langle v, w\rangle=\langle v_\ell,w_\ell\rangle+\langle v_\ell,w_r\rangle+\langle v_r,w_\ell\rangle + \langle v_r,w_r\rangle=0.
    \]
    Here, $\langle v_\ell,w_r\rangle=\langle v_r,w_\ell\rangle=0$, since these vectors have disjoint supports. Moreover, if $v_\ell\neq 0$, $\langle v_\ell, w_\ell\rangle$ can take any real value $c\in\mathbb{R}$ by choosing $w_\ell=c\frac{v_\ell}{\langle v_\ell,v_\ell\rangle}$.
    % CHECK: this only works if v_\ell\neq 0 (division by <v_\ell,v_\ell>). Simpler: take w=v_\ell (so w_r=0), which gives <v_\ell,v_\ell>=0.
    Therefore, the equality holds for every $w$ if and only if $v_\ell=0$. This leaves $\langle v_r, w_r\rangle=0$, which holds if and only if $v_r\in\ker B_{q,2}^{s,t}\cap\left(\ker \begin{bsmallmatrix}
        B_{q,2}^{s,t-1} & 0\\
        0 & I_{n_{q+1}^t-n_{q+1}^{t-1}}
    \end{bsmallmatrix}\right)^\perp$.

    Finally, the columns of $V$ form a basis of $\widebar{C_{q+1}^{s,t}}$, and $B_{q,1}^{s,t}V$ describes the action of the boundary operator on this basis, therefore $B_{q,1}^{s,t}V=[\bar\partial_{q+1}^{s,t}]$.
\end{proof}

\subsection{Proof of Theorem \ref{thm:vertical_operator}}
\begin{figure}
    \centering
    \includegraphics[width=0.75\linewidth]{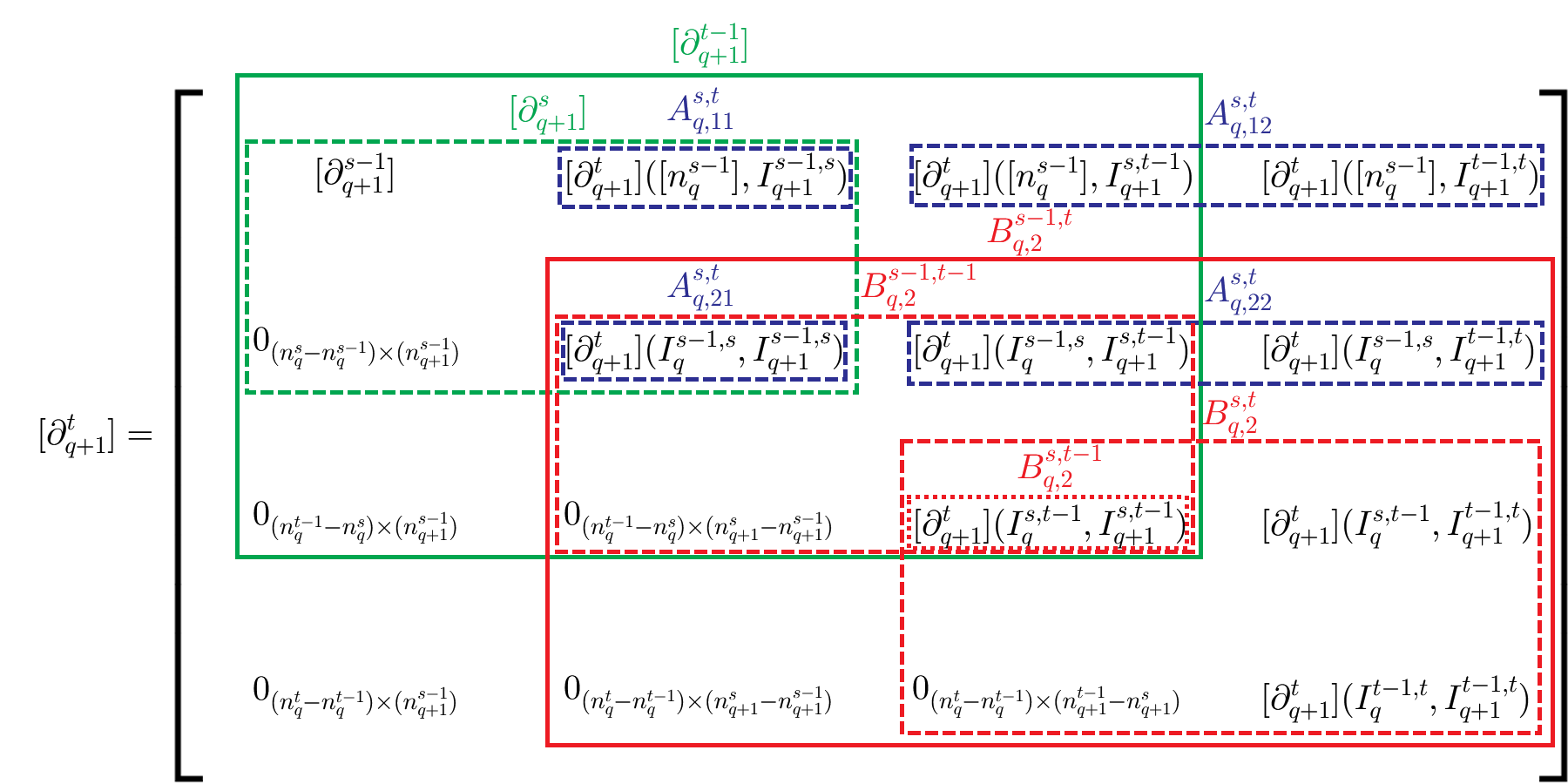}
    \caption{The boundary matrix $[\partial_{q+1}^t]$ split into the different submatrices. Each zero block is denoted by $0_{m\times n}$, where $m$ and $n$ are its numbers of rows and columns, respectively.}
    \label{fig:matrices_in_one}
\end{figure}

We use the following lemma, which is analogous to Lemma~\ref{lemma:B22_horizontal}.
\begin{lemma}\label{lemma:B22_vertical}
    In a simplicial filtration $\{K_t\}_{0\leq t\leq T}$, for $1\leq s\leq t$, we have
    \begin{equation}\label{eq:B22_vertical}
        \begin{bmatrix}
            0_{n_{q+1}^s-n_{q+1}^{s-1}}&0\\
            0&(B_{q,2}^{s,t})^\dagger B_{q,2}^{s,t}
        \end{bmatrix}=\begin{bmatrix}
            0_{n_{q+1}^s-n_{q+1}^{s-1}}&0\\
            0&(B_{q,2}^{s,t})^\dagger B_{q,2}^{s,t}
        \end{bmatrix}(B_{q,2}^{s-1,t})^\dagger B_{q,2}^{s-1,t}=(B_{q,2}^{s-1,t})^\dagger B_{q,2}^{s-1,t}\begin{bmatrix}
            0_{n_{q+1}^s-n_{q+1}^{s-1}}&0\\
            0&(B_{q,2}^{s,t})^\dagger B_{q,2}^{s,t}
        \end{bmatrix}.
    \end{equation}
\end{lemma}
In the proof, we use the following notation.
\begin{align}\label{eq:A_partition}
    A_q^{s,t} &:=[\partial_{q+1}^t]([n_q^s],I_{q+1}^{s-1,t})=\begin{bmatrix}
        A^{s,t}_{q,11} & A^{s,t}_{q,12}\\
        A^{s,t}_{q,21} & A^{s,t}_{q,22}
    \end{bmatrix},\\
    A^{s,t}_{q,11}&:=[\partial_{q+1}^t]([n_q^{s-1}], I_{q+1}^{s-1,s}), & A^{s,t}_{q,12}&:=[\partial_{q+1}^t]([n_q^{s-1}], I_{q+1}^{s,t}),\nonumber\\
    A^{s,t}_{q,21}&:=[\partial_{q+1}^t](I_q^{s-1,s}, I_{q+1}^{s-1,s}), & A^{s,t}_{q,22}&:=[\partial_{q+1}^t](I_q^{s-1,s}, I_{q+1}^{s,t}),\nonumber\\
    B_{q,1}^{s-1,t}&=\begin{bmatrix}
        A^{s,t}_{q,11} & A^{s,t}_{q,12}
    \end{bmatrix}, & B_{q,1}^{s,t}&=\begin{bmatrix}
        A^{s,t}_{q,12}\\
        A^{s,t}_{q,22}
    \end{bmatrix}.\nonumber
\end{align}
Figure~\ref{fig:matrices_in_one} shows how these matrices are obtained from the boundary matrix $[\partial_{q+1}^t]$.

We now prove Theorem~\ref{thm:vertical_operator}.
\begin{proof}
    Again, we prove the identity in the form $[\Delta_{q,+}^{s,t}]-[\iota\Delta_{q,+}^{s-1,t}\iota^*]=[\tilde\partial_{q+1}^{s,t}(\tilde\partial_{q+1}^{s,t})^*]$. First, note that
    \[
        [\iota\Delta_{q,+}^{s-1,t}\iota^*]=\begin{bmatrix}
        I_{n_q^{s-1}}\\
        0_{ (n_q^s-n_q^{s-1})\times n_q^{s-1}}
    \end{bmatrix}[\Delta_{q,+}^{s-1,t}]\begin{bmatrix}
        I_{n_q^{s-1}} & 0_{n_q^{s-1}\times (n_q^s-n_q^{s-1})}
    \end{bmatrix}=\begin{bmatrix}
        [\Delta_{q,+}^{s-1,t}] & 0\\
        0 & 0_{(n_q^s-n_q^{s-1})\times (n_q^s-n_q^{s-1})}
    \end{bmatrix}.
    \]

    Using the notation~\eqref{eq:A_partition} and Lemma~\ref{lemma:up_laplacian_boundary}, we can rewrite $[\Delta_{q,+}^{s,t}]$ as
    \begin{align*}
        [\Delta_{q,+}^{s,t}]%&=[\Delta_{q,+}^s]+B_{q,1}^{s,t}(I-(B_{q,2}^{s,t})^\dagger B_{q,2}^{s,t})(B_{q,1}^{s,t})^T\\
        &=[\partial_{q+1}^s][\partial_{q+1}^s]^T+B_{q,1}^{s,t}(B_{q,1}^{s,t})^T-B_{q,1}^{s,t}(B_{q,2}^{s,t})^\dagger B_{q,2}^{s,t}(B_{q,1}^{s,t})^T\\
        &=\begin{bmatrix}
        [\partial_{q+1}^{s-1}][\partial_{q+1}^{s-1}]^T+A^{s,t}_{q,11}\left(A^{s,t}_{q,11}\right)^T & A^{s,t}_{q,11}\left(A^{s,t}_{q,21}\right)^T\\
        A^{s,t}_{q,21}\left(A^{s,t}_{q,11}\right)^T & A^{s,t}_{q,21}\left(A^{s,t}_{q,21}\right)^T
    \end{bmatrix} + \begin{bmatrix}
        A^{s,t}_{q,12}\left(A^{s,t}_{q,12}\right)^T & A^{s,t}_{q,12}\left(A^{s,t}_{q,22}\right)^T\\
        A^{s,t}_{q,22}\left(A^{s,t}_{q,12}\right)^T & A^{s,t}_{q,22}\left(A^{s,t}_{q,22}\right)^T
    \end{bmatrix}-A_q^{s,t}\begin{bmatrix}
        0&0\\
        0&(B_{q,2}^{s,t})^\dagger B_{q,2}^{s,t}
    \end{bmatrix}(A_q^{s,t})^T\\
    &=\begin{bmatrix}
        [\Delta_{q,+}^{s-1}] & 0\\
        0 & 0
    \end{bmatrix}+A_q^{s,t}\left(I-\begin{bmatrix}
        0&0\\
        0&(B_{q,2}^{s,t})^\dagger B_{q,2}^{s,t}
    \end{bmatrix}\right)(A_q^{s,t})^T.
    \end{align*}

    Furthermore,
    \begin{equation*}
        B_{q,2}^{s-1,t}=\begin{bmatrix}
        A_{q,21}^{s,t} & A_{q,22}^{s,t}\\
        0 & B_{q,2}^{s,t}
    \end{bmatrix}.
    \end{equation*}
    Hence, $\ker B_{q,2}^{s-1,t}\subseteq\ker \begin{bmatrix} A_{q,21}^{s,t} & A_{q,22}^{s,t}\end{bmatrix}$, and Equation~\eqref{eq:pinv_subset_ker} gives
    \begin{equation}\label{eq:B_sm1_subset_A}
        \begin{bmatrix} A_{q,21}^{s,t} & A_{q,22}^{s,t}\end{bmatrix}=\begin{bmatrix} A_{q,21}^{s,t} & A_{q,22}^{s,t}\end{bmatrix} (B_{q,2}^{s-1,t})^\dagger B_{q,2}^{s-1,t}.
    \end{equation}

    Using Lemma~\ref{lemma:up_laplacian_boundary} again, we can rewrite $[\iota\circ\Delta_{q,+}^{s-1,t}\circ\iota^*]$ as follows, where the last step uses~\eqref{eq:B_sm1_subset_A}:
    \begin{align*}
        [\iota\circ\Delta_{q,+}^{s-1,t}\circ\iota^*] &= \begin{bmatrix}
            [\Delta_{q,+}^{s-1}]+B_{q,1}^{s-1,t}(I-(B_{q,2}^{s-1,t})^\dagger B_{q,2}^{s-1,t})(B_{q,1}^{s-1,t})^T & 0\\
            0 & 0
        \end{bmatrix}\\
        &= \begin{bmatrix}
        [\Delta_{q,+}^{s-1}] & 0\\
        0 & 0
    \end{bmatrix}+\begin{bmatrix}
        A_{q, 11}^{s,t} & A_{q, 12}^{s,t}\\
        0 & 0
    \end{bmatrix}(I-(B_{q,2}^{s-1,t})^\dagger B_{q,2}^{s-1,t})\begin{bmatrix}
        \left(A_{q, 11}^{s,t}\right)^T & 0\\
        \left(A_{q, 12}^{s,t}\right)^T & 0
    \end{bmatrix}\\
    &=\begin{bmatrix}
        [\Delta_{q,+}^{s-1}] & 0\\
        0 & 0
    \end{bmatrix}+\begin{bmatrix}
        A_{q, 11}^{s,t} & A_{q, 12}^{s,t}\\
        A_{q,21}^{s,t} & A_{q,22}^{s,t}
    \end{bmatrix}(I-(B_{q,2}^{s-1,t})^\dagger B_{q,2}^{s-1,t})\begin{bmatrix}
        \left(A_{q, 11}^{s,t}\right)^T & \left(A_{q,21}^{s,t}\right)^T\\
        \left(A_{q, 12}^{s,t}\right)^T & \left(A_{q,22}^{s,t}\right)^T
    \end{bmatrix}.
    \end{align*}

    Hence, the difference is
    \begin{align*}
        [\Delta_{q,+}^{s,t}]-[\iota\circ\Delta_{q,+}^{s-1,t}\circ\iota^*]&=A_q^{s,t}\left((B_{q,2}^{s-1,t})^\dagger B_{q,2}^{s-1,t}-\begin{bmatrix}
        0&0\\
        0&(B_{q,2}^{s,t})^\dagger B_{q,2}^{s,t}
    \end{bmatrix}\right)(A_q^{s,t})^T.
    \end{align*}
    By Lemma~\ref{lemma:B22_vertical}, we obtain
    \begin{align*}
        [\Delta_{q,+}^{s,t}]-[\iota\circ\Delta_{q,+}^{s-1,t}\circ\iota^*]&=A_q^{s,t}(B_{q,2}^{s-1,t})^\dagger B_{q,2}^{s-1,t}\left(I-\begin{bmatrix}
            0&0\\0&(B_{q,2}^{s,t})^\dagger B_{q,2}^{s,t}
        \end{bmatrix}\right)(A_q^{s,t})^T.
    \end{align*}
    By the same reasoning as in the proof of Theorem~\ref{thm:horizontal_operator}, we obtain, for $\Psi_v:=\ker \begin{bmatrix}
            0&0\\0&B_{q,2}^{s,t}
        \end{bmatrix}\cap\left(\ker B_{q,2}^{s-1,t}\right)^{\mathrel\bot}$ and a matrix $V_v$ whose columns form an orthonormal basis of $\Psi_v$,
    \[
    [\Delta_{q,+}^{s,t}]-[\iota\circ\Delta_{q,+}^{s-1,t}\circ\iota^*]=A_q^{s,t}V_vV_v^T(A_q^{s,t})^T.
    \]
    We split $C_{q+1}^s$ into $C_{q+1}^{s-1}$ and its orthogonal complement, $C_{q+1}^s=C_{q+1}^{s-1}\oplus \left(C_{q+1}^s\cap(C_{q+1}^{s-1})^\perp\right)$. The orthogonal complement $C_{q+1}^s\cap(C_{q+1}^{s-1})^\perp$ is spanned by the $(q+1)$-simplices of $K_s\setminus K_{s-1}$, so
    \[
    C_{q+1}^s\cap(C_{q+1}^{s-1})^\perp\cong \ker \begin{bmatrix}
        0_{n_{q+1}^{s}-n_{q+1}^{s-1}}.
    \end{bmatrix}
    \]
    Combining this with Lemma~\ref{lemma:B22_direct_sum}, we obtain
    \begin{align*}
        C_{q+1}^{s,t}&\cong C_{q+1}^s\oplus\ker B_{q,2}^{s,t}\\
        &\cong C_{q+1}^{s-1}\oplus\ker \begin{bmatrix}
        0_{n_{q+1}^{s}-n_{q+1}^{s-1}}
        \end{bmatrix}\oplus\ker B_{q,2}^{s,t}\\
        &= C_{q+1}^{s-1}\oplus\ker \begin{bmatrix}
        0_{n_{q+1}^{s}-n_{q+1}^{s-1}} & 0\\
        0 & B_{q,2}^{s,t}
        \end{bmatrix}.
    \end{align*}
    This allows for the comparison between the consecutive horizontal chain spaces using again Lemma~\ref{lemma:B22_direct_sum}. 
    \begin{align*}
        \widetilde{C_{q+1}^{s,t}}&=C_{q+1}^{s,t}\cap (C_{q+1}^{s-1,t})^\perp\\
        &\cong \left(C_{q+1}^{s-1}\oplus\ker \begin{bmatrix}
        0_{n_{q+1}^{s}-n_{q+1}^{s-1}} & 0\\
        0 & B_{q,2}^{s,t}
        \end{bmatrix}\right)\cap\left(C_{q+1}^{s-1}\oplus\ker B_{q,2}^{s-1,t}\right)^\perp.
    \end{align*}
    Now, using the same argument as in the proof of Theorem~\ref{thm:horizontal_operator}, we conclude that $A_q^{s,t}V_v=[\tilde\partial_{q+1}^{s,t}]$.
\end{proof}

\section{Proof of Lemma \ref{lemma:horizontal_vertical_distinct}}\label{sec:proof_lemma_distinct}
\begin{proof}
    (i) For $s=t$, the left-hand side is $\{0\}$ by the convention $C_{q+1}^{t,t-1}=\{0\}$, so let $s\leq t-1$. Let $c\in C_{q+1}^{s-1,t}\cap C_{q+1}^{s,t-1}$. Since $c\in C_{q+1}^{s,t-1}$, the chain $c$ lies in $C_{q+1}^{t-1}$, and since $c\in C_{q+1}^{s-1,t}$, its boundary lies in $C_q^{s-1}$. Hence, $c\in C_{q+1}^{s-1,t-1}$. Conversely, $C_{q+1}^{s-1,t-1}$ is contained in both $C_{q+1}^{s-1,t}$ and $C_{q+1}^{s,t-1}$, which gives equality.
 
    (ii) By definition, $\widebar{C_{q+1}^{s-1,t}}\subseteq C_{q+1}^{s-1,t}$ and $\widetilde{C_{q+1}^{s,t-1}}\subseteq C_{q+1}^{s,t-1}$, and both spaces are orthogonal to $C_{q+1}^{s-1,t-1}$. By (i), their intersection is therefore contained in $C_{q+1}^{s-1,t-1}\cap\big(C_{q+1}^{s-1,t-1}\big)^\perp=\{0\}$.
\end{proof}

\section{Proof of Theorem \ref{thm:new_chain_decomposition}}\label{sec:proof_chain_decomposition}
\begin{proof}
    Define
    \[
    \mathcal{S}_{q+1}^{s,t}:= \bigoplus_{\substack{0\leq s^-\leq s,\ t^-\leq t,\\ s^-\leq t^-}} N_{q+1}^{s^-,t^-}.
    \]
    The goal is to prove $\mathcal{S}_{q+1}^{s,t}=C_{q+1}^{s,t}$, where this statement includes that the sum defining $\mathcal{S}_{q+1}^{s,t}$ is direct. We do this by induction on $t$ and, for fixed $t$, on $s$.

    \textbf{Base case $t$:} Let $t=0$. Then $s\leq t=0$ implies $s=0$, and everything simplifies:
    \[
        \mathcal{S}_{q+1}^{0,0}=N_{q+1}^{0,0}=C_{q+1}^{0,0}\cap(\{0\})^\perp=C_{q+1}^{0,0}.
    \]

    \textbf{Induction step $t$:} Assume that $\mathcal{S}_{q+1}^{s,t-1}=C_{q+1}^{s,t-1}$ for all $s\leq t-1$. We have to show that $\mathcal{S}_{q+1}^{s,t}=C_{q+1}^{s,t}$ for all $s\leq t$. To do this, we use induction on $s$.

    \textbf{Base case $s$:} For $s=0$, we have $\mathcal{S}_{q+1}^{0, t}=\bigoplus_{t^-\leq t}N_{q+1}^{0, t^-}=N_{q+1}^{0,t}\oplus \mathcal{S}_{q+1}^{0,t-1}$. Here, $N_{q+1}^{0, t}=C_{q+1}^{0,t}\cap (C_{q+1}^{0,t-1})^\perp$ by the convention $C_{q+1}^{-1,t}=\{0\}$. Using the induction assumption, we obtain $\mathcal{S}_{q+1}^{0, t}=\left(C_{q+1}^{0,t}\cap (C_{q+1}^{0,t-1})^\perp\right)\oplus C_{q+1}^{0,t-1}=C_{q+1}^{0,t}$ by the orthogonal decomposition theorem. This sum is direct, since $N_{q+1}^{0,t}$ is orthogonal to $C_{q+1}^{0,t-1}$.

    \textbf{Induction step $s$:} Let $1\leq s\leq t$, and assume that $\mathcal{S}_{q+1}^{s-1,t}=C_{q+1}^{s-1,t}$. It remains to show $\mathcal{S}_{q+1}^{s,t}=C_{q+1}^{s,t}$, which proves both inductions.

    First, note that we can decompose $\mathcal{S}_{q+1}^{s,t}$ using the induction hypothesis on $t$:
    \begin{align*}
        \mathcal{S}_{q+1}^{s,t}&=\mathcal{S}_{q+1}^{s-1,t-1}\oplus\bigoplus_{s\leq t^-\leq t-1}N_{q+1}^{s,t^-}\oplus \bigoplus_{0\leq s^-\leq s-1}N_{q+1}^{s^-,t}\oplus N_{q+1}^{s,t}\\
        &=C_{q+1}^{s-1,t-1}\oplus\bigoplus_{s\leq t^-\leq t-1}N_{q+1}^{s,t^-}\oplus \bigoplus_{0\leq s^-\leq s-1}N_{q+1}^{s^-,t}\oplus N_{q+1}^{s,t}.
    \end{align*}
    Here, we have only regrouped the summands of $\mathcal{S}_{q+1}^{s,t}$. At the end of the proof, we show that this sum is direct. We write
    \[
        \mathcal{N}_{v}:=\bigoplus_{s\leq t^-\leq t-1}N_{q+1}^{s,t^-},
        \qquad
        \mathcal{N}_{h}:=\bigoplus_{0\leq s^-\leq s-1}N_{q+1}^{s^-,t}.
    \]
    If $s=t$, the first sum is empty, so $\mathcal{N}_v=\{0\}$.
    We use the decomposition~\eqref{eq:new_chain_decomposition}, which by the definition of $N_{q+1}^{s,t}$ reads
    \begin{equation*}
        C_{q+1}^{s,t}=\widebar{C_{q+1}^{s-1,t}}\oplus \widetilde{C_{q+1}^{s,t-1}}\oplus C_{q+1}^{s-1,t-1}\oplus N_{q+1}^{s,t}.
    \end{equation*}
    Since $\widebar{C_{q+1}^{s-1,t}}$ and $\widetilde{C_{q+1}^{s,t-1}}$ are orthogonal to $C_{q+1}^{s-1,t-1}$, while the spaces $N_{q+1}^{s^-,t^-}$ in general are not, the spaces $\mathcal{N}_v\oplus\mathcal{N}_h$ and $\widebar{C_{q+1}^{s-1,t}}\oplus \widetilde{C_{q+1}^{s,t-1}}$ are in general not equal. Instead, we show that they agree up to $C_{q+1}^{s-1,t-1}$.

    By the induction hypotheses and the orthogonal decomposition theorem, we obtain
    \[
        C_{q+1}^{s-1,t-1}\oplus\mathcal{N}_h=\mathcal{S}_{q+1}^{s-1,t-1}\oplus\mathcal{N}_h=\mathcal{S}_{q+1}^{s-1,t}=C_{q+1}^{s-1,t}=C_{q+1}^{s-1,t-1}\oplus\widebar{C_{q+1}^{s-1,t}}.
    \]
    Similarly, if $s\leq t-1$,
    \[
        C_{q+1}^{s-1,t-1}\oplus\mathcal{N}_v=\mathcal{S}_{q+1}^{s-1,t-1}\oplus\mathcal{N}_v=\mathcal{S}_{q+1}^{s,t-1}=C_{q+1}^{s,t-1}=C_{q+1}^{s-1,t-1}\oplus\widetilde{C_{q+1}^{s,t-1}}.
    \]
    For $s=t$, both $\mathcal{N}_v$ and $\widetilde{C_{q+1}^{t,t-1}}$ are $\{0\}$. In these decompositions, the sums on the left are direct but not necessarily orthogonal, while those on the right are orthogonal. Consequently,
    \[
        C_{q+1}^{s-1,t-1}+\mathcal{N}_v+\mathcal{N}_h=C_{q+1}^{s-1,t-1}+\widetilde{C_{q+1}^{s,t-1}}+\widebar{C_{q+1}^{s-1,t}},
    \]
    so that, by the decomposition~\eqref{eq:new_chain_decomposition},
    \[
        \mathcal{S}_{q+1}^{s,t}=C_{q+1}^{s-1,t-1}+\mathcal{N}_v+\mathcal{N}_h+N_{q+1}^{s,t}=C_{q+1}^{s,t}.
    \]

    It remains to show that the sum is direct. Let $x_0+x_v+x_h+x_N=0$ with $x_0\in C_{q+1}^{s-1,t-1}$, $x_v\in\mathcal{N}_v$, $x_h\in\mathcal{N}_h$ and $x_N\in N_{q+1}^{s,t}$. Since $x_0+x_v+x_h$ lies in $C_{q+1}^{s-1,t-1}+\widetilde{C_{q+1}^{s,t-1}}+\widebar{C_{q+1}^{s-1,t}}$, which is orthogonal to $N_{q+1}^{s,t}$, we obtain $x_N=0$.

    Next, the vector $x_v=-(x_0+x_h)$ lies in $C_{q+1}^{s-1,t}$, and since $x_v\in\mathcal{N}_v\subseteq C_{q+1}^{s,t-1}$, Lemma~\ref{lemma:horizontal_vertical_distinct}(i) gives $x_v\in C_{q+1}^{s-1,t-1}$. Hence, $x_v\in C_{q+1}^{s-1,t-1}\cap\mathcal{N}_v=\{0\}$: for $s=t$, this holds since $\mathcal{N}_v=\{0\}$, and for $s\leq t-1$, since the sum $C_{q+1}^{s-1,t-1}\oplus\mathcal{N}_v=\mathcal{S}_{q+1}^{s,t-1}$ is direct.

    It follows that $x_0+x_h=0$, and since the sum $C_{q+1}^{s-1,t-1}\oplus\mathcal{N}_h=\mathcal{S}_{q+1}^{s-1,t}$ is direct, $x_0=x_h=0$. Finally, $\mathcal{S}_{q+1}^{s-1,t-1}$, $\mathcal{N}_v$ and $\mathcal{N}_h$ are themselves direct sums by the induction hypotheses, since they are parts of $\mathcal{S}_{q+1}^{s,t-1}$ and $\mathcal{S}_{q+1}^{s-1,t}$. Therefore, the sum $\mathcal{S}_{q+1}^{s,t}$ is direct, which completes both inductions.
\end{proof}

\section{Proofs of the remaining lemmas}
\label{ch:lemma_proofs}
\subsection{Proof of Lemma \ref{lemma:up_laplacian_boundary}}
\begin{proof}
    We first consider the case when $n_{q+1}^K=n_{q+1}^L$. As stated in \cite{laplacians_algorithms}, we then have $[\Delta_{q,+}^{K,L}]=[\Delta_{q,+}^{K}]$, which proves the lemma.
    
    Second, let $n_q^K=n_q^L$. We can then write the matrix representation of the boundary operator as
    \begin{equation*}
        [\partial_{q+1}^L]=\begin{bmatrix}
        [\partial_{q+1}^K]&B_{q,1}^{K,L}
    \end{bmatrix}.
    \end{equation*}
    In this setting all $(q+1)$-simplices in $L$ have a boundary in $K$. Therefore $\partial_{q+1}^{K,L}=\partial_{q+1}^L$ and we can compute the persistent up-Laplacian.
    \begin{equation*}
        [\Delta_{q,+}^{K,L}]=[\partial_{q+1}^L][\partial_{q+1}^L]^T=[\partial_{q+1}^K][\partial_{q+1}^K]^T+B_{q,1}^{K,L}\left(B_{q,1}^{K,L}\right)^T=[\Delta_{q,+}^{K}]+B_{q,1}^{K,L}\left(B_{q,1}^{K,L}\right)^T.
    \end{equation*}
    
    Finally, when $n_q^K<n_q^L$, the matrix of the boundary operator $\partial_{q+1}^L$ can be written as
    \begin{equation*}
        [\partial_{q+1}^L]=\begin{bmatrix}
        [\partial_{q+1}^K]&B_{q,1}^{K,L}\\
        0&B_{q,2}^{K,L}
    \end{bmatrix}.
    \end{equation*}
    Hence, the combinatorial up-Laplacian is \\
    \[[\Delta_{q,+}^L]=[\partial_{q+1}^L][\partial_{q+1}^L]^T=\begin{bmatrix}
        [\partial_{q+1}^K][\partial_{q+1}^K]^T+B_{q,1}^{K,L}(B_{q,1}^{K,L})^T&B_{q,1}^{K,L}(B_{q,2}^{K,L})^T\\
        B_{q,2}^{K,L}(B_{q,1}^{K,L})^T&B_{q,2}^{K,L}(B_{q,2}^{K,L})^T
    \end{bmatrix}.
    \]
    Applying the Schur complement representation of Theorem~\ref{thm:laplacian_schur}, we get
    \begin{align*}
        [\Delta_{q,+}^{K,L}]&=[\Delta_{q,+}^{L}]/([\Delta_{q,+}^{L}](I_q^{K,L},I_q^{K,L}))\\
        &=[\partial_{q+1}^K][\partial_{q+1}^K]^T+B_{q,1}^{K,L}(B_{q,1}^{K,L})^T-B_{q,1}^{K,L}(B_{q,2}^{K,L})^T\left[B_{q,2}^{K,L}(B_{q,2}^{K,L})^T\right]^\dagger B_{q,2}^{K,L}(B_{q,1}^{K,L})^T\\
        &=[\Delta_{q,+}^K]+B_{q,1}^{K,L}\left(I-(B_{q,2}^{K,L})^T\left[B_{q,2}^{K,L}(B_{q,2}^{K,L})^T\right]^\dagger B_{q,2}^{K,L}\right)(B_{q,1}^{K,L})^T.
    \end{align*}
    Using the property~\eqref{eq:pinv_AAT} of the pseudoinverse and writing $Q=(B_{q,2}^{K,L})^\dagger B_{q,2}^{K,L}$, we get
    \begin{align*}
        [\Delta_{q,+}^{K,L}]&=[\Delta_{q,+}^K]+B_{q,1}^{K,L}\left(I-(B_{q,2}^{K,L})^T((B_{q,2}^{K,L})^\dagger)^T(B_{q,2}^{K,L})^\dagger B_{q,2}^{K,L}\right)(B_{q,1}^{K,L})^T\\
        &=[\Delta_{q,+}^K]+B_{q,1}^{K,L}\left(I-Q^TQ\right)(B_{q,1}^{K,L})^T.
    \end{align*}
    Finally, since $Q$ is symmetric and idempotent, see Equations~\eqref{eq:Q_symmetric} and~\eqref{eq:Q_idempotent}, we obtain
    \begin{align*}
        [\Delta_{q,+}^{K,L}]&=[\Delta_{q,+}^K]+B_{q,1}^{K,L}\left(I-Q\right)(B_{q,1}^{K,L})^T\\
        &=[\Delta_{q,+}^K]+B_{q,1}^{K,L}(I-(B_{q,2}^{K,L})^\dagger B_{q,2}^{K,L})(B_{q,1}^{K,L})^T.
    \end{align*}
\end{proof}

\subsection{Proof of Lemma \ref{lemma:B22_direct_sum}}
\begin{proof}
    Let $c\in C_{q+1}^L$, and write $c=\sum_{i=1}^{n_{q+1}^L}a_i\sigma_i^{q+1}$ in the simplex basis, with $a_i\in\mathbb{R}$ and $\sigma_i^{q+1}\in S^L_{q+1}$ the $(q+1)$-simplices of $L$. Let $c_K\in C_{q+1}^K$ denote the part of $c$ in $C_{q+1}^K$, that is, $c_K=\sum_{i=1}^{n_{q+1}^K}a_i\sigma_i^{q+1}$, and let $c_{K,L}$ denote the remaining part, $c_{K,L}=\sum_{i=n_{q+1}^K+1}^{n_{q+1}^L}a_i\sigma_i^{q+1}$. It remains to show that $c\in C_{q+1}^{K,L}$ if and only if $c_{K,L}\in\ker B_{q,2}^{K,L}$.

    For $n_q^K<j\leq n_q^L$ and $n_{q+1}^K<i\leq n_{q+1}^L$, the entries of $B_{q,2}^{K,L}$ are
    \begin{equation*}
        (B_{q,2}^{K,L})_{j-n_q^K,i-n_{q+1}^K}=
        \left\langle\partial_{q+1}^L\sigma_i^{q+1}, \sigma_j^q\right\rangle.
    \end{equation*}
    Hence, the entries of $B_{q,2}^{K,L}[c_{K,L}]$ are
    \begin{equation*}
        (B_{q,2}^{K,L}\cdot[c_{K,L}])_{j-n_q^K}=
        \sum_{i=n_{q+1}^K+1}^{n_{q+1}^L}a_i\left\langle\partial_{q+1}^L\sigma_i^{q+1}, \sigma_j^q\right\rangle.
    \end{equation*}
    Consequently, $c_{K,L}\in\ker B_{q,2}^{K,L}$ if and only if
    $\sum_{i=n_{q+1}^K+1}^{n_{q+1}^L}a_i\left\langle\partial_{q+1}^L\sigma_i^{q+1}, \sigma_j^q\right\rangle=0$
    for all $n_q^K<j\leq n_q^L$.

    We now show that $c\in C_{q+1}^{K,L}$ is also equivalent to this. Applying the boundary operator to $c$ yields
    \begin{align*}
        \partial_{q+1}^L(c)&=\sum_{i=1}^{n_{q+1}^L}a_i\partial_{q+1}^L\sigma_i^{q+1}\\
        &=\sum_{i=1}^{n_{q+1}^L}a_i\sum_{j=1}^{n_q^L}\left\langle\partial_{q+1}^L\sigma_i^{q+1},\sigma_j^{q}\right\rangle\sigma_j^{q}\\
        &=\sum_{j=1}^{n_q^L}\sum_{i=1}^{n_{q+1}^L}a_i\left\langle\partial_{q+1}^L\sigma_i^{q+1},\sigma_j^{q}\right\rangle\sigma_j^{q}.
        % &=\sum_{i=1}^{n_{q+1}^L}a_i\sum_{j=0}^{q+1} (-1)^jd_je_i^{q+1}\\
        % &=\dots=\sum_{j=1}^{n_q^L}\left(\sum_{i:e_i^{q+1}\supset e_j^q\text{, pos oriented}}a_i-\sum_{i:e_i^{q+1}\supset e_j^q\text{, neg oriented}}a_i\right)e_j^q.
    \end{align*}
    Therefore, $c\in C^{K,L}_{q+1}$, or equivalently $\partial_{q+1}^L(c)\in C^K_q$, if and only if $\sum_{i=1}^{n_{q+1}^L}a_i\left\langle\partial_{q+1}^L\sigma_i^{q+1},\sigma_j^{q}\right\rangle=0$ for all $n_q^K<j\leq n_q^L$. Finally, whenever $j>n_q^K$, we have $\left\langle\partial_{q+1}^L\sigma_i^{q+1},\sigma_j^{q}\right\rangle=0$ for $i\leq n_{q+1}^K$, since the boundary of a $(q+1)$-simplex of $K$ lies in $K$. We therefore obtain the following equivalences:
    \begin{equation*}
        c_{K,L}\in\ker B_{q,2}^{K,L}\iff \sum_{i=n_{q+1}^K+1}^{n_{q+1}^L}a_i\left\langle\partial_{q+1}^L\sigma_i^{q+1}, \sigma_j^q\right\rangle=0,\,\forall j>n_q^K \iff c\in C^{K,L}_{q+1}.
    \end{equation*}
    % CHECK: \sigma_j should be \sigma_j^q. Also consider ending with one sentence that concludes the lemma: since c_K is arbitrary,
    % C_{q+1}^{K,L} consists exactly of the chains c_K+c_{K,L} with c_{K,L}\in\ker B_{q,2}^{K,L}.
\end{proof}

\subsection{Proof of Lemma \ref{lemma:B22_horizontal}}
\begin{proof}
    In a simplicial filtration, we can write
    \begin{equation*}
        B_{q,2}^{s,t}=\begin{bmatrix}
            [\partial_{q+1}^{t}](I_q^{s,t-1}, I_{q+1}^{s,t-1})&[\partial_{q+1}^t](I_q^{s,t-1},I_{q+1}^{t-1,t})\\
            0&[\partial_{q+1}^t](I_{q}^{t-1,t},I_{q+1}^{t-1,t})
        \end{bmatrix}=\begin{bmatrix}
            B_{q,2}^{s,t-1}&[\partial_{q+1}^t](I_q^{s,t-1},I_{q+1}^{t-1,t})\\
            0&[\partial_{q+1}^t](I_{q}^{t-1,t},I_{q+1}^{t-1,t})
        \end{bmatrix}.
    \end{equation*}
    Let $v\in\ker    \begin{bmatrix}    (B_{q,2}^{s,t-1})^T B_{q,2}^{s,t-1} & 0\\    0 & I_{n_{q+1}^t-n_{q+1}^{t-1}}    \end{bmatrix}$, then we can write $v=\begin{bmatrix} v_1\\ v_2 \end{bmatrix}$ where $v_1\in \ker (B_{q,2}^{s,t-1})^T B_{q,2}^{s,t-1}=\ker B_{q,2}^{s,t-1}$ and $v_2=0\in\mathbb{R}^{n_{q+1}^t-n_{q+1}^{t-1}}$. This allows for the computation
    \begin{align*}
        B_{q,2}^{s,t}v&=\begin{bmatrix}
            B_{q,2}^{s,t-1}&[\partial_{q+1}^t](I_q^{s,t-1},I_{q+1}^{t-1,t})\\
            0&[\partial_{q+1}^t](I_{q}^{t-1,t},I_{q+1}^{t-1,t})
        \end{bmatrix}\begin{bmatrix} v_1\\ v_2 \end{bmatrix}=0.
    \end{align*}
    Therefore,
    \[
        \ker\begin{bmatrix}
            (B_{q,2}^{s,t-1})^T B_{q,2}^{s,t-1} & 0\\
            0 & I_{n_{q+1}^t-n_{q+1}^{t-1}}
        \end{bmatrix}\subseteq\ker B_{q,2}^{s,t}.
    \]
    Finally, Equations~\eqref{eq:pinv_subset_ker} and~\eqref{eq:pinv_def_A_inv} give
    \begin{align*}
        (B_{q,2}^{s,t})^\dagger B_{q,2}^{s,t}&=(B_{q,2}^{s,t})^\dagger B_{q,2}^{s,t}\begin{bmatrix}
            (B_{q,2}^{s,t-1})^TB_{q,2}^{s,t-1} & 0\\
            0 & I_{n_{q+1}^t-n_{q+1}^{t-1}}
        \end{bmatrix}^\dagger \begin{bmatrix}
            (B_{q,2}^{s,t-1})^TB_{q,2}^{s,t-1} & 0\\
            0 & I_{n_{q+1}^t-n_{q+1}^{t-1}}
        \end{bmatrix}\\
        &=(B_{q,2}^{s,t})^\dagger B_{q,2}^{s,t}\begin{bmatrix}
            (B_{q,2}^{s,t-1})^\dagger B_{q,2}^{s,t-1} & 0\\
            0 & I_{n_{q+1}^t-n_{q+1}^{t-1}}
        \end{bmatrix}.
    \end{align*}
    The final equality in~\eqref{eq:B22_horizontal} follows by taking the transpose of both sides, since the matrices involved are symmetric.
\end{proof}

\subsection{Proof of Lemma \ref{lemma:B22_vertical}}
\begin{proof}
    As in the proof of Lemma~\ref{lemma:B22_horizontal}, we write
    \begin{equation*}
        B_{q,2}^{s-1,t}=\begin{bmatrix}
        [\partial_{q+1}^t](I_q^{s-1,s},I_{q+1}^{s-1,s})&[\partial_{q+1}^t](I_q^{s-1,s},I_{q+1}^{s,t})\\
        0&[\partial_{q+1}^t](I_q^{s,t},I_{q+1}^{s,t})
    \end{bmatrix}=\begin{bmatrix}
        [\partial_{q+1}^t](I_q^{s-1,s},I_{q+1}^{s-1,s})&[\partial_{q+1}^t](I_q^{s-1,s},I_{q+1}^{s,t})\\
        0&B_{q,2}^{s,t}
    \end{bmatrix}.
    \end{equation*}
    Hence for $v\in \ker B_{q,2}^{s-1,t}$, we can write $v=\begin{bmatrix} v_1\\ v_2 \end{bmatrix}$, where we know $v_2\in\ker B_{q,2}^{s,t}=\ker \left(B_{q,2}^{s,t}\right)^TB_{q,2}^{s,t}$ and $v_1\in\mathbb{R}^{n_{q+1}^{s}-n_{q+1}^{s-1}}$ is the remaining part of the vector. For any $v_1$, we then have that $v\in\ker \begin{bmatrix}
        0_{n_{q+1}^{s}-n_{q+1}^{s-1}}&0\\
        0&(B_{q,2}^{s,t})^TB_{q,2}^{s,t}
    \end{bmatrix}$.
    Therefore, $\ker B_{q,2}^{s-1,t}\subseteq\ker\begin{bmatrix}
        0&0\\
        0&(B_{q,2}^{s,t})^TB_{q,2}^{s,t}
    \end{bmatrix}$, and Equations~\eqref{eq:pinv_def_A_inv} and~\eqref{eq:pinv_subset_ker} give
    \begin{align*}
    \begin{bmatrix}
            0&0\\
            0&(B_{q,2}^{s,t})^\dagger B_{q,2}^{s,t}
        \end{bmatrix} &=\begin{bmatrix}
        0&0\\
        0&(B_{q,2}^{s,t})^TB_{q,2}^{s,t}
    \end{bmatrix}^\dagger \begin{bmatrix}
        0&0\\
        0&(B_{q,2}^{s,t})^TB_{q,2}^{s,t}
    \end{bmatrix}\\
    &=\begin{bmatrix}
        0&0\\
        0&(B_{q,2}^{s,t})^TB_{q,2}^{s,t}
    \end{bmatrix}^\dagger \begin{bmatrix}
        0&0\\
        0&(B_{q,2}^{s,t})^TB_{q,2}^{s,t}
    \end{bmatrix}(B_{q,2}^{s-1,t})^\dagger B_{q,2}^{s-1,t}\\
    &=\begin{bmatrix}
            0&0\\
            0&(B_{q,2}^{s,t})^\dagger B_{q,2}^{s,t}
        \end{bmatrix}(B_{q,2}^{s-1,t})^\dagger B_{q,2}^{s-1,t}.
    \end{align*}
    Again, the final equality in~\eqref{eq:B22_vertical} follows by taking the transpose, since the matrices involved are symmetric.
\end{proof}

%% file: reference_clean.bib
@inproceedings{chen_persistent_2011,
  author = {Chen, Chao and Kerber, Michael},
  title = {Persistent homology computation with a twist},
  booktitle = {Proceedings of the 27th {European} {Workshop} on {Computational} {Geometry} ({EuroCG})},
  year = {2011},
}

@inproceedings{bauer_clear_2014,
  author = {Bauer, Ulrich and Kerber, Michael and Reininghaus, Jan},
  editor = {Bremer, Peer-Timo and Hotz, Ingrid and Pascucci, Valerio and Peikert, Ronald},
  title = {Clear and {Compress}: {Computing} {Persistent} {Homology} in {Chunks}},
  booktitle = {Topological {Methods} in {Data} {Analysis} and {Visualization} {III}},
  pages = {103--117},
  publisher = {Springer International Publishing},
  address = {Cham},
  year = {2014},
  doi = {10.1007/978-3-319-04099-8_7},
}

@article{pv_tests,
  author = {Vittorietti, Martina and Kok, Piet J. J. and Sietsma, Jilt and Li, Wei and Jongbloed, Geurt},
  title = {General framework for testing Poisson-Voronoi assumption for real microstructures},
  journal = {Applied Stochastic Models in Business and Industry},
  volume = {36},
  number = {4},
  pages = {604--627},
  year = {2020},
  doi = {10.1002/asmb.2517},
}

@article{persistence_landscapes,
  author = {Bubenik, Peter},
  title = {Statistical topological data analysis using persistence landscapes},
  journal = {Journal of Machine Learning Research},
  volume = {16},
  pages = {77--102},
  year = {2015},
}

@article{laplacians_survery,
AUTHOR = {Wei, Xiaoqi and Wei, Guo-Wei},
TITLE = {Persistent Topological Laplacians—A Survey},
JOURNAL = {Mathematics},
VOLUME = {13},
YEAR = {2025},
NUMBER = {2},
ARTICLE-NUMBER = {208},
URL = {https://www.mdpi.com/2227-7390/13/2/208},
ISSN = {2227-7390},
DOI = {10.3390/math13020208}
}

@article{laplacians_algorithms,
  author = {M\'{e}moli, Facundo and Wan, Zhengchao and Wang, Yusu},
  title = {Persistent Laplacians: Properties, Algorithms and Implications},
  journal = {SIAM Journal on Mathematics of Data Science},
  volume = {4},
  number = {2},
  pages = {858--884},
  year = {2022},
  doi = {10.1137/21M1435471},
}

@article{laplacians_spectral_graph,
  author = {Wang, Rui and Nguyen, Duc Duy and Wei, Guo-Wei},
  title = {Persistent spectral graph},
  journal = {International Journal for Numerical Methods in Biomedical Engineering},
  volume = {36},
  number = {9},
  pages = {e3376},
  year = {2020},
  doi = {10.1002/cnm.3376},
}

@article{pinv_properties,
  author = {Hung, Ching-Hsiang and Markham, Thomas L.},
  title = {The {Moore}--{Penrose} inverse of a partitioned matrix $M=\left(\begin{smallmatrix}A&0\\ B&C\end{smallmatrix}\right)$},
  journal = {Czechoslovak Mathematical Journal},
  volume = {25},
  number = {3},
  pages = {354--361},
  year = {1975},
}

@inproceedings{pers_lap_mnist,
  author = {Davies, Thomas and Wan, Zhengchao and Sanchez-Garcia, Ruben J.},
  title = {The {Persistent} {Laplacian} for {Data} {Science}: {Evaluating} {Higher}-{Order} {Persistent} {Spectral} {Representations} of {Data}},
  booktitle = {Proceedings of the 40th {International} {Conference} on {Machine} {Learning}},
  pages = {7249--7263},
  publisher = {PMLR},
  year = {2023},
  url = {https://proceedings.mlr.press/v202/davies23c.html},
}

@article{harmonische_funktionen,
  author = {Eckmann, Beno},
  title = {Harmonische {Funktionen} und {Randwertaufgaben} in einem {Komplex}},
  journal = {Commentarii Mathematici Helvetici},
  volume = {17},
  number = {1},
  pages = {240--255},
  year = {1944},
  doi = {10.1007/BF02566245},
}

@article{computing_ph,
  author = {Zomorodian, Afra and Carlsson, Gunnar},
  title = {Computing {Persistent} {Homology}},
  journal = {Discrete \& Computational Geometry},
  volume = {33},
  number = {2},
  pages = {249--274},
  year = {2005},
  doi = {10.1007/s00454-004-1146-y},
}

@article{pers_lap_covid,
  author = {Jiahui Chen and Yuchi Qiu and Rui Wang and Guo-Wei Wei},
  title = {Persistent Laplacian projected Omicron BA.4 and BA.5 to become new dominating variants},
  journal = {Computers in Biology and Medicine},
  volume = {151},
  pages = {106262},
  year = {2022},
  doi = {10.1016/j.compbiomed.2022.106262},
}

@article{pers_lap_protein_ligand,
  author = {Meng, Zhenyu and Xia, Kelin},
  title = {Persistent spectral–based machine learning ({PerSpect} {ML}) for protein-ligand binding affinity prediction},
  journal = {Science Advances},
  volume = {7},
  number = {19},
  pages = {eabc5329},
  year = {2021},
  doi = {10.1126/sciadv.abc5329},
}

@article{hermes,
title = {HERMES: Persistent spectral graph software},
journal = {Foundations of Data Science},
volume = {3},
number = {1},
pages = {67-97},
year = {2021},
issn = {},
doi = {10.3934/fods.2021006},
url = {https://www.aimsciences.org/article/id/51d06f63-bf7f-44ac-a3c5-c3bc54c03bf9},
author = {Rui Wang and Rundong Zhao and Emily Ribando-Gros and Jiahui Chen and Yiying Tong and Guo-Wei Wei}
}

@article{projection_matrices,
  author = {Piziak, R. and Odell, P.L. and Hahn, R.},
  title = {Constructing projections on sums and intersections},
  journal = {Computers \& Mathematics with Applications},
  volume = {37},
  number = {1},
  pages = {67--74},
  year = {1999},
  doi = {10.1016/S0898-1221(98)00242-9},
}

@article{tda,
  author = {Wasserman, Larry},
  title = {Topological data analysis},
  journal = {Annual Review of Statistics and Its Application},
  volume = {5},
  pages = {501--532},
  year = {2018},
  doi = {10.1146/annurev-statistics-031017-100045},
}

@misc{computationalefficiency,
  author = {Jones, Benjamin and Wei, Guo-Wei},
  title = {{PETLS}: {PE}rsistent {T}opological {L}aplacian {S}oftware},
  howpublished = {arXiv preprint arXiv:2508.11560},
  year = {2025},
}

@article{combinatorial_laplacians,
  author = {Horak, Danijela and Jost, J{\"u}rgen},
  title = {Spectra of combinatorial {Laplace} operators on simplicial complexes},
  journal = {Advances in Mathematics},
  volume = {244},
  pages = {303--336},
  year = {2013},
  doi = {10.1016/j.aim.2013.05.007},
}

@article{carlsson_topology_data,
  author = {Carlsson, Gunnar},
  title = {Topology and data},
  journal = {Bulletin of the American Mathematical Society},
  volume = {46},
  number = {2},
  pages = {255--308},
  year = {2009},
}

@article{ghrist_barcodes,
  author = {Ghrist, Robert},
  title = {Barcodes: the persistent topology of data},
  journal = {Bulletin of the American Mathematical Society},
  volume = {45},
  number = {1},
  pages = {61--75},
  year = {2008},
}

@article{chazal_tda,
  author = {Chazal, Fr{\'e}d{\'e}ric and Michel, Bertrand},
  title = {An introduction to topological data analysis: fundamental and practical aspects for data scientists},
  journal = {Frontiers in Artificial Intelligence},
  volume = {4},
  pages = {667963},
  year = {2021},
}

@article{edelsbrunner_persistent_homology,
  author = {Edelsbrunner, Herbert and Letscher, David and Zomorodian, Afra},
  title = {Topological persistence and simplification},
  journal = {Discrete \& Computational Geometry},
  volume = {28},
  number = {4},
  pages = {511--533},
  year = {2002},
  doi = {10.1007/s00454-002-2885-2},
}

@inproceedings{cohen_steiner_stability,
  author = {Cohen-Steiner, David and Edelsbrunner, Herbert and Harer, John},
  title = {Stability of persistence diagrams},
  booktitle = {Proceedings of the 21st Annual Symposium on Computational Geometry},
  pages = {263--271},
  year = {2005},
}

@article{persistence_images,
  author = {Adams, Henry and Emerson, Tegan and Kirby, Michael and Neville, Rachel and Peterson, Chris and Shipman, Patrick and Chepushtanova, Sofya and Hanson, Eric and Motta, Francis and Ziegelmeier, Lori},
  title = {Persistence images: A stable vector representation of persistent homology},
  journal = {Journal of Machine Learning Research},
  volume = {18},
  number = {8},
  pages = {1--35},
  year = {2017},
}

@article{persistent_homology_transform,
  author = {Turner, Katharine and Mukherjee, Sayan and Boyer, Doug M},
  title = {Persistent homology transform for modeling shapes and surfaces},
  journal = {Information and Inference: A Journal of the IMA},
  volume = {3},
  number = {4},
  pages = {310--344},
  year = {2014},
}

@misc{persistent_laplacian_diagrams,
  author = {Jung, Inkee and Kang, Wonwoo and Park, Heehyun},
  title = {Persistent {L}aplacian diagrams},
  howpublished = {arXiv preprint arXiv:2512.05463},
  year = {2025},
}

@misc{effective_resistance,
  author = {Black, Mitchell and Maxwell, William and Nayyeri, Amir},
  title = {An incremental span-program-based algorithm and the fine print of quantum topological data analysis},
  howpublished = {arXiv preprint arXiv:2307.07073},
  year = {2023},
}

@misc{botnan_cheeger_2026,
  author = {Botnan, Magnus Bakke and Dong, Rui},
  title = {Cheeger inequalities for the persistent {L}aplacian},
  howpublished = {arXiv preprint arXiv:2606.02846},
  year = {2026},
}

@article{jost_cheeger_2024,
   title={Cheeger inequalities on simplicial complexes},
   ISSN={0391-173X},
   DOI={10.2422/2036-2145.202307_009},
   journal={Annali della Scuola Normale Superiore di Pisa, Classe di Scienze},
   publisher={Scuola Normale Superiore - Edizioni della Normale},
   author={Jost, Jürgen and Zhang, Dong},
   year={2024},
   month=Dec, pages={30} }

@article{nicolau2011topology,
  author = {Nicolau, Monica and Levine, Arnold J and Carlsson, Gunnar},
  title = {Topology based data analysis identifies a subgroup of breast cancers with a unique mutational profile and excellent survival},
  journal = {Proceedings of the National Academy of Sciences},
  volume = {108},
  number = {17},
  pages = {7265--7270},
  year = {2011},
}

@article{cang2017topologynet,
  author = {Cang, Zixuan and Wei, Guo-Wei},
  title = {TopologyNet: Topology based deep convolutional and multi-task neural networks for biomolecular property predictions},
  journal = {PLoS Computational Biology},
  volume = {13},
  number = {7},
  pages = {e1005690},
  year = {2017},
}

@article{hiraoka2016hierarchical,
  author = {Hiraoka, Yasuaki and Nakamura, Takenobu and Hirata, Akihiko and Escolar, Emerson G and Matsue, Kaname and Nishiura, Yasumasa},
  title = {Hierarchical structures of amorphous solids characterized by persistent homology},
  journal = {Proceedings of the National Academy of Sciences},
  volume = {113},
  number = {26},
  pages = {7035--7040},
  year = {2016},
}

@article{carlsson2008local,
  author = {Carlsson, Gunnar and Ishkhanov, Tigran and De Silva, Vin and Zomorodian, Afra},
  title = {On the local behavior of spaces of natural images},
  journal = {International Journal of Computer Vision},
  volume = {76},
  number = {1},
  pages = {1--12},
  year = {2008},
}

@article{bendich2016persistent,
  author = {Bendich, Paul and Marron, J. S. and Miller, Ezra and Pieloch, Alex and Skwerer, Sean},
  title = {Persistent homology analysis of brain artery trees},
  journal = {The Annals of Applied Statistics},
  volume = {10},
  number = {1},
  pages = {198--218},
  year = {2016},
  doi = {10.1214/15-AOAS886},
}

@misc{dong_faster_2024,
  title = {A faster algorithm of up persistent {Laplacian} over non-branching simplicial complexes},
  url = {http://arxiv.org/abs/2408.16741},
  doi = {10.48550/arXiv.2408.16741},
  publisher = {arXiv},
  author = {Dong, Rui},
  month = aug,
  year = {2024},
  note = {arXiv:2408.16741 [math.AT] version: 1},
}

@article{botnan_nonbranching_2026,
title = {The persistent Laplacian of non-branching complexes},
journal = {Foundations of Data Science},
volume = {10},
number = {0},
pages = {53-89},
year = {2026},
issn = {},
doi = {10.3934/fods.2025020},
url = {https://www.aimsciences.org/article/id/68ea17bacb5dde21e753d644},
author = {Magnus Bakke Botnan and Rui Dong}
}
